\documentclass{amsart}

\usepackage[parfill]{parskip}
\usepackage{graphicx} 
\usepackage{color} 
\usepackage{hyperref} 
\usepackage{graphicx}  
\usepackage{verbatim}  
\hypersetup{pdfpagemode=FullScreen,colorlinks=true}
\usepackage{appendix, url}
\usepackage{amsmath, amssymb, mathrsfs, amsfonts, amsthm} 
\usepackage{tikz-cd}
\DeclareMathOperator{\lie}{lie}
\DeclareMathOperator{\qobs}{q-obs}

\DeclareMathOperator{\area}{area}
\DeclareMathOperator{\carea}{carea}

\DeclareMathOperator{\Aut}{Aut}

\newtheorem {theorem} {Theorem} [section]

\newtheorem{lemma}[theorem] {Lemma}

\newtheorem {definition} [theorem] {Definition} 
\newtheorem {proposition}  [theorem]{Proposition} 

\newtheorem {corollary}[theorem]  {Corollary} 
\newtheorem {example} [theorem]  {Example}
  
\newtheorem {notation}[theorem] {Notation}
\newtheorem {terminology}[theorem]{Terminology}
\newtheorem {remark} [theorem] {Remark}
\numberwithin {equation} {section}

\DeclareMathOperator{\colim}{colim} 
\begin{document} 
\href{http://yashamon.github.io/web2/papers/fukayaII.pdf}{Direct link to author's version}
\author {Yasha Savelyev} 
\today
\address{University of Colima, CUICBAS, Bernal Díaz del
Castillo 340,
Col. Villas San Sebastian,
28045, Colima, Colima,
Mexico}
\email{yasha.savelyev@gmail.com} 
\title [Global Fukaya Category II]{Global Fukaya
category II:  singular connections, quantum obstruction theory and other applications}  
\begin{abstract} 
In part I, using the theory of $\infty$-categories, we
constructed a natural ``continuous action'' of $\operatorname {Ham} (M,  \omega) $ on the Fukaya category of a closed
monotone symplectic manifold. Here we show that this action
is generally homotopically non-trivial, i.e  implicitly the
main part of a conjecture of Teleman. We use this to give
various applications. For example, we find new curvature
constraint phenomena for smooth and singular
$\mathcal{G}$-connections on principal $\mathcal{G}$-bundles over $S
^{4}$, where  $\mathcal{G}$ is  $\operatorname {PU} (2)$ or $\operatorname {Ham} (S ^{2} )$.  Even for the classical group $\operatorname {PU} (2)$, these phenomena are invisible to Chern-Weil theory,
and are inaccessible to known Yang-Mills theory and quantum
characteristic classes techniques. So this can be understood as one application of Floer theory and the theory of $\infty$-categories in basic differential
geometry. We also develop, based on this $\infty$-categorical
Fukaya theory, some basic new integer valued invariants  of smooth manifolds, called quantum obstruction.  On the way we also construct what we call quantum Maslov classes, which are higher degree variants of the relative Seidel morphism. This also leads to new applications in Hofer geometry of
the space of Lagrangian  equators in $S^2$. 
\end{abstract}
\maketitle
\section{Introduction}
Let $Ham(M,\omega)$ denote the group of
Hamiltonian symplectomorphisms of a symplectic
manifold $(M,\omega) $, understood as a Frechet
Lie group, with its $C ^{\infty}$ topology.  A \emph{Hamiltonian bundle} is a smooth fiber
bundle $$M \hookrightarrow P \xrightarrow{p} X,$$ with
structure group $ \text {Ham}(M, \omega) $. Given such a
bundle $P$ with $M$ monotone, in Part I
\cite{citeSavelyevGlobalFukayaCategoryI}
we have constructed a continuous ``classifying map''
\begin{equation*}
   f _{P}: X \to |\mathbb{S}|,
\end{equation*}
where $|\mathbb{S}|$ denotes the space of
$\infty$-categories. Moreover, $f _{P}$ maps into 
the component of $ \operatorname {NFuk} (M,
\omega))$, with the latter denoting the $A _{\infty}
$ nerve of $\operatorname {Fuk} (M,  \omega).  
$ We also denote this component by
$(|\mathbb{S}|,  \operatorname {NFuk} (M,  \omega)) $.  From here on we just refer to \cite{citeSavelyevGlobalFukayaCategoryI} as Part I.

This extends to the universal level, so that there is
a universal (continuous) classifying map: $$\operatorname {BHam} (M,
\omega)  \to (|\mathbb{S}|, \operatorname {NFuk} (M,  \omega) ).$$ As explained in Part I, this is
interpreted as a ``continuous'' (homotopy coherent) action of $\operatorname {Ham}
(M,  \omega) $ on $\operatorname {Fuk} (M,  \omega) $.
Existence of such an action in a somewhat weaker form (on
the level $E _{2}$ algebras rather then space level) has
been conjectured by Teleman ICM 2014.

The construction also induces a certain kind of
simplicial fibration over the smooth singular set $X
_{\bullet }$ of $X$, called categorical fibration:
\begin{equation*}
\operatorname {NFuk} (M,  \omega)
\hookrightarrow {\operatorname {Fuk} _{\infty} (P)} \xrightarrow{p _{\bullet }} X _{\bullet
}.
\end{equation*}
This is called the global Fukaya category of $P$.   We show that for $P$ a non-trivial Hamiltonian $S ^{2} $ fibration over $S ^{4} $, the maximal Kan sub-fibration of $\operatorname {Fuk} _{\infty} (P) $ is non-trivial.
In particular, $\operatorname {Fuk} _{\infty} (P) $ is
non-trivial as a categorical fibration and so $f _{P} $ is homotopically non-trivial. 
In particular, this gives:
\begin{theorem} \label{thm:main}
   The natural homomorphism as constructed in Part I,
	 $$\mathbb{Z} =  \pi _{4} (\operatorname {BHam} (S ^{2},
	 \omega)) \xrightarrow{k} \pi _{4} (\mathbb{S},
	 \operatorname {NFuk} (S ^{2}, \omega )), 
$$
is injective.
\end{theorem}
Thus, we conclude that the natural ``continuous action'' of
$\operatorname {Ham} (M,  \omega) $ on $\operatorname {Fuk}
(M,  \omega) $ is generally homotopically non-trivial.  This of course is implicitly part of Teleman's conjecture mentioned above.

None of the homotopy groups of $\mathbb{S}$ are
known, so the above theorem is also in a sense an application of
geometry to algebraic topology.  Such an application is possible because geometry forces a priori $A _{\infty} $ associativity of certain structures, which then produces needed generators for the relevant homotopy groups. 

Note that the above requires certain chain level
calculations in Fukaya $A _{\infty }$ categories. To this
end, we relate such computations to the computations of certain quantum Maslov classes. The latter are certain higher dimensional analogues of the relative Seidel element in
\cite{citeHuLalondeArelativeSeidelmorphismandtheAlbersmap}.
The calculation of these quantum Maslov classes uses
a regularization technique based on ``virtual Morse theory''
for the Hofer length functional
\cite{citeSavelyevVirtualMorsetheory}. However, given the
work of Chow ~\cite{citeJimmyQuantClasses},  linking quantum
characteristic classes to pseudo-holomorphic quilts, more
algebraic-geometric computations should also be possible in the future.

The arguments of the paper are quiet general,
so that a Hamiltonian $S ^{2}$ fibration over $S ^{4}$ can
be replaced by more general Hamiltonian fibrations with
monotone fiber, obtaining results similar to \ref{thm:main}.
However, as this is a first computation of the kind, we
focus on the ideas and a concrete example.

\begin{remark} \label{remark:}
It is likely that $k$ is surjective. Surjectivity is
in a sense the statement that up to equivalence there are no
exotic categorical fibrations over $S ^{4} $, with fiber
equivalent to $\operatorname {NFuk} (S ^{2},  \omega) $ - they all come from Hamiltonian
$ S ^{2} $ fibrations, via the global Fukaya category.
\end{remark}

\subsection{An application in basic Riemannian geometry}
As one less expected application, we can use the computation
of Theorem \ref{thm:main} to obtain lower bounds for
the curvature of certain types of singular connections.  
\begin{definition} \label{def:basicSingular}
Let ${G} \hookrightarrow P \to X$ be a principal $G$ bundle, where $G$ is a Frechet Lie group. A \textbf{\emph{singular $G$-connection}} on $P$ is a  closed subset $C \subset X$, and a smooth Ehresmann $G$-connection 
 $\mathcal{A}$ on $P| _{X-C} $.
\end{definition}
The above definition is basic, as one often puts additional conditions, see for instance \cite{citeReeseHarverSingular}, \cite{citeSingularSibnerSibner}.

\subsubsection {A non-metric measure of curvature} \label{sec:non-metric}

Let $G$ as above be a Frechet Lie group,
we denote by $\lie G$ its Lie algebra and let
$$\mathfrak n: \lie G   \to \mathbb{R}$$ be an $Ad$ invariant Finsler norm.  For a principal $G$-bundle $P$ over a Riemann surface $(S,j)$, and given a $G$ connection $\mathcal{A}$ on $P$ define a 2-form $\alpha _{\mathcal{A}} $ on $S$ by: 
\begin{equation*}
   \alpha _{\mathcal{A}}  (v,jv) = \mathfrak n (R _{\mathcal{A}} (v, jv)),
\end{equation*}
where $R _{\mathcal{A}} $ is the curvature 2-form of
$\mathcal{A}$. More specifically, the latter form has the
properties: $$R _{\mathcal{A}} (v, w) \in \lie \Aut P_{z},
$$ for $z \in S$, $v,w \in T _{z} S  $, $P _{z}$ the fiber
of $P$ over $z$, $\Aut P _{z} \simeq G $ the group of $G$-torsor
automorphisms of $P _{z} $, and where $\simeq$ means non-canonical group isomorphism.

Define
\begin{equation} \label{eq:areasurface}
\area _{\mathfrak n}  (\mathcal{A}) := \int _{S} \alpha _{\mathcal{A}}.
\end{equation}
In the case $\mathcal{A}$ is singular with singular set $C$, $\alpha _{\mathcal{A}} $  is defined on $S -C$ so we 
define 
\begin{equation*}
\area _{\mathfrak n} (\mathcal{A}) := \int _{S - C} \alpha _{\mathcal{A}},
\end{equation*}
with the right-hand side now being an extended integral. This $\area$ is a non-metric measurement meaning that no Riemannian metric on $S$ is needed.

It is possible to extend the functional above to
a functional on the space $\mathcal{C}$ of $G$-connections
on principal $G$ bundles $P \to \Delta ^{n} $. It may seem
that $\Delta ^{n} $ has no connection to Riemann surfaces,
but in fact there is an intriguing such connection.  
Let ${\mathcal {S}} _{d} $ denote the universal curve over
$\overline{\mathcal{R}} _{d} $ - the moduli space of complex
structures on the disk with $d+1$ punctures on the boundary.
And let ${\mathcal {S}} ^{\circ} _{d} $ denote ${\mathcal {S}} _{d} $, with nodal points of the fibers removed. 
Then it is shown in Part I that there are certain axiomatized systems of maps:
\begin{equation*}
{u}: {\mathcal {S}} ^{\circ}  _{d}  \to \Delta
^{n}, \text{$d,n$ varying}.  
\end{equation*}
Such a system is uniquely determined up to suitable homotopy, and is referred to as $\mathcal{U}$. 

There is then a natural functional:
\begin{align} \label{eq:areaG}
   \area _{\mathcal{U}}: \mathcal{C} \to \mathbb{R} _{\geq 0},
\end{align}
defined with respect to a choice of $\mathcal{U}$, see
Section \ref{sec:qcurvature}.  When $n=2$ it is just the $\area$ functional as previously defined.
\subsubsection {Abstract resolutions of singular
connections} Avoiding generality, suppose that
$\mathcal{A}$ is a singular $G$-connection on a
principal $G$-bundle $P \to S ^{n}$, with a single
singularity $x _{0} $. We will show that it is
possible to control the curvature of the singular
connection $\mathcal{A}$ if we impose a certain
structure on the singularity of $\mathcal{A}$.
The simplest way to do this is to ask for
existence of a certain kind of abstract
resolution.  

First, a simplicial $G$-connection 
on $\mathcal{D}$ on $P$, as defined in Section
\ref{section:simplicialconnections}, is basically a functorial assignment of a smooth $G$ connection $\mathcal{D}  _{\Sigma}  $ on $\Sigma ^{*} P$ for each smooth $$\Sigma: \Delta ^{n} \to S ^{n} . $$
\begin{definition} \label{def:resolution}
For $\mathcal{A}, P$ as above a
\textbf{\emph{simplicial resolution}} of
$\mathcal{A} $ is a simplicial $G$ connection
$\mathcal{A} ^{res} $ on $P$, with
the following property. Let $\Sigma _{0}: \Delta ^{n} \to
S ^{n}  $ represent the generator of $\pi _{n} (S ^{n},
x _{0}) $ (cf. Appendix \ref{appendix:Kan}), then  $$\Sigma ^{*} _{0}   | _{interior \, \Delta ^{n}}  \mathcal{A} = \mathcal{A} ^{res} _{\Sigma _{0} }| _{interior \, \Delta ^{n}}. $$  
\end{definition}

In the following theorem $G=PU (2), n=4$ and the
norm $\mathfrak n$ on $\lie PU (2)$ will be taken to be
the operator norm, normalized so that the Finsler
length of the shortest one parameter subgroup from
$id$ to $-id$ is $\frac{1}{2}$. We will omit $\mathfrak n$  in notation. We also impose an additional constraint on $\mathcal{A} ^{res} $, 
so that the curvature at ``$\infty$'' is bounded by a threshold, which means the following. Let $$\Sigma _{\infty}: \Delta ^{4} \to x _{0} $$ be the constant map. Suppose that:
\begin{equation*}
   \area _{\mathcal{U}} \mathcal{A} ^{res} _{\Sigma _{\infty}} < 1/2,
\end{equation*}
and suppose for simplicity that $\mathcal{A} ^{res} _{\Sigma
_{\infty}}$ is trivial along the edges of $\Delta ^{4} $,
later on this condition is relaxed, see Proposition
\ref{prop:alternative}.  (This condition can be completely
removed but at the cost of significant additional
complexity.) 

%
We say in this case that $\mathcal{A} ^{res} $ is a \textbf{\emph{sub-quantum resolution}}. The following is proved in Section \ref{sec:qcurvature}.

\begin{theorem}  \label{thm:lowerboundsingular}
Let $P \to S ^{4} $ be a non-trivial principal $PU (2)$ bundle.
Let $\mathcal{A}$ be a singular $PU (2)$-connection on $P$ with a single singularity at $x _{0} $. Then for any sub-quantum resolution $\mathcal{A} ^{res}  $ 
of $\mathcal{A}$  and for any $\mathcal{U}$ as above 
$$\area _{\mathcal{U}}   (\mathcal{A} ^{res}  _{\Sigma _{0}} ) \geq 1/2. $$     
%
\end{theorem}
The theorem has certain extensions to Hamiltonian
singular connections $\mathcal{A}$, understanding
$P$ as a principal $Ham (S ^{2} )$ bundle, 
Section \ref{sec:qcurvature}. 

Here is one basic class
of examples. 
\begin{example}
   \label{example:singularconnection} Let
   $P$ be as above, and $\mathcal{A}' $   be an
   ordinary
   smooth $PU (2) $ connection on  $P$. Express
   $S ^{4} $ as a union of sub-balls $D ^{4}
   _{\pm} \subset S ^{4}$,
    intersecting only in the boundary. Suppose
    that that we have the 
   property that $\area  _{\mathcal{U} }(\mathcal{A} ') \,|\, _{D ^{4}
   _{-}} < \frac{1}{2}.$  Let $\mathcal{A} $ be
   the singular connection on $P$ obtained as  the
   push-forward of $\mathcal{A} '$ by the bundle
   map $\widetilde{q}: P \to P$ over the singular     smooth map
   $q: S ^{4} \to S ^{4}$ taking $D ^{4} _{-}$ to a
   single point $\infty \in S ^{4}$, with $q|
   _{interior \, D ^{4} _{-}}$ an
   immersion.  
   Then
   $\mathcal{A} $ has a sub-quantum
   resolution $\mathcal{A} ^{res} $ essentially by
   construction, and we will not elaborate.
   In this case, the theorem above
    simply yields that  $\area _{\mathcal{U}}   (\mathcal{A}' _{D ^{4}_+}   ) \geq 1/2. $ 
 \end{example} 
Let us summarize the above example as the
following basic differential geometric result.
It can be formally understood as a
corollary of Theorem \ref{thm:lowerboundsingular},
as partially explained above, but it is
more elementary to see it as a corollary of  
Theorem \ref{prop:alternative}, which appears
later.    
\begin{corollary} [Of Theorem \ref{thm:lowerboundsingular}
  and of Theorem \ref{prop:alternative}] 
   \label{corol:example}  Let $P$ be a non-trivial
   $\operatorname {PU} (2) $ (or $\operatorname {Ham} (S ^{2}) $)  bundle $P
   \to S ^{4}$, let $D ^{4}
   _{\pm}$ be as above and 
let $\mathcal{A} $ be a smooth $PU (2) $ or $\operatorname
{Ham} (S ^{2},  \omega)  $ connection 
on $P$.  Suppose that $\area  _{\mathcal{U}
   }(\mathcal{A} ) \,|\, _{D ^{4} 
   _{-}} < \frac{1}{2},$  then $\area  _{\mathcal{U} }(\mathcal{A} ) \,|\, _{D ^{4}
   _{+}} \geq \frac{1}{2}.$
\end{corollary}  
The proof of even the above corollary 
traverses the entirety of the theory here. As this
is a very elementary result we may hope for a
simpler argument. It is quite
non-obvious how to do this even for
$PU (2) $. In particular the computation of the
quantum Maslov class, in Section
\ref{section:computeSeidel}, alone is
insufficient, as we also have to do some kind of
algebraic topological gluing of the Floer theory data. In the case of $PU (2) $, one idea
might be to replace Floer theory, used here, by the technically
simpler mathematical Yang-Mills theory over surfaces
~\cite{citeAtiyahBottTheYang-MillsequationsoverRiemannsurfaces}. If we want to mimic the argument presented in this paper, then we should first
extend Yang-Mills theory to work with $G$-bundles
over surfaces with corners and holonomy
constraints over boundary. This might be possible,
but beyond this things are unclear, since, as
mentioned, we also
use certain abstract algebraic topology to glue
the data, and it is not clear how this would work
for Yang-Mills theory.  

We may use the same idea as in the example above to ``push
forward'' simplicial, (not just smooth) 
connections to singular connections with 
more complicated singularities, in such a way
that we again
by construction would have         
sub-quantum resolutions. In this case  Theorem
\ref{thm:lowerboundsingular} no longer has an elementary
interpretation as in corollary above.

There are possible physical interpretations for singular
connections, as appearing in the context here.   A $PU (2) =PSU
(2) $
connection $\mathcal{A} $ on $P$ in physical terms represents
a Yang-Mills field on the space-time $S ^{4}$.
When the space-time has a black hole singularity,
the fields solving the Einstein-Yang-Mills
equations (mathematically connections as above) likewise
develop singularities.  There is a wealth of
physics literature on this subject, and I don't know what
has the highest priority, but here 
is one reference ~\cite{citeYangMillsBlackHoles}. 
 As quantum gravity is often related to simplicial
ideas, it is not inconceivable that the
mathematical sub-quantum
resolution condition above also has a (quantum gravity
theoretic) physical interpretation.

At this point the reader may be curious why
Theorem \ref{thm:main} has something to do with
Theorem \ref{thm:lowerboundsingular}. We cannot give the full story,
but the idea is that the categorical fibration $Fuk
_{\infty} (P) $ only sees the principal bundle $P$
(and its curvature) by the behavior of certain
holomorphic curves. When one has the sub-quantum
condition on the curvature of $\mathcal{A} ^{res}
_{\Sigma _{\infty}}$, certain holomorphic curves
are ruled out so that from the view point of $Fuk
_{\infty} (P) $, $\mathcal{A} ^{res} _{\Sigma
_{\infty}} $ is the trivial connection, (its
curvature is undetectable)  but $Fuk _{\infty} (P) $
is non-trivial as a fibration so that the aforementioned holomorphic curves and consequently curvature must
appear elsewhere.

\subsection {First quantum obstruction and smooth
invariants} \label{sec:firstquantumobstruction}

It is very tempting to use the theory of the
global Fukaya category to find new invariants of
smooth manifolds. One such invariant is already
discussed in Part I, as the homotopy class of the
classifying map $X \to \mathbb{S}$ of the
projectivized, complexified tangent bundle of a
smooth manifold $X$. This by itself is not a very
practical invariant, but we may try to extract
more manageable invariants from this. We present
here a construction of an integer valued invariant which
is based on our theory. This is probably just the
beginning of the story for invariants of smooth
manifolds based on Floer-Fukaya theory. 

Let $M \hookrightarrow P \xrightarrow{p} X$ be a Hamiltonian $M$-bundle, as
previously. Let $$\operatorname {Fuk} _{\infty} (P)
\xrightarrow{p _{\bullet }} X _{\bullet }$$
be the associated categorical fibration, and let  $$K
(P) \to X _{\bullet}$$ be its maximal Kan sub-fibration as in
Lemma \ref{lemma:kanfib}. Then $|K (P)| \to X$ is
a Serre fibration, where $|K (P) |$ is the geometric
realization.

Define 
\begin{equation*} \qobs (P) \in \mathbb{N}
\sqcup \{\infty\}, 
\end{equation*} to be the degree of the first obstruction to
a section of $|K (P)|$. That is $\qobs (P)$ is the smallest integer $n$ such that there is no section of $|K (P)|$ over
the $n$ skeleton of $X$, with respect to some chosen 
CW structure. This is independent of the choice
of the CW structure, as any pair of CW structures on $X$ are
filtered (using cellular filtration)  homotopy equivalent up
to a wedge sum with some collection of $D ^{n}$, $n \in \mathbb{N} $  (with its canonical CW structure),
see Faria~\cite[Theorem 2.4]{citeFilteredHomotopyTypeFaria}.

When no such $n$ exists we set $\qobs (P) = \infty$.
\begin{theorem}
   \label{thm:firstquantumobstruction} Let $S ^{2}
   \hookrightarrow P \to S ^{4} $ be a non-trivial
   Hamiltonian fibration then:
\begin{equation*}
   \qobs (P)= 4. 
\end{equation*} 
\end{theorem}
Indeed the proof of Theorem \ref{thm:main} can be
understood as showing that the associated
obstruction class in $$H ^{4} (S ^{4}, \pi _{3}
(NFuk (S ^{2}))) $$ is non-trivial.
\subsubsection {First quantum obstruction as a
manifold invariant} \label{sec:manifoldinvariant}
Let $X$ be a smooth manifold, and let $P (X)$
denote the fiber-wise projectivization of $TX
\otimes \mathbb{C}$. 
We then define 
\begin{equation*}
   \qobs (X):= \qobs (P (X)) \in \mathbb{N} \sqcup \{\infty\},
\end{equation*}
which is then an invariant of the smooth manifold
$X$. Either this invariant is expressible in terms
of classical invariants, which would be
fascinating since the construction is in terms
pseudo-holomorphic curves or this invariant is
new, that is not expressible in classical terms,
which would also be interesting. There are of
course gauge theory based invariants of smooth
(3,4)-folds, like Donaldson and Seiberg-Witten
invariants. I do not see any connections of the
above to these invariants at the moment, even in
dimension 4. It should be noted that this ``first
quantum obstruction'' invariant is only sensitive
to the tangent bundle, whereas for example
Donaldson invariants can see finer aspects of the
smooth structure. In fact the ``quantum Novikov
conjecture'' of Part I would immediately imply
that the first quantum obstruction is only a
topological invariant of $X$.

\subsection {Hamiltonian rigidity vs flexibility} \label{section:introHamrigidity} By way of the calculation we also obtain an application in Hofer geometry. It can be understood as a relative analogue of a result in \cite{citeSavelyevBottperiodicityandstablequantumclasses}. 

Let $Lag (M, L _{0} )$ denote the space of oriented Lagrangian  submanifolds of a symplectic manifold $(M,\omega)$, Hamiltonian isotopic to $L _{0} $, we may also just write $Lag (M)$. Let $\Omega  _{L _{0}} Lag (M) $ denote the space of based smooth loops in $Lag (M)$, constant near end points, and let $\Omega ^{taut}   _{L _{0} } Lag (M) \subset \Omega  _{L _{0} } Lag (M) $ be the subspace of loops taut concordant to the constant loop at $L _{0} $. The notion of taut concordance is defined in more generality in Definition \ref{definition:concordance}. In the case above, two loops $$p _{1}, p _{2} \in \Omega  _{L _{0} } Lag (M)  $$ are said to be \emph{taut concordant} if the following holds: 
\begin{itemize}
	\item There is a Lagrangian sub-fibration $$\mathcal{L}
	\subset Cyl \times M, \quad Cyl=S ^{1} \times [0,1], $$
	such that $\mathcal{L}$ over the boundary circles
	corresponds, in the natural sense, to the pair $p _{1},
	p _{2} $.
	\item There is a Hamiltonian connection $\mathcal{A}$ on
	$M \times [0,1]$ preserving $\mathcal{L}$, such that
	the coupling form $\Omega _{\mathcal{A}} $ of $\mathcal{A}$ vanishes on $\mathcal{L}$. See Section \ref{sec:couplingform} for the definition of coupling forms.
\end{itemize}
Note that of course $Lag (S ^{2} )$ is homotopy equivalent
to $Lag ^{eq}  (S ^{2} ) \simeq S ^{2} $ where $Lag ^{eq} (S
^{2} ) $ denotes the space of oriented equators in $S ^{2}
$.  Moreover, there is an embedding $$i: (\Omega Lag ^{eq} (S ^{2}) \simeq \Omega S ^{2}) \hookrightarrow  \Omega ^{taut} _{L _{0}} Lag ( S ^{2}),$$ 
as two loops $p _{1}, p _{2} \in \Omega Lag ^{eq} (S ^{2} ) \simeq \Omega S ^{2} $ are taut concordant iff they are homotopic in $Lag ^{eq} (S ^{2} )$, see Lemma \ref{lemma:tautLagS2}. 
\begin{theorem} \label{thm:Hofer}
   Let $L _{0} \subset S ^{2} $ be the equator. And let $$f: S ^{2} \to
   \Omega _{L _{0} } ^{taut}  Lag (S ^{2}),$$  represent $i_*g$, for $g$ the generator of $$\pi _{2}(\Omega (S ^{2})) \simeq \pi _{3} (S ^{2})) \simeq \mathbb{Z},$$ and $i$ as above. Then we have identity for the systole with respect to $L ^{+} $:
   \begin{equation*}
  \min_{f', [f'] = [f]} \max _{s \in S ^{2} } L ^{+} (f' (s)) = 1/2 \cdot \area (S ^{2},
  \omega),
  \end{equation*} where $L ^{+} $ denotes the positive Hofer length functional, as defined in Section \ref{section:hoferlength}. The minimum is attained on a cycle of equators in $S ^{2} $.
\end{theorem}
Even though everything is now smooth, this is not obvious.
For suppose by contrast we measure a related quantity of the ``girth'' (infimum of the diameter of a representative) of the generator $[g]$ of $\pi _{2} Lag (S ^{2} ),$ 
as in \cite{citeItamar}. Then there is an upper bound for
this girth, which is smaller then the lower bound for girth considered in the subspace of $Lag (S ^{2} )$ consisting of equators. 
In other words, if we generalize from equators to general
oriented $S ^{1} $ Lagrangians in $S ^{2} $ we may reduce
the girth to less than the classically expected quantity. By ``classical'' we mean for the classical objects: great circles.
Indeed, it may be that girth of the generator $$[g] \in \pi _{2} Lag (S ^{2} )$$ 
is actually 0. (This would rather astonishing however.)
On the other hand, our theorem says that this kind of
squeezing cannot happen at all for the systole we
consider. In other words whereas our systole exhibits Hamiltonian rigidity, the girth in \cite{citeItamar} while closely related, exhibits flexibility. 

Theorem \ref{thm:Hofer} is proved in Section
\ref{section:applicationHofer}. On the way in
Section \ref{section:SeidelMorph} we construct the
quantum Maslov classes. We show their  non-triviality in Section \ref{section:computeSeidel}.
The Sections \ref{section:applicationHofer},
\ref{section:SeidelMorph}, \ref{section:computeSeidel} are mostly logically independent of the $\infty$-categorical and even the $A _{\infty} $ setup and may be read independently.
Theorems \ref{thm:main}, \ref{thm:firstquantumobstruction} are proved in Section \ref{section:immediateconsequences}, they are basic consequences of the main technical lemma.
\section {Acknowledgements} I am grateful to RIMS institute
at Kyoto university and Kaoru Ono for the invitation,
financial assistance and a number of discussions which took
place there. Much thanks also ICMAT Madrid and Fran Presas
for providing financial assistance, and a lovely research
environment during my stay there. I have also benefited from 
conversations with (in no particular order) Hiro Lee Tanaka,
Mohammed Abouzaid, Kevin Costello, Bertrand Toen and Paul Seidel.
\tableofcontents
\section{Outline} \label{section.application} 
In what follows, when we say Part I we shall mean
\cite{citeSavelyevGlobalFukayaCategoryI}.
We will mostly follow the notation and setup of Part
I.  The reader may review the basics of simplicial
sets, as used by us, in Section 3 of Part I. For a more detailed introduction, which also includes some theory of
quasi-categories, we recommend Riehl
~\cite{citeRiehlAmodelstructureforquasi-categories}.    Here are some specific summary points.
\begin{notation}  
We use notation $\Delta ^{n} $ to denote the standard topological $n$-simplex.
For the standard representable $n$-simplex as a simplicial set we use the notation $\Delta ^{n}
_{\bullet}  $. When $X$ is a smooth manifold $X _{\bullet}
$ will denote the smooth singular set of $X$. That is $X _{\bullet } ([n])   $ is the set of
smooth maps  $\Delta^{n}  \to X$.
If $p: X \to Y$ is a map of spaces, $p _{\bullet}: X _{\bullet} \to Y _{\bullet}  $ will mean the induced simplicial map.  $X
_{\bullet}$ can also denote an abstract
simplicial set when there is no possibility of
confusion. We will denote abstract Kan complexes or
quasi-categories by calligraphic letters e.g.
$\mathcal{X}, \mathcal{Y}$.
\end{notation}
Let us briefly review what we do in Part I. Let $M \hookrightarrow P \xrightarrow{p} X$ be a Hamiltonian fibration. Denote by $\Delta (X):=\Delta/X _{\bullet}$ the smooth simplex category of $X$, with objects smooth maps $\Sigma: \Delta ^{n} \to X $ and morphisms commutative diagrams:
\begin{equation*}
\begin{tikzcd}
   \Delta ^{n} \ar [rd, "\Sigma _{0}"] \ar[r, "mor"] & \Delta ^{m} \ar[d, "\Sigma _{1}" ] \\
& X,
 \end{tikzcd}
\end{equation*}
where $mor: \Delta ^{n} \to \Delta ^{m}  $  is a simplicial
map, that is an affine map taking vertices to vertices, preserving the order.

As in Part I, an auxiliary perturbation data $\mathcal{D}$
for $P$, (in particular) involves: 
\begin{itemize}
	\item A choice of a natural system $\mathcal{U} $, see
	Section \ref{sec:constructionSmallData}, consisting of
	certain maps
	\begin{equation*}
{u}: {\mathcal {S}} ^{\circ}  _{d}  \to \Delta
^{n}, \text{with $d,n$ varying},
\end{equation*} as already discussed in Section
\ref{sec:non-metric} of the Introduction.
	\item Choices of certain Hamiltonian connections, on
	Hamiltonian bundles associated to the maps $u$.
  (Oversimplified for this outline.) 
\end{itemize}
Given such a $\mathcal{D} $, we construct in Part I a functor $$F: \Delta (X) \to A _{\infty}-Cat,$$ 
where $A _{\infty}-Cat $ denotes the category of $A _{\infty} $ categories. 
The properties of this functor are such that we may
algebraically get an induced functor 
$$F ^{unit}: \Delta (X) \to   
 A _{\infty}-Cat ^{unit},  $$ with $A _{\infty}-Cat ^{unit}$ denoting the category of unital $A _{\infty} $ categories,
by taking unital replacements. In what follows we rename $F$
by $F ^{raw} $ and $F ^{unit} $ by $F$, as $F$ is the main
object here.

We then define $$\operatorname {Fuk} _{\infty} (P) = \colim
_{\Delta (X)} NF,  $$ which is shown to be an
$\infty$-category whose equivalence class (under
concordance, see Definition \ref{def:concordancefibration}) is independent of all choices.
This also has the structure of a categorical fibration:
\begin{equation*}
\operatorname {NFuk} (M,  \omega)
\hookrightarrow \operatorname {Fuk} _{\infty} (P)
\xrightarrow{p _{\bullet }} X _{\bullet}, 
\end{equation*}
where $\operatorname {NFuk} (M,  \omega) $ is the $A
_{\infty} $ nerve of the Fukaya category of $M$. We will
extract from the above fibration a Kan fibration and work with that, since then we can just use standard tools of topology.

To this end we have the following elementary lemma.
\begin{lemma} \label{lemma:kanfib}
Suppose we have a categorical fibration $p:
   \mathcal{Y}
\to \mathcal{X}$, where $\mathcal{X} $ is a Kan complex. Let
   $K (\mathcal{Y})$
denote the maximal Kan sub-complex of $\mathcal{Y} $ then $p:
K (\mathcal{Y}) \to \mathcal{X}$ is a Kan fibration.
\end{lemma}
The proof is given in Appendix \ref{appendix:Kan}. 
In particular by the above lemma $$K (P):=K (\operatorname
{Fuk} _{\infty}(P)) \xrightarrow{p _{\bullet }} X _{\bullet }$$ is a Kan fibration. 
\begin{notation}
\label{notation:} In what follows $p _{\bullet }$ will
refer to this projection unless specified otherwise.
\end{notation}

\begin{definition} \label{def:concordancefibration}
We say that a Kan fibration or a categorical
fibration $\mathcal{P} $ over a Kan complex
$\mathcal{X} $ is \textbf{\emph{non-trivial}}
if it is not null-concordant. Here $\mathcal{P}
$ is \textbf{\emph{null-concordant}} means that
there is a Kan respectively categorical
fibration $$\mathcal{Y}  \to \mathcal{X}
\times \Delta ^{1} _{\bullet} ,$$ whose
pull-back by $i _{0}: \mathcal{X}   \to
\mathcal{X}  \times \Delta ^{1} _{\bullet}   $
is trivial and by $i _{1}: \mathcal{X}  \to
\mathcal{X}  \times \Delta ^{1} _{\bullet}   $
is $\mathcal{P} $. Here the two maps $i_{0}, i_{1}$ correspond to the two vertex inclusions $\Delta ^{0} _{\bullet}  \to \Delta ^{1} _{\bullet}   $.
\end{definition}
\begin{theorem} \label{thm:noSection}
   Suppose that $p: P \to S ^{4} $ is a non-trivial Hamiltonian
	 $S ^{2} $ fibration then $p _{\bullet}: K (P) \to S ^{4}
	 _{\bullet}  $ does not admit a section. In particular $K
	 (P)$ is a non-trivial Kan fibration over $S ^{4}
	 _{\bullet}  $ and so $Fuk _{\infty} (P) $ is
	 a non-trivial categorical fibration over $S ^{4} _{\bullet}  $.
\end{theorem} 
This is the main technical result of the paper.
Although in a sense we just are just deducing
existence of a certain holomorphic curve, for this
deduction we need a global compatibility condition
involving multiple moduli spaces, involved in
multiple local datum's of Fukaya categories, so
that this computation will not be straightforward.

The proof will be aided by constructing suitable
perturbation data, and will be split into a number of
sections.  


\section{Qualitative description of the perturbation data}
Let $\operatorname {Fuk} (S ^{2}, \omega  )$ denote the $\mathbb{Z} _{2} $-graded $A _{\infty} $ category over $\mathbb{Q}$, with objects oriented spin Lagrangian
submanifolds Hamiltonian isotopic to the equator. Our
particular construction of $\operatorname {Fuk} (M, \omega
)$ is presented in Part I. In particular, we use the
language of perturbation systems $\mathcal{D}$,  see Section
5, and 6.1 Part I. The data $\mathcal{D} $ is
generally associated to a Hamiltonian fibration, and uses
the language of connections. As a symplectic manifold is a Hamiltonian fibration over a point, we write $\mathcal{D} _{pt}$ for this restricted data, needed for construction of
$\operatorname {Fuk} (S ^{2}, \omega  )$. 

Denote by $\operatorname {Fuk ^{eq}}  (S ^{2}, \omega
) \subset \operatorname {Fuk} (S ^{2}, \omega  )$ the full
sub-category obtained by restricting our objects to be
equators in $S ^{2} $. We take our perturbation data
$\mathcal{D} _{pt} $ so that the following is satisfied. 
\begin{itemize}
	\item All the connections $\mathcal{A} (L, L') $ for $L, L' \in \operatorname {Fuk ^{eq}} (S ^{2})$ are $PU
(2)$-connections.
\item For $L$ intersecting $L'$ transversally, the $PU (2)$ connection $\mathcal{A} (L, L')$ is the trivial flat connection. 
\item For $L=L'$ the corresponding connection is generated by an autonomous Hamiltonian.
\end{itemize}


The associated cohomological Donaldson-Fukaya category $DF  (S ^{2} )$ is equivalent as a linear category over $\mathbb{Q} $ to $FH (L_0, L _{0} )$ (considered as a linear category with one
object) for $L _0 \in \operatorname {Fuk} (S ^{2}, \omega ) $.

It is easily verified that a morphism (1-edge) $f$ is an
isomorphism in the nerve $\operatorname {NFuk}(S
^{2} )$, see Part I for definitions, if and only if it corresponds, under the nerve construction $N$, to a morphism in 
$\operatorname {Fuk} (S ^{2}, \omega ) $ that induces an isomorphism in $DF (S ^{2})$. Such a morphism will be called a \emph{$c$-isomorphism}.

Consequently the maximal Kan subcomplex $K (S ^{2} )$ of
$\operatorname {NFuk}  (S
^{2} ) $ is characterized as the maximal subcomplex with 1-simplices the images by $N$ of $c$-isomorphisms in $Fuk (S ^{2} )$.  

\begin{remark} \label{remark:localization}
   It would be interesting (and likely not too difficult) to identify the homotopy type of   $K
   (S ^{2} ) $. 
\end{remark} 
\subsection{Extending $\mathcal{D} _{pt} $ to higher dimensional simplices}
\label{section:model}
\begin{terminology} A bit of possibly non-standard terminology: we say that $A$ is a \emph{model}
for $B$ in some category, with weak equivalences, if there is a morphism
$mod: A \to B$ which is a weak-equivalence.
The map $mod$ will be called a \emph{modelling map}. 
In our
context the modeling map $mod$ always turns out to be a monomorphism.
\end{terminology}
Let us model $D ^{4} _{\bullet}  $ and $S ^{3} _{\bullet} $  as follows. Take the standard representable 3-simplex $\Delta ^{3} _{\bullet}$, and the standard
representable 0-simplex $\Delta ^{0} _{\bullet}  $. Then collapse all faces of $\Delta
^{3}_{\bullet}  $ to a point, that is take the colimit of the following diagram:
\begin{equation}
 \begin{tikzcd}
   & & \Delta ^{0} _{\bullet}     &  \\
  \Delta ^{2} _{\bullet} \ar [drr, "i_0"] \ar [urr] &  \Delta ^{2} _{\bullet} \ar [dr, "i_1"] \ar [ur]   & \Delta ^{2} _{\bullet}
    \ar [d, "i_2"]  \ar [u] & \Delta ^{2}
   _{\bullet} \ar [dl, "i _{3}"]  \ar [ul]   \\
 & &   \Delta ^{3} _{\bullet}   &  \\  
\end{tikzcd}
\end{equation}
Here $i _{j} $ are the inclusion maps of the non-degenerate 2-faces.
This gives a simplicial set $S _{\bullet} ^{3, mod}  $ modelling the simplicial set $S ^{3} _{\bullet} $, in other words there is a natural a weak-equivalence $$S ^{3,mod} _{\bullet} \to S ^{3} _{\bullet}.  $$


Now take the cone on $S _{\bullet} ^{3, mod}  $, denoted by $C (S ^{3,mod}
_{\bullet})$, and collapse the one non-degenerate 1-edge.
The resulting simplicial set $D
^{4,mod}_{\bullet}  $ is our model for $D ^{4} _{\bullet}  $, it may be identified with a subcomplex of $D ^{4} _{\bullet}  $ so that the inclusion map $mod: D
^{4,mod}_{\bullet} \to D ^{4} _{\bullet}   $ induces a weak homotopy equivalence of pairs 
\begin{equation} \label{eq:homotopyequivalence}
(D^{4,mod}_{\bullet}, S ^{3,mod}  _{\bullet} ) \to  (D ^{4} _{\bullet}, S ^{3}  _{\bullet}). 
\end{equation}
We set $b _{0} \in D ^{4} _{\bullet} $ to be the vertex which is the image by $mod$ of the unique 0-vertex in $D ^{4,mod} _{\bullet}   $.

Suppose we have a commutative diagram:
\begin{equation*}
\begin {tikzcd}
    D ^{4}     \ar[r,"h_+"] &  S ^{4}  &  \arrow {l} [above] {h_-} D ^{4} \\
    & S ^{3} \ar [ul, "i"] \arrow {ur} [below] {i} &  
  \end{tikzcd}
\end{equation*}
where $i: S ^{3} \to D ^{4} $ is the natural boundary
inclusion, and s.t. the following is satisfied.
\begin{itemize}
	\item $h _{\pm}:D ^{4} \to S ^{4} $ are smooth, and their
	images cover $S ^{4} $.
	\item $$h _+ ({D} ^{4})  \cap h _- ({D ^{4} })  $$ is contained in the image $E$ of $$h _{\pm} \circ i: S^{3}  \to S ^{4}. $$
	\item $h _{\pm} $ takes $b _{0} $ to $x _{0} $.
\end{itemize}


For example, we may just let $h _{-}$ represent the
generator of $\pi _{4} (S ^{4}, x _{0}) $ and $h _{+} $ to
be the constant map to $x _{0} $. 
We call such a pair $h _{\pm} $ a \emph{complementary pair}.

We set $$D _{\pm}:= h _{\pm} (D ^{4,mod} _{\bullet}) \subset S ^{4} _{\bullet}   $$ and we set $\Sigma _{\pm} \in S ^{4} _{\bullet} $ to be the image by $h _{\pm} $ of the sole non-degenerate 4-simplex of $D ^{4,mod} _{\bullet}  $. We also set $$\partial D _{\pm} := h _{\pm}(\partial D ^{4,mod} _{\bullet}), $$ where $\partial D ^{4,mod} _{\bullet}$ is the image of the natural inclusion $S ^{3,mod} _{\bullet} \to  D ^{4,mod} _{\bullet}$.

Fix a Hamiltonian frame for the fiber $P_{x_{0}} $ of $P$
over $x _{0} $, in other words a Hamiltonian bundle
diffeomorphism 
\begin{equation*}
\begin{tikzcd}
S ^{2} \ar[r, ""] \ar [d, ""] &  P \ar [d,"p"] \\
pt \ar [r, "x _{0}"]   & X.
\end{tikzcd}.
\end{equation*}
In particular, this allows us to identify
$\operatorname {Fuk} (S
^{2}, \omega  )$ with $F ^{raw} (x _{0})$, using the
analytic perturbation data
$\mathcal{D} _{pt}$ for both. Denote by $x
_{0,\bullet}$ the image of the map $$\Delta ^{0} _{\bullet}
\to S ^{4,mod} _{\bullet},  $$ induced by the inclusion of
the 0-simplex $x _{0}$. 

We continue with the description of the data $\mathcal{D}
= \mathcal{D} (P) $. This must associate certain data
$\mathcal{D} _{\Sigma}$ for each singular simplex $\Sigma: \Delta^{n}
\to S ^{4}$. Recall from Section 8 Part I, that given the data
$\mathcal{D} _{\Sigma}$  for a non-degenerate simplex
$\Sigma $, we assigned extended perturbation data
$\mathcal{D} _{\widetilde{\Sigma}}$  for all
degeneracies $\widetilde{\Sigma} $ of this simplex. 
So by this discussion, our chosen data $\mathcal{D} _{pt} $ induces perturbation data for all degeneracies of $x _{0} $, that is 
for all simplices of $x _{0,\bullet}$, this data will again
be denoted by $\mathcal{D} _{pt} $, for simplicity.

Fix an object $L _{0} \in \operatorname {Fuk ^{eq}} (S ^{2} ) \subset F ^{raw} (x _{0})$. Denote by $\gamma \in \hom _{F ^{raw}  (x_0)} (L _{0}, L _{0}) $ the
generator of $FH _{1}  (L_0, L _{0}) $,  i.e. the fundamental chain, so that it corresponds to the identity in $DF  (L _{0}, L _{0}  )$. This $\gamma$ is uniquely determined by our conditions and corresponds to a single geometric section. Denote by $L ^{i}  _{0} $ the image of $L _{0} $ by the embedding $$F ^{raw}  (x
_{0}) \to F ^{raw} (\Sigma _{+}),$$ corresponding to the $i$'th vertex inclusion into
$\Delta ^{4} $, $i=0, \ldots ,4$. 

Let $m _{i} $ be the edge between $i-1,i$ vertices and set $$\overline{m}_{i}:=\Sigma _{+}  \circ m _{i}.$$ Let $\Sigma ^{0} _{i}  $ denote the
0-simplex obtained by restriction of $\Sigma ^{4} $ to the $i$'th vertex.
Note that each $\overline {m} _{i} $ is degenerate by
construction, so we have an induced morphism $$F ^{raw}
(pr): F ^{raw}  (\overline{m} _{i}) \to F ^{raw}  (x _{0} ),$$
for $pr$ the degeneracy morphism in $\Delta (S ^{4})$:  $$pr: \overline{m} _{i}  \to \Sigma ^{0} _{i}.    $$
Finally, for each $L _{0} ^{i-1}, L _{0} ^{i} $ we have
a $c$-isomorphism $$\gamma _{i}: L_{0} ^{i-1} \to L _{0}
^{i}    $$ in $F ^{raw}  (\overline{m} _{i}) \subset
F ^{raw} (\Sigma _{+} ) $, which corresponds to $\gamma$, meaning that the
fully-faithful projection $F ^{raw}  (pr)$ takes
$\gamma _{i} $ to $\gamma$. We will denote by $\gamma _{i,j} $ the analogous
$c$-isomorphisms $L _{0} ^{i} \to L _{0} ^{j}    $.

\begin{notation}
 Let us denote from now on, the morphism spaces  $hom _{F
 ^{raw} (\Sigma  _{\pm})} (L _{0}, L _{1}) $ by $hom
 _{\Sigma  _{\pm}} (L _{0}, L _{1}) $. And
 denote the $A _{\infty}$ composition maps $\mu ^{d} _{F
 ^{raw} (\Sigma _{\pm})}$, in the $A _{\infty }$ category $F
 ^{raw} (\Sigma _{\pm}) $, by $\mu ^{d} _{\Sigma _{\pm}}$. 
\end{notation}
\begin{definition} \label{def:smalldata} We call perturbation data $\mathcal{D}$ for $P$  \textbf{\emph{small}} if it is extends the data $\mathcal{D} _{pt} $ as above, and if with respect to $\mathcal{D}$ 
 \begin{equation} \label{eq:gammaunit+}
\mu ^{d} _{\Sigma _{+}} (\gamma ^{1}, \ldots,
   \gamma ^{d} )=0, \text{ for $2 <d <4$},  
\end{equation}
   where $(\gamma ^{1}, \ldots, \gamma ^{d}) $ is
	 a composable chain, and each $\gamma ^{k} $ is of the
	 form $\gamma _{i,j} $  as above.
\end{definition}
%
We will see further on how to construct such small data, assume for now that it is constructed.

Let $\{f _{J} \} $,  corresponding to an $n$-simplex, be as
in the definition of the $A _{\infty} $ nerve in Appendix
A.4 Part I, where $J$ is a subset of $[n] = \{0, \ldots, n\}$.

\begin{lemma} \label{lemma:conditionsfJ}
 Let $\mathcal{D}$ be small as above, then there is a 
   a 4-simplex $\sigma \in NF ^{raw}(\Sigma _{+})$ 
with faces determined by the conditions:
\begin{itemize} 
   \item $f _{J} =0$, for $J$ any subset of $[4]$ with at least $3$ elements.
   \item $f _{\{i-1,i\}} = \gamma _{i}$ for $\gamma _{i} $ as before.
\end{itemize}  
\end{lemma}
\begin{proof} 
This follows by \eqref{eq:gammaunit+} and by the identity $\mu ^{2} _{\Sigma _{+}  } (\gamma, \gamma) = \gamma $.
   \end{proof}
If we take our unital replacements so that $\gamma$
corresponds to the unit, then $\sigma$ induces (by the
construction) a section of $K (P _{+}) \to D _{+} $, where $K (P _{\pm} )$ will be shorthand for $K (P)$ restricted over ${D _{\pm}}$.

%
Let
\begin{equation*}
   i: \left(K (P _{+})|_{\partial D_ {+}}:= p _{\bullet} ^{-1} (\partial D_ {+}) \right)   \to K (P_{-}),
\end{equation*}
be the natural inclusion map. 
Set $$sec = i
\circ \sigma \circ h _{+}| _{\partial D ^{4,mod} _{\bullet}   }.  $$


\subsection {The main lemma and immediate consequences} \label{section:immediateconsequences}
\begin{lemma} \label{lemma:nontrivialelement} 
Suppose that $P$ is a non-trivial Hamiltonian fibration and
$\mathcal{D}$ is small data for $P$ as above, then $sec$ as
above does not extend to a section of $K (P _{-})$.
Moreover, small data $\mathcal{D} $ exists.
\end{lemma} 
This lemma involves all the ingredients of our theory,  its
proof that will be broken up in parts, and will follow shortly.
\begin{proof} [Proof of Theorem \ref{thm:firstquantumobstruction}]
Clearly $\qobs (P) \geq 4$, since the $3$-skeleton of $S ^{4} $ is trivial. 
By Lemma \ref{lemma:nontrivialelement} above, $K (P)$ does
not have a section over the $4$-skeleton. 
\end{proof}
\begin{remark}
 When $P$ is obtained by clutching with a generator of $\pi _{3} (PU (2)) $, and when $h _{\pm} $ are embeddings, the class $[sec]$ in $\pi _{3} (K (P _{-} )) \simeq \pi _{3} (K (S ^{2} )) $ can be thought of as ``quantum'' analogue of the class of the classical Hopf map.
\end{remark}
\begin{proof}[Proof of Theorem \ref{thm:noSection}]
It might be helpful to first review Appendix
\ref{appendix:Kan} before reading the following.
If we take any small perturbation data $\mathcal{D}$ for $P$, then the first part follows immediately by Lemma \ref{lemma:nontrivialelement}.
So $K (P)$ is non-trivial as a Kan fibration.
This then implies that $\operatorname {Fuk} _{\infty} (P)$ is non-trivial as
a categorical fibration, which means in particular that its
classifying map $$f _{P}: S ^{4} _{\bullet} \to (\mathbb{S},
\operatorname {NFuk} (S ^{2}, \omega  )) $$ is not null-homotopic.

To see this, suppose otherwise that we have 
a categorical fibration
$$\widetilde{\mathcal{P} } \to S ^{4}
_{\bullet}  \times I _{\bullet},  $$
restricting to $\operatorname {Fuk} _{\infty} (P) $ over $S
^{4} _{\bullet}  \times 0 _{\bullet}  $ and to
$\operatorname {NFuk} (S ^{2}, \omega  ) \times S ^{4} _{\bullet}   $
over the other end $S ^{4} _{\bullet}  \times 1
_{\bullet}  $. 
Here $0 _{\bullet} $,
respectively $1
_{\bullet} $ are notation for the  images of $i _{j, \bullet }: \Delta ^{0} _{\bullet}
\to \Delta^{1}  _{\bullet }$, $j=0,1$, where $i
_{j, \bullet }$ are induced by the pair of
boundary point inclusions.  

Now take the maximal Kan sub-fibration of
   $\widetilde{\mathcal{P} } $, then by Lemma \ref{lemma:kanfib} we obtain a trivialization of $K (P)$ which is a contradiction. 
\end{proof}
\begin{proof} [Proof of Theorem \ref{thm:main}] 
Theorem \ref{thm:noSection} implies that the group
homomorphism $$\mathbb{Z} \simeq \pi _{4} (\operatorname
{BHam}(S ^{2}, \omega ) \xrightarrow{} \pi _{4} (\mathbb{S},
\operatorname {NFuk} (S ^{2}, \omega  )),$$ 
   has vanishing kernel, so that the result follows.
\end{proof}

\section {Towards the proof of Lemma \ref{lemma:nontrivialelement}} \label{section:partI}

We will denote by $L _{0, \bullet} $ the image of the map $\Delta ^{0} _{\bullet} \to K (P _{-}),$ induced by the inclusion of $L _{0}$ into $K (S ^{2} )$ as a 0-simplex.  Suppose that $sec$ extends to a section of $K (P _{-})$, so we have map $$e: D ^{4,mod} _{\bullet} \to K (P_{-}) $$ extending $sec$ over $\partial D ^{4,mod} _{\bullet}  $. 
We may assume WLOG that $e$ lies over $h_{-}$, meaning $$p _{\bullet} \circ e = h _{-}.$$ Since it can be homotoped to have this property. To see this, first take a relative homotopy of 
\begin{equation*}
   p _{\bullet} \circ e:  (D ^{4,mod} _{\bullet}, \partial D ^{4,mod} _{\bullet})   \to (D _{-}, \partial D _{-})
\end{equation*}
to $h _{-} $, using that we have a homotopy equivalence of pairs \eqref{eq:homotopyequivalence},
and then lift the homotopy to a relative homotopy upstairs using the defining lifting property of Kan fibrations.

And so we have a 4-simplex $$T = e (\Sigma ^{4}) \in K (P
_{-} ) $$ projecting to $\Sigma _{-} \in D _{-} $ by $p
_{\bullet} $. Since $T$ is in the image of $e$, all but one
3-faces of $T$ are totally degenerate with image in $L _{0,
\bullet } $. The exceptional 3-face is the sole non-degenerate 3-face of $sec$, (of $sec (\partial D ^{4,mod} _{\bullet})$).

Let $m _{i,j}, \gamma _{i,j} $ be as in the previous
section, but corresponding now to $\Sigma _{-} $ rather then
$\Sigma _{+} $. Then by the boundary condition on $e$, the edges
of $T$ (which are all edges of $sec$)  correspond, under the
nerve construction, to the
generators $\gamma _{i,j}$. As this is the condition for the
edges of $sec$.

\begin{lemma} \label{lemma:corellator}
For $\mathcal{D}$ small as above, and for the unital replacement $F$ of $F ^{raw} $ as above, the simplex $T$ exists if and only if
$\mu ^{4} _{\Sigma _{-}} (\gamma _{1},  \ldots,  \gamma _{4} )$
 is exact. 
\end{lemma}
\begin{proof} 
 The following argument will be over
   $\mathbb{F}_{2} $ as opposed to $\mathbb{Q}$ as the signs will not matter. Recall that we take the unital replacement so that $\gamma \in hom_ {F^{raw} (P _{x _{0}})} (L _{0}, L _{0})$ corresponds to the unit in the unital replacement.

Now if $T \in K (P _{-} )$ as above exists, then it
corresponds under unital replacement (see Remark 7.5 in Part
I) to a $4$-simplex $T' \in NF ^{raw} (\Sigma _{-})$
satisfying the following condition on its $4$-face.
Recalling the nerve construction, the morphism $f _{[4]} \in hom _{\Sigma  _{-}} (L_{0} ^{0},
L _{0} ^{4} ) $, figuring in the definition of the $4$-face,
satisfies:
\begin{equation} \label{eq:simplex}
\mu ^{1} _{\Sigma  _{-}} f _{[4]} = \sum _{1 < i < 4} f _{[4] - i} + \sum _{s} \sum _{(J _{1},
   \ldots J _{s} ) \in decomp
_{s} } \mu ^{s} _{\Sigma  _{-} }  (f_{J _{1}}, \ldots , f _{J _{s} }  ).
\end{equation}
By our conditions on the boundary of $T$, by the condition on the unital replacement, and by the conditions in Lemma \ref{lemma:conditionsfJ}, 
we must have $f_{J} =0 $, for every proper subset $J \subset [4]$, in some length $s$ decomposition of $[4]$, unless $J=\{i,j\}$ in which case  $f _{i,j}=\gamma _{i,j} $.
Given this \eqref{eq:simplex} holds if and only if $\mu ^{4} _{\Sigma _{-}}
(\gamma _{1} , \ldots , \gamma _{4} ) $ is exact. 

\end{proof}
We are going to show that for  small $\mathcal{D}$,   $ \mu ^{4} _{\Sigma _{-}} (\gamma _{1} ,  \ldots,  \gamma
   _{4} )$  does not vanish in homology, which will finish the proof of the Lemma up to construction of small $\mathcal{D} $.  However the calculation  will require significant setup.
\section {Hamiltonian fibrations and taut structures, holomorphic sections and area bounds}
We collect here some preliminaries on moduli spaces of
holomorphic sections of fibrations with Lagrangian boundary
constraints, and the closely related curvature bounds.
There is an apparently new theory here of taut Hamiltonian
structures,  but aside from that much of this material has previously appeared elsewhere, perhaps in less generality. We will eventually need all that is presented in this section, but the reader may only skim on the first reading.
\subsection {Coupling forms} \label{sec:couplingform}
We refer the reader to  \cite[Chapter
6]{citeMcDuffSalamonIntroductiontosymplectictopology} for more details on what
follows. A Hamiltonian fibration is a smooth fiber bundle $$M \hookrightarrow P \to X,$$
with structure group $ \operatorname  {Ham}(M, \omega) $ with its $C ^{\infty} $ Frechet topology. A \textbf{\emph{Hamiltonian connection}} is just an Ehresmann connection for a Hamiltonian fibration. 

Given that $M$ is closed, a \emph{coupling form},
originally appearing in 
\cite{citeGuilleminLermanEtAlSymplecticfibrationsandmultiplicitydiagrams}, for a Hamiltonian
fibration $M \hookrightarrow P \xrightarrow{p} X$, is a closed 2-form $
{\Omega} $ on $P$ whose restriction to fibers coincides with $\omega $  and
which has the property: 
\begin{equation*}  \int _{M}  {\Omega} ^{n+1} =0,
\end{equation*}
with integration being integration over the fiber operation.
Such a 2-form determines a Hamiltonian connection $ \mathcal {A} _{
{\Omega}} $, by declaring horizontal spaces to be $ {\Omega}
$-orthogonal spaces to the vertical tangent spaces. A coupling form generating a
given connection $\mathcal{A}$ is unique. A Hamiltonian connection $
\mathcal {A} $ in turn determines a coupling form $ {\Omega} _{
\mathcal {A}} $ as follows. First we ask that $ {\Omega} _{ \mathcal
{A}} $ induces the connection $ \mathcal {A} $ as above. This determines $
{\Omega} _{ \mathcal {A}} $ up to values on $ \mathcal {A}
$-horizontal lifts $ \widetilde{v},  \widetilde{w} \in T _{p} P $ of $v,w \in T
_{x} X $. We specify these values by the formula
\begin{equation} \label{eq:couplingvalue}
{\Omega} _{
\mathcal {A}} ( \widetilde{v}, \widetilde{w}) = R _{ \mathcal {A}} (v, w) (p),
\end{equation}
where $R _{\mathcal {A}}$ is the lie algebra valued curvature 2-form of $\mathcal{A}$. Specifically, for each $x$, $R _{\mathcal {A}}| _{x} $ is 
a 2-form valued in $C
^{\infty} _{norm} (p ^{-1} (x))$ - the space of 0-mean normalized smooth
functions on $p ^{-1} (x) $.

\subsection {Hamiltonian structures on fibrations} \label{section:exact}
Let $S$ be a Riemann surface with boundary, with punctures
in the boundary, and a fixed structure of strip end charts
at ends, (positive or negative), i.e. a strip end structure as in Part I. 

Let $M
\hookrightarrow \widetilde{{S}} \xrightarrow{pr} {S} $ be
a Hamiltonian fiber bundle, with model fiber a monotone symplectic manifold $(M,\omega)$, with distinguished Hamiltonian bundle trivializations $$[0,1] \times (0, \infty) \times M \to  \widetilde{{S}} $$ at the positive ends,
and with distinguished Hamiltonian bundle trivializations $$[0,1] \times (-\infty, 0) \times M \to \widetilde{{S}}, $$ at the negative ends.
These are collectively called called \textbf{\emph{strip end
charts}}, (slightly abusing terminology).  Given the structure of such bundle trivializations we say that $\widetilde{{S}} $ has
   \textbf{\emph{end structure}}.
\begin{definition} \label{def:respectsendstructure}
    Let  $$\mathcal {L} \subset (\widetilde{S}| _{\partial
		S}= pr ^{-1} (\partial S))   \to \partial S$$ be
		a Lagrangian sub-bundle, with model fiber an object, in
		the sense of Part I, (in particular a spin oriented
		Lagrangian submanifold). 
		We say that $\mathcal{L}$ \textbf{\emph{respects the end
   structure}} if $\mathcal{L}$ is a constant sub-bundle in
	 the strip end chart trivializations above.
\end{definition}
For $\mathcal{L} $ as above, in the strip end chart coordinates at the end $e
_{i} $, let $L ^{j} _{i}  $ denote the fibers (which are by assumption $t$ independent) of
$\mathcal{L}$ over 
$$\{j\} \times \{t\}, \,j=0,1.  $$

We say that a Hamiltonian connection $\mathcal {A}$ on $\widetilde{{S}} $ is \textbf{\emph{compatible}}
with the connections $\{\mathcal{A} _{i} \}$ on $[0,1]
\times M$ at each end $e _{i} $, if in the strip coordinate chart at the $e _{i} $ end,
$\mathcal{A}$ is flat and $\mathbb{R}$-translation invariant and has the form $\overline{\mathcal{A}} _{i}$ 
where $\overline {\mathcal{A}} _{i}   $ denotes its
$\mathbb{R}$-translation  invariant extension of
$\mathcal{A} _{i} $ to $(0, \pm \infty) \times \mathbb{R}
\times M$, depending  on whether the end is positive or
negative. We say that a Hamiltonian connection, $\mathcal
{A}$ on $\widetilde{{S}} $  is
$\mathcal{L}$-\textbf{\emph{exact}} if $\mathcal{A}$ preserves
$\mathcal{L}$ (this means that the horizontal spaces of $\mathcal{A}$ are tangent to $\mathcal{L}$). 


For $\mathcal{A}$ compatible with $\{\mathcal{A} _{i} \}$ as above, a family $\{j _{z} \}$ of fiber wise $\omega$-compatible almost complex structures on
$\widetilde{S} $ will be said to \textbf{\emph{respect the
end structure}} if at each end $e _{i} $, in the strip end chart above, the family $\{j _{z}\} $ is
$\mathbb{R}$-translation invariant and is admissible with
respect to $\mathcal{A} _{i}$, in the sense of Part I,
Definition 5.3. The data $\Theta =(\widetilde{S}, S,
\mathcal{L}, \mathcal{A}, \{j _{z} \} )$, with $\mathcal{A}$
compatible with $\{\mathcal{A} _{i} \}$, $\{j _{z} \}$,
respecting the end structure,  will be called a \textbf{\emph{Hamiltonian structure}}.

We will normally suppress $\{j _{z} \}$ in the notation and elsewhere for simplicity, 
as it will be purely in
the background in what follows, (we do not need to manipulate it explicitly).
%
%
\begin{definition} Let $(\widetilde{S}, S, \mathcal{L}, \mathcal{A})$ be a Hamiltonian structure, we say that a smooth section $\sigma$ of $\widetilde{S} \to S $ is \textbf{\emph{asymptotically flat}} if at each end $e _{i} $ of $S$, $\sigma$ $C ^{1} $-converges to an $\mathcal{A}$-flat section.
   Specifically, in the strip end chart at a positive end, this means that there is a $\mathcal{A}$-flat section $$\widetilde{\sigma}: [0,1] \times (0,\infty) \to [0,1] \times (0,\infty) \times M,$$ so that for every $\epsilon>0$ there is a $t>0$ so that $$d _{C ^{1} } (\widetilde{\sigma}, {\sigma}| _{[0,1] \times [t,\infty)} ) < \epsilon. $$
Similarly for a negative end. 
\end{definition} 
Note that the above definition implies that $$\lim _{s
\mapsto \infty}  {\sigma| _{[0,1] \times
\{s\} }} = \gamma ^{i}, $$ for some $\mathcal{A}
_{i}$-flat sections of $[0,1] \times M $,  where the limit
is the $C ^{1}$ limit. (Similarly for negative ends.) 
So we can say that $\sigma$ 
is \textbf{\emph{asymptotic}} $\gamma ^{i}$ at the $e _{i}$
end, and that $\gamma ^{i}$ is the
\textbf{\emph{asymptotic constraint}}  of $\sigma$ at the $e
_{i}$ end.
\begin{definition}\label{def:relativeclass}
	Given a pair of asymptotically
flat sections $\sigma _{1}, \sigma _{2}  $, with boundary in
$\mathcal{L}$, we say that they have the same
\textbf{\emph{relative class}} if:
\begin{itemize}
	\item They are asymptotic to the same flat sections at
	each end. (In the sense above.) 
	\item They are homologous relative to the boundary
	conditions and relative to the asymptotic constraints at the ends.
\end{itemize}
  The set of relative classes will be denoted by $H _{2}
	^{sec} (\widetilde{S}, \mathcal{L}) $. Since a class $A
	\in H _{2}
	^{sec} (\widetilde{S}, \mathcal{L}) $  is represented by
	a section with determined asymptotic constraints. We can
	say that $A$ has asymptotic constraints.
\end{definition}
\subsubsection{Families of Hamiltonian structures.} \label{section:Family}
\begin{definition} \label{def:HamiltonianStructure}
 A \textbf{\emph{family Hamiltonian structure}} or henceforth just Hamiltonian structure, consists of the following:
\begin{enumerate}
      \item A  smooth, connected, compact, oriented manifold $\mathcal{K}$ with boundary, (or corners).  
      \item For
each $r \in \mathcal{K}$ a Hamiltonian structure
$(\widetilde{S} _{r}, S _{r}, \mathcal{L} _{r}, \mathcal{A}_{r}  )$, 
 such that there are smooth fibrations $$ \widetilde{S}
 \hookrightarrow \widetilde{\textbf{S}} \xrightarrow{p _{1}
 } \mathcal{K},  \quad S \hookrightarrow \textbf{S} \xrightarrow{p} \mathcal{K},  $$
and $\{\widetilde{S} _{r} \}, \{S _{r} \}$ correspond to the
fibers of the first and second fibration, respectively and
such that the following holds:
\begin{itemize}
	\item The second fibration has fiber a Riemann surface, so that $\{S _{r} \}= \{p  ^{-1} r  \}$.
	\item The first fibration is a fibration whose fibers $p
	_{1} ^{-1} (r)$  are themselves the total spaces of smooth
	Hamiltonian fibrations $M \hookrightarrow \widetilde{S}
	_{r} \to S _{r} $, ($\widetilde{S} _{r} \simeq
	\widetilde{S} $),
	such that the structure group of $\widetilde{S} \hookrightarrow  \widetilde{\textbf{S}}
	\xrightarrow{p _{1} } \mathcal{K}$ can be reduced to
	smooth Hamiltonian bundle maps (of $\widetilde{S}$).
\end{itemize}
     To elaborate further, let $$M \hookrightarrow \widetilde{S} \to S$$ be a Hamiltonian $M$-fibration over a Riemann surface $S$. Let $Aut$ denote the group of Hamiltonian $M$-bundle automorphisms of $\widetilde{S} $. Then $\widetilde{\textbf{S}} \xrightarrow{p _{1} } \mathcal{K}$ is the associated bundle $P \times _{Aut} \widetilde{S}$ for some principal $Aut$ bundle $P$ over $\mathcal{K}$.
      \item The strip end charts $$e _{i,r}:  [0,1] \times (0, \infty) \times M \to 
\widetilde{{{S}  }} _{r},  $$ for the positive ends, fit into a Hamiltonian $M$-bundle diffeomorphism onto the image:
\begin{equation} \label{eqE}
\widetilde{e}  _{i}:  [0,1]
   \times (0, \infty) \times \mathcal{K} \times M \to \widetilde{\textbf{S}},
\end{equation}
similarly for the negative ends.
\item \label{item:rinvariant} In case of positive ends, we then have an induced smooth $r$-family of connections $\{e _{i,r} ^{*}  \mathcal{A} _{r} \}$
on $[0,1] \times (0, \infty) \times M$, and an induced smooth $r$-family of Lagrangian subfibrations $\{e _{i,r} ^{-1} \mathcal{L}_{r} \}   $ over $\partial [0,1]
   \times (0,  \infty)$.
We ask that $$\forall r: \{e _{i,r} ^{-1}
   \mathcal{L}_{r} \} = \{0\} \times (0,  \infty) \times L ^{0} _{i} \cup \{1\} \times (0,  \infty) \times L ^{1} _{i},        $$
where $L ^{j} _{i}  $ are as following the Definition
\ref{def:respectsendstructure}.
Furthermore, we ask that $$\forall r: \{e _{i,r} ^{*}
\mathcal{A} _{r} \} = \overline {\mathcal{A}}  _{i} $$ for
$\mathcal{A} _{i}, \overline{\mathcal{A}}_{i}    $ as
previously.  (Similarly for negative ends.) 
   \item \label{property:A} There is a Hamiltonian
	 connection $\mathcal{A}$ on $\widetilde{\textbf{S}} \to
	 \textbf{S} $ that extends all the connections
	 $\mathcal{A} _{r} $ (in the natural sense), and preserves $\textbf{L} := \cup _{r} \mathcal{L} _{r}  $. 
\end{enumerate}  
\end{definition}
We will write  
$\{\widetilde{{S}}_{r}, {S}_{r}, \mathcal{L}_{r},
\mathcal{A} _{r} \} _{\mathcal{K}} $  for this data, $\mathcal{K}$ may be omitted from notation when it is implicit.

Let $\{\widetilde{{S}}_{r}, {S}_{r}, \mathcal{L}_{r},
\mathcal{A} _{r} \} $ be a Hamiltonian
structure. In the notation above, if in addition there exists
a Hamiltonian connection $\mathcal{A}$ on
$\widetilde{\textbf{S}} \to \textbf{S} $ as in Property \ref{property:A}, so that $\Omega _{\mathcal{A}}$ vanishes on $\textbf{L}$, we will say that $\{\widetilde{{S}}_{r}, {S}_{r}, \mathcal{L}_{r}, \mathcal{A} _{r} \}  $ is a \textbf{\emph{hyper taut Hamiltonian structure}}.

%
\subsubsection {Moduli spaces of sections of Hamiltonian
structures} \label{sec:ModuliSpacesHamStructures}
Let $\Theta = (\widetilde{{S}}, S,  \mathcal{L},
\mathcal{A} )$ be a Hamiltonian structure.
For a section $\sigma$ of $\widetilde{S} $ define its
vertical $L ^{2} $ energy or \textbf{\emph{Floer energy}}  by $$e (\sigma):= \int _{S} |\pi _{vert}  \circ d\sigma| ^{2}, $$ $$\pi _{vert}: T \widetilde{S} \to T ^{vert} \widetilde{S}    $$ is the $\mathcal{A}$-projection, for $T ^{vert} \widetilde{S}$ the vertical tangent bundle of $\widetilde{S} $, that is the kernel of the projection $T \widetilde{S} \to TS$. 

As in Part I, let $J (\mathcal{A} ) $ denote the almost
complex structure on $\widetilde{S} $ determined by
$\mathcal{A} $ and $\{j _{z}\}$ naturally as follows.
\begin{itemize}
	\item $J (\mathcal{A} ) $ preserves the $\mathcal{A}
	$-horizontal distribution of $\widetilde{S} $.
	\item The projection map $\widetilde{S} \to S$ is $J (A)
	$-holomorphic.
	\item  The restriction of $J (\mathcal{A} ) $ to each
	fiber $M _{z}$ of $\widetilde{S} $ over $z \in S$ is $j
	_{z}$.
\end{itemize}
We say that $J (\mathcal{A} ) $ is \textbf{\emph{induced}}
by $\mathcal{A} $ $\{j _{z}\}$.

Define  $\overline{\mathcal{M}} (\Theta)$ to be the
Gromov-Floer compactification of the space of  $J
(\mathcal{A})$-holomorphic sections $\sigma$ of
$\widetilde{\mathcal{S}} $, with finite Floer energy, and with
boundary on $\mathcal{L}$. Note that for any
$J({\mathcal{A}})$-holomorphic $\sigma$ we have an identity:
\begin{equation*}
e (\sigma) = \int _{S} \sigma ^{*} \Omega _{\mathcal{A}},   
\end{equation*}
and $\Omega _{\mathcal{A}}$ vanishes on $\mathcal{L}$ by the condition that $\mathcal{A}$ preserves $\mathcal{L}$, so that the standard energy controls apply, to deduce the standard Gromov-Floer compactification structure.

More generally, if $\{\Theta _{r} \}=\{\widetilde{{S}}_{r}, S _{r},  \mathcal{L}_{r}, \mathcal{A}_{r}\} _{\mathcal{K}}$ is a Hamiltonian structure, let $$\overline{\mathcal{M}} (\{\Theta _{r}\})$$ be the Gromov-Floer compactification of the space of pairs $(\sigma, r)$, $r \in \mathcal{K} $ with $\sigma$  a $J (\mathcal{A} _{r})$-holomorphic, finite Floer energy section of $\widetilde{\mathcal{S}}_{r} $, with
boundary on $\mathcal{L}
_{r}$. 

We also denote by $$\overline{\mathcal{M}} ( \{\Theta _{r}
\},A) \subset \overline{\mathcal{M}} (\{\Theta _{r}  \})$$
the subset corresponding to relative class $A \in H _{2}
^{sec} (\widetilde{S}, \mathcal{L}) $ curves, where the
latter is as defined above, and where 
$\widetilde{S} \simeq \widetilde{S} _{r} $, $\forall r \in
\mathcal{K} $.

Let $\{\Theta _{r} =  (\widetilde{{S}}_{r}, S _{r},  \mathcal{L}_{r},
\mathcal{A}_{r})\}$ be a Hamiltonian structure, then for each end $e _{i} $ of $S _{r}$ we have a Floer chain complex  $$
CF (\mathcal{A} _{i}):= CF (L ^{0} _{i}, L ^{1} _{i}, \mathcal{A} _{i}, \{j _{z} \} ),$$ 
(independent of $r$ by part \ref{item:rinvariant} of
Definition \ref{def:HamiltonianStructure}) generated over
$\mathbb{Q}$  by $\mathcal{A} _{i} $-flat sections of $[0,1]
\times M$, with boundary on $L ^{0} _{0}, L ^{0} _{1}  $.
This chain complex is defined as in Section 6.1 of Part I.  
\begin{definition}
We say that $\{\Theta _{r} \}$ is
$A$-\textbf{\emph{regular}} if:
\begin{itemize}
	\item The pairs  $(\mathcal{A} _{i}, \{j _{z} \} )$ are
	regular so that the Floer chain complexes $CF (\mathcal{A}
	_{i})$ are defined.
	\item ${\mathcal{M}} (\{\Theta _{r} \},  A)$  is regular,
(transversely cut out).
\end{itemize}
And we say that $\{\Theta _{r} \}$ is \textbf{\emph{regular}} if it is $A$-regular for all $A$. We say that $\{\Theta _{r} \}$ is $A$-\textbf{\emph{admissible}} if 
there are no elements $$(\sigma, r) \in \overline{\mathcal{M}} (\{\Theta _{r} \},  A),$$  for $r$ in a neighborhood of the boundary of $\mathcal{K}$. 

\end{definition}

\begin{definition}  \label{definition:concordance} 
Given a pair $\{\Theta _{r} ^{i}  \} = \{\widetilde{{S}} _{r} ^{0} , {S}_{r} ^{0},
      \mathcal{L}_{r} ^{0}, \mathcal{A} _{r} ^{0}\} _{\mathcal{K}} $, $i=1,2$, of Hamiltonian structures
we say that they are \textbf{\emph{concordant}} if 
there is a Hamiltonian structure $$\{\mathcal{T} _{r} \}= \{\widetilde{{T}}  _{r},
{T}_{r} , \mathcal{L}'_{r}, \mathcal{A}' _{r}
   \} _{\mathcal{K} \times [0,1]} ,$$ with
an oriented diffeomorphism (in the natural sense, preserving all structure)
$$  \{\widetilde{{S}} _{r} ^{0} , {S}_{r} ^{0} ,
      \mathcal{L}_{r} ^{0}, \mathcal{A} _{r} ^{0} 
      \} _{ \mathcal{K} ^{op}}     \sqcup \{\widetilde{{S}} _{r} ^{1}, {S}_{r} ^{1} ,
      \mathcal{L}_{r} ^{1}, \mathcal{A} _{r} ^{1}
      \} _{ \mathcal{K}} 
      \to  \{\widetilde{{T}} _{r},
{T} _{r}, \mathcal{L}'_{r}, \mathcal{A}' _{r
      } \} _{\mathcal{K} \times \partial I} ,$$ 
      where $op$ denotes the opposite orientation for $\mathcal{K}$.

\end{definition}
\begin{definition} We say that a Hamiltonian structure 
   $\{\Theta _{r} \} = \{\widetilde{{S}} _{r}, {S}_{r},
   \mathcal{L}_{r}, \mathcal{A} _{r} \} _{\mathcal{K}} $ is \textbf{\emph{taut}} if for any pair $r _{1}, r _{2} \in \mathcal{K}  $, $\Theta _{r _{1} }$ is concordant to $\Theta _{r _{2}} $ by a concordance $\{\widetilde{{T}}  _{\tau},
{T}_{\tau}, \mathcal{L}'_{\tau}, \mathcal{A}' _{\tau}
   \} _{[0,1]}$ which is a hyper taut Hamiltonian structure.
%
%
\end{definition}

\begin{definition}  \label{definition:isotopy} 
Given an $A$-admissible pair $\{\Theta _{r} ^{i}  \} $, $i=1,2$, of Hamiltonian structures, we say that they are $A$-\emph {\textbf{admissibly concordant}} if there
is an $A$-admissible Hamiltonian structure $$\{\widetilde{{T}}  _{r},
{T}_{r} , \mathcal{L}'_{r}, \mathcal{A}' _{r}
   \} _{\mathcal{K} \times [0,1]}, $$ which furnishes a concordance. If this concordance is in addition a taut Hamiltonian structure, then we say that these pairs are \textbf{\emph{$A$-admissibly taut concordant}}. 
\end{definition}

\begin{lemma} \label{lemma:gluing} 
Let $\{\Theta _{r}\} = \{\widetilde{{S}} _{r}, {S}_{r}, \mathcal{L}_{r}, \mathcal{A} _{r} \}$ be $A$-regular and $A$-admissible, with $S _{r} $ having one distinguished
negative end $e _{0} $, and let $\gamma _{0} $ be the
asymptotic constraint of $A$ at the $e _{0} $ end. 
Define $$ev _{A} = ev (\{\Theta _{r} \}, A)  = \# {\mathcal{M}} (\{\Theta _{r} \},
A) \cdot \gamma _{0} \in CF (\mathcal{A} _{0}), $$ 
where $\# {\mathcal{M}} (\{\Theta _{r} \},  A)$ means signed count of elements when the dimension is 0, and is otherwise set to be zero.
Furthermore, suppose that $CF (\mathcal{A} _{0} ) $ is perfect.
%
Then  $ev _{A}$ is a cycle  and its homology class  depends only on the $A$-admissible concordance class of $\{\Theta _{r} \}$.
\end{lemma}
\begin{proof}

Suppose we are given an $A$-admissible concordance (which we may assume to be regular)
$$\mathcal{T} = \{\widetilde{{T}}  _{r}, {T}_{r},
\mathcal{L}'_{r}, \mathcal{A}' _{r} \} _{\mathcal{K} \times
[0,1]},  $$ 
   between Hamiltonian structures $\{\Theta _{r} ^{0}\}$ and $\{\Theta _{r} ^{1}\}$.  
Then we get a one dimensional compact moduli space
$\overline{{\mathcal{M}}} (\{\widetilde{{T}}  _{r}, {T}_{r},
\mathcal{L}'_{r}, \mathcal{A}' _{r} \}, A)$. By assumption
on the perfection of $CF (\mathcal{A} _{0})$, boundary contributions from Floer degenerations cancel out, so that the 
boundary is:
\begin{equation*}
\partial {\overline{\mathcal{M}}} (\{\widetilde{{T}}  _{r}, {T}_{r},
\mathcal{L}'_{r}, \mathcal{A}' _{r} \}, A) = {\mathcal{M}}
(\{\Theta ^{0} _{r} \} ^{op}, A)  \sqcup  {\mathcal{M}}
(\{\Theta ^{1} _{r} \}, A)
\end{equation*}
where $op$ denotes opposite orientation. From which the result follows.
\end{proof}
\subsection {Area of fibrations} \label{section:areafib}
\begin{definition} \label{def:area}
For a Hamiltonian connection $\mathcal{A} $ on a bundle $M
\hookrightarrow \widetilde{S}  \to S$, with $S$ a Riemann
surface, define a 2-form $\alpha _{\mathcal{A} }$ on $S$ by:
\begin{equation} \label{eq:alphaA} \alpha _{\mathcal{A}}  (v, jv) :=  |R _{\mathcal{A}} (v,jv)| _{+},
\end{equation}
where $v \in T _{z} S $, $R _{\mathcal{A}} (v,w)$ as before
identified with a zero mean smooth function on the fiber
$\widetilde{S} _{z}$ over $z$ and where $|\cdot| _{+} $ is
operator: $|H| _{+} = \max _{\widetilde{S} _{z}} H,$ i.e.
the ``positive Hofer norm''.
\end{definition} 
And define
\begin{align} \label{eqDefArea}
\area  (\mathcal{A})  :=  \int _{S} \alpha _{\mathcal{A}}.
\end{align} 
Note that if $\Omega _{\mathcal{A}} $ is the coupling form
of $\mathcal{A}$, as before, then
${\Omega} _{\mathcal{A}}   + \pi ^{*} (\alpha _{\mathcal{A}})$ is \emph{nearly
symplectic}, meaning that
\begin{equation*} 
   \forall z \in S \, \forall v  \in T _{z} S: ({\Omega} _{\mathcal{A}}   + \pi ^{*} (\alpha)) ( \widetilde{v},
\widetilde{jv}) \geq 0,
\end{equation*}
where $ \widetilde{v}, \widetilde{jv} $ are the
$\mathcal{A}$-horizontal lifts of  $v, jv \in T _{z} S $.

Note that $\area (\mathcal{A})$ could be infinite if there are no constraints on $\mathcal{A}$ at the ends.

\begin{lemma} \label{lemma:energypositivity1}
Let $(\widetilde{S}, S, \mathcal{L}, \mathcal{A} )$ be
a Hamiltonian structure. For $\sigma \in \overline{\mathcal{M}} (\widetilde{{S}}, S,  \mathcal{L},
\mathcal{A})$
 we have  
  \begin{equation*}
    - \int _{S} \sigma ^{*}
    {\Omega} _{\mathcal{A}}     \leq \area  
  (\mathcal{A}).
  \end{equation*}
 \end{lemma}
\begin{proof}
We have $$\int _{S} \sigma ^{*}
    ({\Omega} _{\mathcal{A}} + \pi ^{*} \alpha) \geq 0, $$ 
whenever ${\Omega} _{\mathcal{A}}   + \pi ^{*} (\alpha)$ is nearly symplectic, by the defining properties of $J _{\mathcal{A}} $ and by $\sigma$ being $J _{\mathcal{A}} $-holomorphic.
From which our conclusion follows.
\end{proof}

\begin{lemma} \label{lemmaInvariantPairing} 
Let $\{(\widetilde{S} _{t}, S _{t},
\mathcal{L}_{t}, 
\mathcal{A} _{t}) \} _{[0,1]} $ be a taut concordance.
Let $\sigma _{j} $, $j=0,1$ be asymptotically flat sections of $\widetilde{S} _{j}$ in relative class $A$. 
Then
%
$$
   - \int _{S _{1} } \sigma _{1}  ^{*}  {\Omega} _{\mathcal{A} _{1} }  =  - \int _{S _{0} } \sigma _{0}  ^{*}  {\Omega} _{\mathcal{A} _{0}},
$$
whenever both integrals are finite. In particular, for a Hamiltonian structure  $(\widetilde{S}, S,
\mathcal{L}, 
   \mathcal{A})$,
$\int _{S} \sigma  ^{*}  {\Omega} _{\mathcal{A}}$ depends only on the relative class of $A$, whenever the integral is finite.
%
%
%
%
 \end{lemma}
\begin{proof}
 By the hypothesis, there is a connection $\mathcal{A}$ on
 $\widetilde{\textbf{S}} $, extending each $\mathcal{A} _{t}
 $ and such that $\Omega _{\mathcal{A}} $ vanishes on $ \textbf{L} \subset \widetilde{\textbf{S}}. $  The first part then follows by Stokes theorem. 
Here are the details. For $\sigma _{j} $ as above and for
each end $e _{i} $, cut off the part of the section $\sigma
_{j} $ lying over $[0,1] \times (t _{\delta _{1}, \delta
_{2} }, \infty)$ in the corresponding strip end chart at the end. Here $t _{\delta _{1}, \delta _{2} } $ is such that $\sigma _{0}| _{[0,1] \times \{t\}}   $ is $C ^{1} $ $\delta _{1} $-close to $\sigma _{1}| _{[0,1] \times \{t\}} $ for all $t > t _{\delta _{1}, \delta _{2} } $ and for each end, and is such that $$\int _{[0,1] \times (t _{\delta _{1}, \delta _{2} }, \infty)} \sigma _{j} ^{*}| _{[0,1] \times (t _{\delta _{1}, \delta _{2} }, \infty)}  {\Omega} _{{\mathcal{A} _{j} }} < \delta _{2}, \, j=1,2,    $$ for each end $e _{i} $.
Call the sections with the ends cut off as above by $\sigma _{j} ^{\delta _{1}, \delta _{2}  }  $, they are sections over the compact surfaces $S ^{cut} _{j} $, with ends correspondingly cut off. Then by Stokes theorem, using that ${\Omega} _{{\mathcal{A}}}  $ is closed and using the vanishing of ${\Omega} _{{\mathcal{A}}} $ on $ \textbf{L}   $:   for each $\epsilon$ there exists $\delta _{1}, \delta _{2}  $ such that
   $$\int _{S ^{cut} _{1}} (\sigma _{1} ^{\delta _{1}, \delta _{2}}) ^{*} {\Omega} _{{ \mathcal{A}}} -   \int _{S ^{cut} _{0}} (\sigma _{0} ^{\delta _{1}, \delta _{2}}) ^{*} {\Omega} _{{ \mathcal{A}}} < \epsilon,$$
and $$\int _{S ^{cut} _{j}} (\sigma _{j} ^{\delta _{1}, \delta _{2}}) ^{*} {\Omega} _{{ \mathcal{A} _{j} }} - \int _{S _{j} } \sigma _{j}  ^{*}  {\Omega} _{\mathcal{A} _{j} }  <\epsilon, \, j=1,2.  $$

The last part of the lemma follows from the first. For if
$\mathcal{A}$ preserves $\mathcal{L}$ then $\Omega
_{\mathcal{A}} $ vanishes on $\mathcal{L}$, and consequently
the corresponding constant concordance: $$\{(\widetilde{S}
_{t}, S _{t}, \mathcal{L}_{t}, \mathcal{A} _{t}) \}
_{[0,1]}, \quad \forall t \in [0,1]: (\widetilde{S}_{t}, S _{t}, \mathcal{L}_{t}, \mathcal{A} _{t}) = (\widetilde{S}, S, \mathcal{L}, \mathcal{A})   $$ is taut.
\end{proof}
\begin{definition}
For $\sigma$ a relative class $A$ section of $\Theta=(\widetilde{{S}}, {S}, \mathcal{L}, \mathcal{A}) $ let us call:
\begin{equation*}
-\int _{S} \sigma  ^{*}  {\Omega} _{\mathcal{A}},
\end{equation*}
the $\mathcal{A}$-\textbf{\emph{coupling area}} of $\sigma$, denoted by $carea (\Theta, \sigma)$,
 we may also write $carea (\Theta, A)$ for the same quantity. By the lemma above this is an invariant of the taut concordance class of $\Theta$.
\end {definition}
%
\begin{definition} \label{def:small}
Given a Hamiltonian structure $\Theta=(\widetilde{S},S, \mathcal{L}, \mathcal{A})$  we will say that 
 $\Theta$ is  $A$-\textbf{\emph{small}} if $$\area (\Theta)< carea (\Theta,A).$$ Similarly, given a taut Hamiltonian structure $\{\widetilde{{S}}_{r}, {S}_{r}, \mathcal{L}_{r}, \mathcal{A} _{r} \} _{\mathcal{K} } $ we say that it is 
is \textbf{\emph{$A$-small near boundary}} if
each $(\widetilde{{S}}_{r}, {S}_{r}, \mathcal{L}_{r}, \mathcal{A} _{r})$ is $A$-small for $r$ in a neighborhood of $\partial \mathcal{K}$. 
\end{definition}

\begin{lemma} \label{lemma:smallempty}
   Suppose that $\Theta=(\widetilde{{S}}, {S}, \mathcal{L}, \mathcal{A})  $ is $A$-small then $\overline{\mathcal{M}} (\Theta, A)$ is empty. Or as a contrapositive, if $\overline{\mathcal{M}} (\Theta, A)$ is non-empty then:
\begin{equation*}
\carea (\Theta,A) \leq \area (\Theta).
\end{equation*}
\end{lemma}
\begin{proof} 
This is just a reformulation of Lemma \ref{lemma:energypositivity1}.
\end{proof}

\begin{lemma} \label{lemma:Aadmissible} Let $\{\widetilde{{S}}_{r}, {S}_{r}, \mathcal{L}_{r}, \mathcal{A} _{r} \} _{ \mathcal{K}} $ be a taut Hamiltonian structure with $\mathcal{K}$ connected, so that in particular, for each $r$, $\Theta _{r}= (\widetilde{{S}}_{r}, {S}_{r}, \mathcal{L}_{r}, \mathcal{A} _{r})$ is taut concordant to a fixed $\Theta$.
Suppose that $\{\widetilde{{S}}_{r}, {S}_{r}, \mathcal{L}_{r}, \mathcal{A} _{r} \} _{ \mathcal{K}} $ is $A$-small near boundary then $\{\widetilde{{S}}_{r}, {S}_{r}, \mathcal{L}_{r}, \mathcal{A} _{r} \} _{ \mathcal{K}} $ is $A$-admissible for all $A$ such that $carea (\Theta,A)>0$.
\end{lemma}
\begin{proof} 
Follows immediately by the lemma above.
\end{proof}
\subsection {Gluing Hamiltonian structures} \label{section:gluing}


Let $\mathcal{D}$ denote the Riemann surface which is topologically  $D ^{2} - z _{0}$, $z _{0} \in \partial D ^{2} $,
endowed with a choice of a strip end chart at the end (positive or negative depending on context). The complex structure $j$ here is as induced from $\mathbb{C}$ under the assumed embedding $D ^{2} \subset \mathbb{C} $.

Let $(\widetilde{S}, S, \mathcal{L}, \mathcal{A})$ be
a Hamiltonian structure. We may cap off some of the open
ends $\{e _{i} \} _{i=0} ^{n} $ of ${S}$, by gluing at the ends copies of $\mathcal{D} $ with oppositely signed end.
More explicitly, in the strip coordinate charts at some, say positive, end 
$e _{i} $ of ${S}$, 
excise $[0,1] \times (t, \infty)$ for some $t>0$, call the resulting surface ${S} - e _{i} $. Likewise excise the negative end of $\mathcal{D} $, call this surface $\mathcal{D} - end$. Then glue ${S} -e _{i}  $ with  $\mathcal{D}-end$, 
along their new smooth boundary components. 
Let us denote the capped off surface by ${S} ^{/i}  $.  

Since $\widetilde{{S}}$ is naturally trivialized at the
ends, we may similarly cap off $\widetilde{\mathcal{S}} _{r}
$ over the $e _{i} $ end by gluing with the bundle
$\mathcal{D} \times M $ at the end, obtaining a Hamiltonian $M$ bundle $\widetilde{{S}}  ^{/i}
$ over ${S} ^{/i}$. 

Moreover, we have a certain gluing operation of Hamiltonian structures. In the case of ``capping off'' as above we glue $\Theta = (\widetilde{S}, S, \mathcal{L}, \mathcal{A})$ with the Hamiltonian structure $\Theta'  = (\mathcal{D} \times M, \mathcal{D}, \mathcal{L}', \mathcal{A}')  $ at the $e _{i} $ end, provided $\mathcal{A}'$ is compatible with the connection $\mathcal{A} _{i} $, in the sense of Section \ref{section:exact}, and provided 
$\mathcal{L}$ is compatible with $\mathcal{L}'$. The latter
means that $L ^{j} _{i} = {L'} ^{j} _{i}   $ where these are
Lagrangians corresponding to the strip end chart trivialization of $\mathcal{L}, \mathcal{L}'$ at the corresponding ends, as in Definition \ref{def:respectsendstructure}.

Let us name the result of this capping off $\Theta \# _{i} \Theta' $. The following is immediate:
\begin{lemma} \label{lemma:immediategluing} Suppose that $\{\Theta _{r} \} _{\mathcal{K}} $, $\{\Theta' _{r} \} _{\mathcal{K}} $ with $\Theta' _{r} = (\mathcal{D} \times M, \mathcal{D}, \mathcal{L}' _{r} , \mathcal{A}' _{r} ) $ are taut Hamiltonian structures. Then:
\begin{equation*}
\{\Theta _{r} \# _{i}   \Theta' _{r}  \} _{\mathcal{K}} 
\end{equation*}
  is taut, whenever the gluing operation is well defined, that is whenever we have compatibility of connections and Lagrangian sub-fibrations at the corresponding end.
\end{lemma}

\begin{definition} \label{defInducedLoop}
Let $\pi: \mathbb{R} \to [0,1]$ denote the continuous
retraction map, sending $(-\infty, 0]$ to $0$, and sending
$[1,\infty)$ to $1$. Assuming the end $e _{0}$  of
$\mathcal{D} $ is positive, and using the coordinates of the
strip end chart $e _{0}: [0,1] \times (0, \infty) \to
\mathcal{D} $, fix the following parametrization $\zeta$ of the
boundary of $\mathcal{D} $.
$\zeta: \mathbb{R} \to \partial \mathcal{D}$, satisfies $\zeta (t) \in \{0\} \times (0, \infty)$ for $t \in (-\infty, 0) $, and 
$\zeta (t) \in \{1\} \times (0, \infty)$ for $t \in 
 (1, \infty) $. Given a smooth path 
$$p: [0,1] \to
Lag (M) $$ constant near $0,1$, let $\mathcal{L} _{p} \subset \partial \mathcal{D} \times M  $ denote the Lagrangian subfibration over $\partial \mathcal{D}$, with fiber over $r \in \partial \mathcal{D}$ given by $p \circ \pi (r)$. We say that a Lagrangian subfibration $\mathcal{L}$ as above is
\textbf{\emph{determined by $p$}} if 
$\mathcal{L} = \mathcal{L} _{p}$, after a fixed choice of parametrization of boundary of $\mathcal{D}$ by $\mathbb{R}$.
(In the case the end of $\mathcal{D} $ is negative, the
above is meant to be analogous.) 
\end{definition}
A Hamiltonian connection $\mathcal{A}$ on $[0,1] \times M$
uniquely corresponds to a choice of a smooth function $H:
[0,1] \times M \to \mathbb{R}$, normalized to have mean zero
at each moment. For the holonomy path of $\mathcal{A}$ over
$[0,1] $ is a path $\phi _{\mathcal{A}}: [0,1] \to Ham
(M,\omega)  $, generated by a Hamiltonian $H:
[0,1] \times M \to \mathbb{R}$, and this uniquely determines
the connection. Conversely, $H$ uniquely  determines
a Hamiltonian connection with holonomy path generated by
$H$. We can say that $H$ \textbf{\emph{generates
$\mathcal{A}$}}. 
\begin{lemma} \label{lemma:areadisk} Let $p$ and
$\mathcal{L} _{p}  \subset \partial \mathcal{D} \times
M  $ be as in definition above with $L ^{\pm}
(\widetilde{p})=\rho$, where $\widetilde{p} $ is some lift
of $p$ to $\operatorname {Ham} (M,\omega)$, that is $p (t)
= \widetilde{p} (t) (p (0)) $. Let $\mathcal{A} _{0}   $ be
a Hamiltonian connection on $[0,1] \times M$, generated by
a Hamiltonian $H: [0,1] \times M \to \mathbb{R} ^{} $ with
$L ^{\pm} $ length $\kappa$, constant for $t$ near $0,1$.
Then there is a Hamiltonian connection
$ \widetilde{{\mathcal{A}}} ^{p}_{0} $ on $\mathcal{D}
\times M$,  preserving $\mathcal{L} _{p} $, compatible with
respect to $\mathcal{A} _{0} $, and satisfying $$area
(\widetilde {\mathcal{A}} ^{p} _{0} ) \leq \kappa + \rho. $$
The construction is natural in the sense that
$(\widetilde{p},\mathcal{A} _{0}) \mapsto
\widetilde{{\mathcal{A}}} ^{p}_{0}  $ can be made into
a smooth map (of Frechet manifolds).
\end{lemma}
\begin{proof} Let $q: [0,1] \to \operatorname {Ham}
(M,\omega)$ be the holonomy path of $\mathcal{A} _{0} $, $q
(0)=id$, generated by $H$. Let $\widetilde{p}  \cdot q$ be
the usual path concatenation in diagrammatic order, and $H'$ be its generating Hamiltonian.

   Define a coupling form ${\Omega}'$ on $D ^{2} \times M$:
\begin{equation*}
{\Omega}' =  \omega - d (\eta (rad) \cdot H' d \theta),
\end{equation*}
for $(rad, \theta)$ the modified angular coordinates on $D ^{2}  $, 
$\theta \in [0,1]$, $0 \leq rad \leq 1$, and 
$\eta: [0,1] \to [0,1]$ is a smooth function satisfying $$0 \leq \eta' (rad),$$
and 
\begin{equation} 
\eta (rad) = \begin{cases} 1 & \text{if } 1 -\delta \leq rad \leq 1 ,\\
rad ^{2}  & \text{if } rad \leq 1-2\delta,
\end{cases}
\end{equation}
for a small $\delta >0$.  
By an elementary calculation  $$\area  (\mathcal{A}') = L ^{+} (p \cdot q) = L ^{+} (p)  + L ^{+} (q) , $$ where $\mathcal{A}'$ is the connection induced by $\Omega'$.
 Set $$arc=\{(1,\theta) \in D ^{2} \,|\, 0 \leq \theta \leq
 1/2\}.$$ Let $arc ^{c} $ denote the complement of $arc$ in
 $\partial D ^{2} $.   Fix a smooth embedding $i: D ^{2}
 \hookrightarrow \mathcal{D}$ such that the following is
 satisfied (see Figure \ref{figure:arc}):
\begin{itemize}
	\item The image of the embedding contains
$\partial \mathcal{D} - end$, where $end$ is the image of the
   distinguished (say positive) strip end chart $$[0,1] \times (0,\infty) \to \mathcal{D}.$$
	 \item $i(arc) \subset end ^{c} $,
	 \item $i (arc ^{c}) \subset end$. 
\end{itemize}
 \begin{figure} [h]
 \includegraphics[width=1.5in]{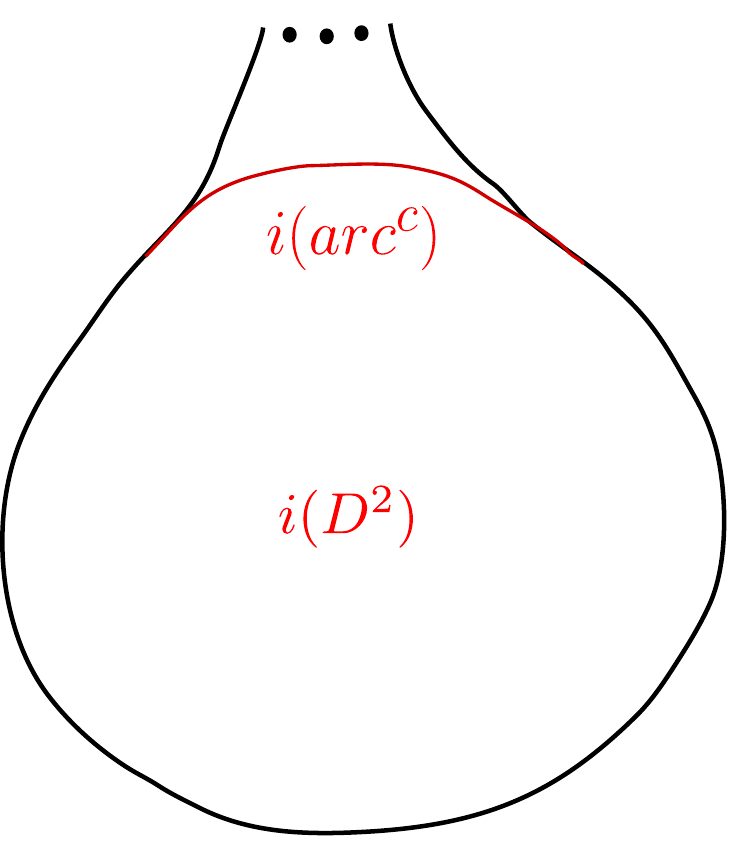}
 \caption {} \label{figure:arc}
\end{figure} 

Next fix a deformation retraction $ret$ of $\mathcal{D}$
onto  $i (D ^{2}) $, so that in the strip end chart above,
for $r \geq 1$ $ret$ is the composition $i \circ param \circ
pr$, where $$pr: [0,1] \times (0,\infty) \to [0,1]$$ the
projection and where $$param: [0,1] \to arc ^{c} \subset
D ^{2} $$ is a diffeomorphism.
   Finally, set ${\Omega} = ret ^{*} {\Omega}'$ on $ \mathcal{D} \times S ^{2} $, and set $\widetilde{\mathcal{A}} _{0} ^{p}  $ to be the induced Hamiltonian connection. As constructed $ \widetilde{{\mathcal{A}}} _{0} $ will be compatible with $\mathcal{A} _{0} $,  when the end of $\mathcal{D}$ is positive.  When the end is negative we take the reverse paths $p ^{-1}, q ^{-1}  $.
\end{proof}
Let us denote by $\Theta (p, \mathcal{A} _{0})
= (\mathcal{D} \times M, \mathcal{D},
\mathcal{L} _{p}, \widetilde{{\mathcal{A}}} ^{p}_{0})   $ the Hamiltonian structure as in the lemma above. When $p$ is the
constant map to $L$  we will instead write 
\begin{equation} \label{eq:ThetaL}
\Theta (L, \mathcal{A} _{0}). 
\end{equation}

The following says that under suitable conditions the
connection of the lemma above can be made to have $area$ 0.
\begin{lemma} \label{lemma:areadisk2} Let $H: M \times [0,1]
\to \mathbb{R} ^{} $ be a smooth time-dependent function
with zero mean at each moment. Let $p: [0,1] \to
\operatorname {Ham} (M,  \omega) $ be the path generated by
$H$. Let $L \in Lag (M, L _{0} )$,
and $p _{L}: [0,1] \to Lag (M, L _{0} )$ be the path $p _{L}
(t) = p (t) (L) $. Let $\mathcal{L} _{p} \in \partial
\mathcal{D} \times M$ be as in the Lemma
\ref{lemma:areadisk}. Let $\mathcal{D} $ have the positive
end $e _{0}$. And let $\mathcal{A} _{0}$ at the $e _{0}$ be
generated by $H$, then there is a Hamiltonian connection
$ \widetilde{{\mathcal{A}}} ^{H} $ on $\mathcal{D}
\times M$,  preserving $\mathcal{L} _{p} $, compatible with
respect to $\mathcal{A} _{0} $, and satisfying $$area
(\widetilde {\mathcal{A}} ^{H} ) = 0.$$
The construction is natural in the sense that
$H \mapsto \widetilde{{\mathcal{A}}} ^{H}  $ can be made into a smooth map (of Frechet manifolds).
\end{lemma}
\begin{proof} We only sketch the proof as the idea is similar to the proof of Lemma \ref{lemma:areadisk}.
In the notation of the proof of Lemma \ref{lemma:areadisk}
let $r: i (D ^{2}) \to i (arc) $ be a smooth retraction.
Set $\mathcal{A}  := r ^{*} \mathcal{A} _{0}  $, and set
$\widetilde{\mathcal{A}} ^{H} = ret ^{*} \mathcal{A} $ (for $ret$ as
before).
\end{proof}
For future use, we denote by 
\begin{equation} \label{eq:ThetaH}
\Theta (H) = (\mathcal{D} \times M, \mathcal{D},
\mathcal{L} _{p}, \widetilde{{\mathcal{A}}} ^{H}),
\end{equation}
the Hamiltonian structure as in the lemma above. 

Now let $\Theta = (\widetilde{S}, S, \mathcal{L},
\mathcal{A})$ be a Hamiltonian structure.  For simplicity,
suppose that $\mathcal{L} $ is trivial with fiber $L _{0}$,
and that $\mathcal{A} $ is trivial over the boundary.
Suppose further that at the end $e _{i} $ the corresponding connection $\mathcal{A} _{i} $ is generated
by $L ^{\pm} $-length $\kappa _{i} $ Hamiltonian $H
_{i}$. By capping each $e _{i}$ end with $\Theta (L _{0},
\mathcal{A} _{i})  $ (keeping in mind that negative-positive
distinction)  we obtain a Hamiltonian structure we
call $\Theta ^{/} = (\widetilde{S} ^{/}, S, \mathcal{L}
^{/}, \mathcal{A} ^{/})  $.
By the lemma above:
\begin{equation} \label{eq:areakappai}
\area (\mathcal{A} ^{/}) \leq \area (\mathcal{A})
+ \sum _{i} \kappa   _{i}.
\end{equation}

\begin{lemma} \label{lemma:gluinglowerbound} 
Let $L _{0} \subset M$ be a monotone Lagrangian submanifold
with monotonicity constant $const>0$. Meaning that for
a relative class $A \in H _{2} (M, L _{0})$: $\omega (A) = const
\cdot \mu (A)$, $\mu$ the Maslov number. Let
$$\Theta:=\{\Theta _{r} \}:=\{\widetilde{S} _{r}, S _{r},
\mathcal{L} _{r} , \mathcal{A} _{r} \} _{\mathcal{K}} $$ be
a hyper taut Hamiltonian structure satisfying:
\begin{itemize}
	\item $\mathcal{K}$ is connected.
	\item $\mathcal{L} _{r} $ is the trivial bundle with fiber
	$L _{0} $ for each $r$.
	\item $\mathcal{A} _{r} $ is the trivial connection over the
	boundary of ${S} _{r} $ for each $r$.
	\item The Floer chain complex  $CF(\mathcal{A} _{i} )$ is perfect for each $i$ and $\mathcal{A} _{i} $ is generated by a time dependent Hamiltonian $H _{i} $ with $L ^{\pm} $ length $\kappa _{i} $.
\end{itemize}
Let $\Theta ^{/} _{r}   = (\widetilde{S} ^{/} _{r}, S ^{/} _{r} , \mathcal{L} ^{/} _{r}, \mathcal{A} ^{/} _{r} ) $ be obtained from $\Theta _{r} $ by capping off each end 
$e _{i} $, so that \eqref{eq:areakappai} is satisfied.
For a given $A \in H _{2} ^{sec} (\widetilde{S}, \mathcal{L})  $, if 
\begin{equation*}
\forall r: \area (\mathcal{A} _{r} ) < carea (\Theta _{r} ^{/}, A ^{/} )  - \sum _{i} \kappa _{i},
\end{equation*} 
where $A ^{/} $ is the capping off of $A$ is described in the proof, then 
$\overline{\mathcal{M}} (\{\Theta _{r}\} , A)$ is empty. 
Moreover, $$\forall r: carea (\Theta _{r} ^{/}, A ^{/}) = -const \cdot Maslov ^{vert}   (A ^{/}),$$
where $Maslov ^{vert} $ is as in Appendix \ref{appendix:maslov}.
\end{lemma}
\begin{proof} 
Suppose otherwise that we have an element $(\sigma _{0},r _{0} ) \in \overline{\mathcal{M}} (\{\Theta _{r}\} , A)$. 
Suppose for the moment that $\overline{\mathcal{M}} (\{\Theta _{r}\} , A)$ is regular. 
There is a morphism (cf. Albers~\cite{citeAlbersPSS})  $$PSS: QH (L) \to FH (L, L),$$ where the right hand side is defined using our construction in terms of flat sections, 
and the left hand side is interpreted for example as the
homology of the Pearl complex,
Biran-Cornea~\cite{citeBiranCorneaLagQuantHomology}. Moreover, as shown by Albers this is an isomorphism in the present monotone context.

We won't give the full construction of this morphism in our
setting, as it just a reformulation of the construction
in \cite{citeAlbersPSS}. Here is a quick sketch. 
Let $$\Theta _{-}  = (\mathcal{D} \times M, \mathcal{D}, \mathcal{L}, \mathcal{A} _{-} ), $$ be the Hamiltonian structure with $e _{0} $ being a negative end, $\mathcal{L}$ trivial with fiber $L$ (which is an object as before), and $$\mathcal{A} _{-}:= \widetilde{\mathcal{A}} ^{p=const}  _{0},  $$ 
with right hand side as in Lemma \ref{lemma:areadisk}, for
$p$ being the constant path at $L$. Suppose that $\Theta _{-}$  is regular. 
Define $PSS ([L]) $ as the homology class of the chain $C
([L]) \in CF(\mathcal{A} _{0} ) $ determined by:
\begin{equation} \label{eq:CL}
 \langle C ([L] ) , \gamma  \rangle = \sum _{A} \#  ev (\Theta _{-}, A),
\end{equation}
where the sum is over all classes $A \in H _{2} ^{sec}
(\mathcal{D} \times M, \mathcal{L}), $ which have asymptotic
constraint $\gamma$, and where $\gamma$ is a geometric
generator of $CF( \mathcal{A} _{0} )$. 

Now, for a general class $a \in QH (L) $, $PSS (a) $ is
defined similarly, but using the moduli space $\mathcal{M}
(\Theta _{-}, a, A) $. The latter can be defined as the subset  of $\mathcal{M} (\Theta _{-}, A) $ consisting of sections
intersecting a fixed smooth pseudocycle, see Zinger
~\cite{citeZingerPseudocycles}, representative of $a$.
More specifically, for $z _{0} \in \partial \mathcal{D}
$ let $\widetilde{S} _{z _{0}}$ be the fiber. Fix
a pseudo-cycle $g: B \to L \subset \widetilde{S}  _{z
_{0}}$ representing $a \in H _{2} (L) $. Then $\mathcal{M}
(\Theta _{-}, a, A) $ consists of elements of $\mathcal{M}
(\Theta _{-}, A) $ intersecting image of $g$. (Although we
use the language of pseudocycles in this outline, for analysis it is technically simpler to use Morse homology and Perl complex language as in ~\cite{citeAlbersPSS}.) 

Now, the PSS morphism is an isomorphism in our monotone context,
and $CF(\mathcal{A} _{i} )$ is perfect for each $i$, by
assumption.  It follows
that the asymptotic constraint $\gamma _{i} $ of $\sigma
_{0} $ at each (positive) end $e _{i} $ satisfies:
\begin{equation*}
\langle \gamma _{i}, PSS (a)  \rangle =1,
\end{equation*}
for some $a$ uniquely determined. 
Moreover, Fredholm index and monotonicity restrictions
insure that only a single class $A^i$ can contribute in
the formula for $PSS (a) $ (analogous to \eqref{eq:CL}). 
Let then $\sigma _{A ^{i}} \in \mathcal{M} (\Theta _{-}, a, A _{i}) $ be some element.

Note that at the negative ends, the above story needs to be
suitably dualized, but we will not elaborate, as this is all
very standard. With this understanding, at each end $e _{i} $, glue $\sigma
_{0}$ with $\sigma _{A ^{i}}$. We then obtain a $J (\mathcal{A} ^{/} _{r _{0} } )$-holomorphic, class $A ^{/} $ section $\sigma ^{/} _{0}$ of $\Theta ^{/} _{r _{0}}  $. 

By Lemma \ref{lemma:smallempty}:
\begin{equation*}
   carea (\Theta _{r _{0} } ^{/}, A ^{/} ) \leq \area (\mathcal{A} _{r_0} ^{/}) \leq \area(\mathcal{A} _{r _{0} })  + \sum _{i} \kappa  _{i},  
\end{equation*}
so 
\begin{equation*}
   carea (\Theta _{r _{0} } ^{/}, A ^{/} ) - \sum _{i} \kappa  _{i}  \leq \area (\mathcal{A} _{r _{0} }),
\end{equation*}
so that we contradict the hypothesis. So in the case $\overline{\mathcal{M}} (\{\Theta _{r}\} , A)$ is regular we are done with the first part of the lemma. When it is not regular instead of gluing just pre-glue to get a holomorphic building $\sigma ^{/} _{0}$, 
and the conclusion follows by the same argument.

To prove the last part of the lemma, note that each $\Theta ^{/} _{r} $ is taut concordant to $$\Theta _{0}= (D ^{2} \times M, D ^{2}, \mathcal{L}, \mathcal{A} ^{tr} ),  $$ with $\mathcal{L}$ trivial with fiber $L _{0} $,  and for  $\mathcal{A} ^{tr} $ the trivial connection. And $$carea (\Theta _{0}, \cdot) = {-}{const}\cdot Maslov ^{vert}  (\cdot)$$  as functionals on $H _{2} ^{sec}  (D ^{2} \times M, \mathcal{L}) $. It follows by Lemma \ref{lemmaInvariantPairing} that
$$carea (\Theta ^{/}, \sigma ^{/}_{0} )= carea (\Theta _{0}, \sigma ^{/} _{0} ) = -{const}\cdot Maslov ^{vert}  (\sigma ^{/} _{0}  ) =  -{const}\cdot Maslov (A ^{/}). $$

\end{proof}

\section {Construction of small data}
\label{sec:constructionSmallData}
To forewarn, we use here notation and notions from Part I. 

Let $(m ^{1}, \ldots, m ^{d}) $  be a composable chain 
of morphisms in $\Pi (\Delta ^{n}) $, which we recall means
that the target of $m ^{i-1} $ is the source of $m ^{i}
$ for each $i$. 
The perturbation data $\mathcal{D} $, in
particular, specifies for each $n$ and for each such
composable chain, certain maps $$ 
u (m ^{1}, \ldots, m ^{d},n): \mathcal{E} _{d} ^{\circ}   \to \Delta ^{n}, 
$$
where $\mathcal{E} _{d}  $ is the universal curve over
$\overline{\mathcal{R}} _{d}$, and $\mathcal{E} _{d}
^{\circ}  $ denotes $\mathcal{E} _{d}  $ with nodal points
of the fibers removed. The collection of these maps,
satisfying certain axioms, is denoted  by $\mathcal{U}$. We have already mentioned this in the introduction.

The restriction of $u (m ^{1}, \ldots, m ^{d},n)$ to the fiber
$\mathcal{S} _{r} $ of $\mathcal{E} _{r} ^{\circ}  $ over
$r$, is denoted by $u (m ^{1}, \ldots, m ^{d},n, r)$ which may also be abbreviated by $u _{r} $.

Let $\Sigma: \Delta ^{n} \to X$ be smooth, denote $\overline {m} ^{i}
:= \Sigma \circ m ^{i}  $, and denote
$\widetilde{\mathcal{S}} _{r}:= (\Sigma \circ u _{r})
^{*}P$.
Then  $\mathcal{D}$ specifies for
each such $\Sigma $, each $r$, each composable chain $(m^1, \ldots, m ^{d}) $ in $\Pi (\Delta ^{n} )$,  and for each chain of objects $L' _{0}, \ldots, L' _{d}  $ with $L' _{i} \subset P| _{\overline{m} _{i} (1)} $, $i \geq 1$, $L' _{0} \subset P | _{\overline{m} _{1} (0)} $,
a certain Hamiltonian connection $\mathcal{A} _{r}$ on 
\begin{equation} 
 \widetilde{\mathcal{S}} _{r} \to \mathcal{S} _{r}.
\end{equation} 
Using terminology of the previous section, we can say that
$\mathcal{D} $ specifies a Hamiltonian structure
$(\widetilde{\mathcal{S} } _{r}, \mathcal{S} _{r},
\mathcal{L} _{r}, \mathcal{A} _{r}) $ where $\mathcal{L}
_{r}$ is trivial over each boundary component, with fiber
the corresponding objects $L' _{i}$.

Now, specialize to the case $n=4$,  $\Sigma = \Sigma
_{+} $, $\forall i: L' _{i} = L _{0}  $. In this case, we
write $\mathcal{A} _{r} ^{+} (m ^{1}, \ldots, m ^{d} ) $ for the corresponding connections,
further abbreviated by $\mathcal{A} _{r} ^{+}$ as $(m ^{1}, \ldots, m ^{d} )$ will usually be implicit. 

Suppose that $\mathcal{D}$ extends $\mathcal{D} _{pt} $ from
before. If $\mathcal{L} _{r} \subset \widetilde{\mathcal{S}}
_{r}| _{\partial \mathcal{S} _{r} }  $ denotes the trivial
Lagrangian sub-bundle with fiber $L _{0} $, then we obtain
a Hamiltonian structure $\Theta ^{+}= \{\Theta ^{+} _{r}
\}= \{\widetilde{\mathcal{S}} _{r}, \mathcal{S} _{r},
\mathcal{L} _{r}, \mathcal{A} ^{+}   _{r}    \}
_{\overline{\mathcal{R}} _{d}  } $. By the properties of
these connections, necessitated by $\mathcal{D}$, at each
end $e _{i} $ of $\mathcal{S} _{r} $, $\mathcal{A} ^{+} _{r}
$ is compatible with the connection $\mathcal{A} _{i}
=  {\mathcal{A}} (L _{0}, L _{0})$,  where $\mathcal{A} (L
_{0}, L _{0})$ is the connection on $[0,1] \times S ^{2}
$ also part of our data $\mathcal{D}$.
Then $\Theta ^{+} $ is trivially taut since for each $r$
$\mathcal{L} _{r} $ is naturally trivial and $\mathcal{A}
_{r} $ is likewise trivial over $\partial \mathcal{S} _{r}
$, for each $r$, by the assumed properties of these connections.

%
Set $$\hbar:= \frac{1}{2} \area (S ^{2},\omega). $$
Let $\kappa$ denote the $L ^{\pm} $-length of the holonomy
path in $\operatorname {Ham} (S ^{2}, \omega) $ of $\mathcal{A} _{0}= \mathcal{A} (L _{0}, L _{0})$. We may suppose that  
\begin{equation} \label{eqAreaAd}
\forall r: \area (\mathcal{A} ^{+}  _{r}  ) <  \hbar- 5 \kappa,
\end{equation}
is satisfied after taking $\kappa$ to be sufficiently small.
(There is no obstruction since the corresponding bundles are
naturally trivializeable, continuously in $r$.) 

%
Fix a complex structure $j _{0} $ on $M$,
and let $\{J _{r} = J (\mathcal{A} ^{+} _{r})  \}$ be the
corresponding induced family of almost complex structures on
$\{\widetilde{\mathcal{S}} _{r}  \}$ as in Section
\ref{sec:ModuliSpacesHamStructures}.
\begin{lemma}
\label{lemma:empty0} 
As in Part I,  let $$ 
   \overline {\mathcal{M}}=\overline{{\mathcal{M}}} (\gamma ^{1}, \ldots, \gamma ^{d}; \gamma ^{0}, \Sigma _{+},
\{J_{r}   \},
   A), 
$$ denote the set of elements of $\overline{{\mathcal{M}}} (\Theta ^{+} 
,A)$ with asymptotic constraints $\gamma ^{i}$ at each $e
_{i} $ end. Here each $\gamma ^{k}   $, $k \neq 0$, is of
the form $\gamma _{i,j} $ where this is as in Section \ref{section:model}. Then whenever the class $A$ is such that ${\mathcal{M}}$ has virtual dimension $0$, and $d$ satisfies $2 < d \leq 4$, $\overline {\mathcal{M}}$ is empty.
\end{lemma}
\begin{proof}
   Let $$\Theta ^{/} := (\Theta ^{+}) ^{/},   $$ and $A ^{/} $ be as in Section \ref{section:gluing}. For a fixed $r$, by the Riemann-Roch (Appendix \ref{appendix:maslov}) we get that the expected dimension of ${\mathcal{M}} (\Theta ^{/}, A ^{/} ) $ 
is 
\begin{equation*}  1+ Maslov ^{vert}  (A ^{/} ).
\end{equation*} 
Consequently, when $\gamma ^{0} =\gamma $, the expected dimension of $\mathcal{M}$ is:
\begin{equation} \label{eqDimension} 1+ Maslov ^{vert}  (A ^{/})  - 1   + (\dim
   \mathcal{R}_{d} = d-2).
\end{equation} 
We need the expected dimension of $\mathcal{M}$ to be 0, and $d \geq 3$, so $Maslov ^{vert} (A ^{/} ) \leq -1$. But $Maslov ^{vert} (A ^{/} )
   = -1$ is impossible as the minimal positive Maslov number
	 is 2.   

Now, note that if $Maslov ^{vert} (A ^{/} ) = -2$ then $-C \cdot Maslov ^{vert} (A ^{/} )
   = \hbar $, for $C$ the monotonicity constant of $(S ^{2},
	 \omega ) $ and the equator $L _{0}$.  Consequently, the result follows by Lemma \ref{lemma:gluinglowerbound} and by the property \eqref{eqAreaAd}.

%
When $\gamma ^{0}  $ is the Poincare dual to $\gamma $, we would get $Maslov ^{vert} (A ^{/} )
   \leq -2$  so for the same reason the conclusion follows. 
\end{proof}

So if we choose our data $\mathcal{D}$ so that the hypothesis of the lemma above are satisfied,
then with respect to this $\mathcal{D}$:
\begin{align} 
   \label{eq:mu:identities2} \mu ^{2} _{\Sigma ^{4} _{\pm}    }
(\gamma _{i,j} , \gamma _{j,k} ) &  = \gamma _{i,k}  \\
   \label{eq:mu:identities1}
   \mu ^{3} _{\Sigma ^{4} _{\pm} }  (\gamma ^{1}, \ldots, \gamma ^{3}) & = 0,    \text{ for $\gamma ^{i} $ as above} \\
\mu ^{4} _{\Sigma ^{4} _{+}}  (\gamma _{1}, \ldots, \gamma _{4}) & = 0.
\end{align}
In particular this $\mathcal{D} $ is small.

\section {The product $\mu ^{4} _{\Sigma ^{4}
_{-}} (\gamma _{1} , \ldots , \gamma _{4} ) $ and
the quantum Maslov classes} \label {section:productAndHigherSeidel}

The product $$\mu ^{4} _{\Sigma ^{4} _{-}} (\gamma _{1}
, \ldots , \gamma _{4} ) $$ a priori depends on various
choices, like the choices of $h _{\pm} $, and then choice of
data $\mathcal{D} _{0} $. However by Lemma
\ref{lemma:gluing}, so long as there is a homotopy of the
choices, together with a homotopy of associated perturbation
data $\{\mathcal{D} _{t} \}$, so that $\mathcal{D} _{t} $ is
small for all $t$, the above product is $t$-invariant. In
particular, for the purpose of computation we may take 
$h _{+} $ to be the constant map to $x _{0} $ and $$h _{-}:
(D ^{4}, \partial D ^{4}) \to (S ^{4}, x _{0})    $$ to be the complementary map, that is representing the generator of $\pi _{4} (S ^{4}, x _{0}) \simeq \mathbb{Z}$. We further suppose that $h _{-} $   is an embedding in the interior of $D ^{4} $.

Let $\Sigma _{-} $ be the 4-simplex of $S ^{4} _{\bullet}
$ corresponding to $h _{-}$  as before in Section \ref{section:model}. We need to study the moduli spaces
\begin{equation} \label{eqModuli}
\overline{\mathcal{\mathcal{M}}} ( \gamma _{1}, \ldots,  \gamma _{4}; \gamma
   ^{0}, \Sigma_{-}, \{\mathcal{A}_{r} \}, A),
\end{equation}
where $\mathcal{A} _{r}$ now denotes the connections on
\begin{equation} 
 \widetilde{\mathcal{S}}   _{r}:= (\Sigma _{-}  \circ u (m _{1}, \ldots,  m _{4}, 4,r)) ^{*}P \to \mathcal{S} _{r},
\end{equation}
part of some small data $\mathcal{D} _{0}  $ as above. We abbreviate $u (m _{1}, \ldots,  m _{4}, 4,r)$ by $u _{r} $ in what follows.

By the dimension formula \eqref{eqDimension}, since we need the expected
dimension of \eqref{eqModuli} to be zero, the class $A ^{/} $ satisfies:
\begin{equation*}
Maslov ^{vert} (A ^{/} ) = -2,
\end{equation*}
and we must have $\gamma ^{0} = \gamma _{0,4}.  $ 

\begin{notation}
 From now on, by slight abuse, $A _{0} $ refers to various section classes of various Hamiltonian structures such that the associated class $A _{0} ^{/}  $ satisfies:
   \begin{equation*}
   Maslov ^{vert} (A _{0} ^{/}  ) = -2.
   \end{equation*}
\end{notation}
\subsection {Constructing suitable $\{\mathcal{A} _{r} \} $}
To get a handle on \eqref{eqModuli} we want to construct
very special, small data $\mathcal{D} _{0} $.

A Hamiltonian $S ^{2} $ fibration over $S ^{4} $ is classified by an element $$[g] \in \pi _{3} (\operatorname {Ham} (S ^{2}, \omega), id) \simeq \pi _{3} (PU (2),id) \simeq \mathbb{Z}.$$ Such an element determines a fibration $P _{g} $
over $S ^{4} $ via the clutching construction:
\begin{equation*}
P _{g} =   D ^{4} _{-} \times S ^{2}    \sqcup  D ^{4} _{+} \times  S ^{2} \sim,   
\end{equation*}
with $D ^{4}_- $, $D ^{4}_+ $ being 2 different names for the standard closed 4-ball $D
^{4} $, and where the equivalence relation $\sim$ is $(d, x) \sim \widetilde{g}
(d,x)$, $$ \widetilde{g}: \partial D ^{4} _{-}
\times S ^{2} \to \partial D ^{4} _{+}  \times S
^{2},  \quad   \widetilde{g} (d,x) = (d, g (d)
^{-1}  (x)).$$  We suppose that the that previously appeared
point $x _{0} = h _{\pm} (b _{0}) $, is in $D ^{4} _{+} \cap
D ^{4} _{-} \subset S ^{4} $.

From now on $P _{g} $ will denote such a
fibration for a non-trivial class $[g]$. Note that the fiber of $P _{g} $ over the base point $x _{0}  \in S ^{3} \subset D ^{4} _{\pm} $ (chosen for definition of the homotopy group $\pi _{3} (\operatorname {Ham} (S ^{2}, \omega), id) $)  has a distinguished, by the construction, identification with $S ^{2} $. 
Take $\mathcal{A}$ to be a connection on $P \simeq P _{g} $ which is trivial in the distinguished trivialization over $D ^{4} _{+}  $. 
This gives connections $$\mathcal{A}' _{r}:= (\widetilde{u}  _{r}) ^{*} \mathcal{A} $$ on $\widetilde{\mathcal{S}} _{r}$, $$\widetilde{u}_{r}  = \Sigma _{-} \circ u _{r}. $$

By the last axiom for the system $\mathcal{U}$ introduced in
Part I, we may choose $\{u _{r} \} $ so that the family
$\{\widetilde{u} _{r} (\mathcal{S} _{r}) \}$ induces
a singular foliation of $S ^{4} $ with the properties:
\begin{itemize}
	\item  The folliation is smooth outside $x _{0} $. Note
	that $x _{0} $ is the image by $\widetilde{u} _{r}  $ of the
	ends (images of $e _{i} $), and the image of the boundary of each $\mathcal{S} _{r} $.
	\item Each $\widetilde{u} _{r} $ is an embedding on the complement of $\widetilde{u} _{r} ^{-1} (x _{0} ) $. 
\end{itemize}
Denote by $E$ the subset $S ^{3} \subset S ^{4} $ bounding $D _{\pm} ^{4}  $.  
We may in addition suppose that each $\widetilde{u} _{r} $ intersects $E$ transversally, again on the complement of $\widetilde{u} _{r} ^{-1} (x _{0} ) $.

\begin{figure} [h]
 \includegraphics[width=2in]{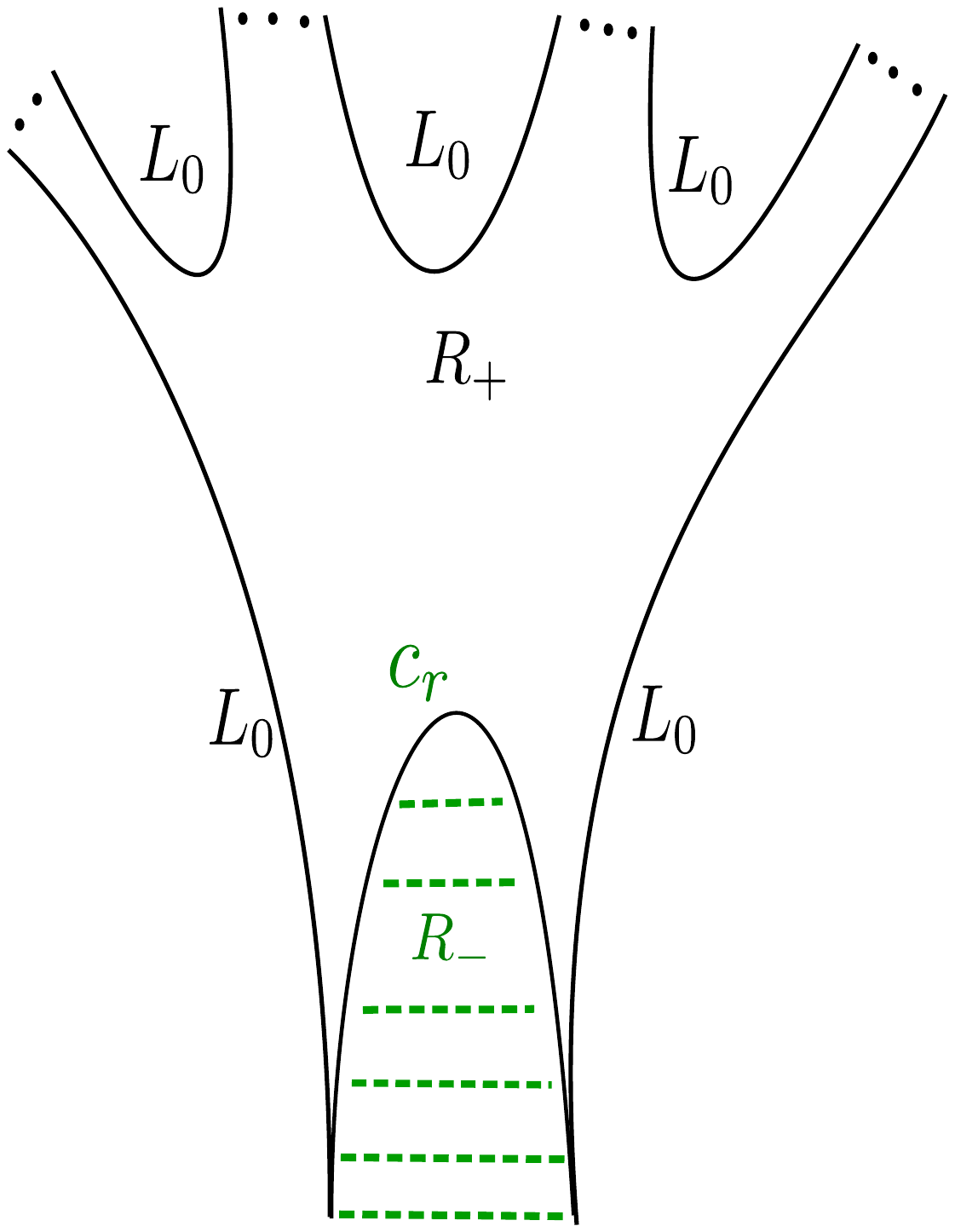}
   \caption {The labels $L _{0} $ indicate that the Lagrangian subbundle is constant with corresponding fiber $L _{0} $. The curve $c _{r} $ bounds $R _{-}$.}
\label{figure:plusminus}
\end{figure}

By the above, the preimage by $\widetilde{u}  _{r} $ of $E$  contains a smoothly embedded curve $c _{r} $ as in Figure \ref{figure:plusminus}, 
and $\widetilde{u} _{r} $ takes $c _{r} $ into $E$. 
This $c _{r} $ not uniquely determined, but we may fix
a family $r \mapsto c _{r} $, with parametrizations $$c
_{r}: \mathbb{R} \to \mathcal{S} _{r},  $$ with the
properties:
\begin{itemize}
	\item $c _{r}$ maps $(-\infty,0)  $ diffeomorphically onto
	$e _{0} (\{0\} \times (-\infty, 0))$.
	\item $c _{r}$ maps $(1,\infty)  $ diffeomorphically onto
	$e _{0}(\{1\} \times (-\infty, 0))$.
	\item  $\{c _{r} \} $ is a $C ^{0} $ continuous family in
	$r$.
\end{itemize}  
 We set:
\begin{equation*}
\widetilde{c}_{r}:= \widetilde{u} _{r} \circ c _{r}.
\end{equation*}
In Figure \ref{figure:plusminus}, the regions $R _{\pm} $ are the preimages by $\widetilde{u}  _{r} $ of $D^{4} _{\pm} \subset S ^{4} $, and $c _{r} $ bounds $R _{-} $. 
It follows that  $\{\widetilde{c}  _{r} \}$ likewise induces a singular foliation of the equator $E \simeq S ^{3} $ that is smooth outside $x _{0} $.

So each $\mathcal{A}' _{r} $ is flat in the region $R _{+} $, in fact is trivial in the distinguished trivialization of $\widetilde{\mathcal{S}} _{r}  $ over $R _{+} $, corresponding to the distinguished trivialization 
of $P$ over $D ^{4} _{+}  $. Likewise we have a distinguished trivialization of $\widetilde{\mathcal{S}} _{r}  $ over $R _{-} $,
corresponding to the distinguished trivialization of $P$ over $D^{4}_{-}$. 
In this latter trivialization let $$\phi _{r}: \mathbb{R} \to \operatorname {Ham} (S ^{2}, \omega) $$ be the holonomy path of $\mathcal{A} _{r}' $ over $c _{r} $. 
Then by construction, $$\phi _{r}| _{(-\infty,0] \sqcup [1, \infty)}  = id, $$ 
so that we may define $$f ({r}) \in \Omega _{L _{0}}  Lag (M,L _{0} ),$$
by
$$ f ({r}) (t)= \phi _{r} (t) (L _{0}), \quad t \in [0,1]$$
where the right hand side means apply an element of $\operatorname {Ham} (S ^{2}, \omega) $ to $L _{0} $ to get a new Lagrangian. We will say that $f (r)$ is \textbf{\emph{generated by $\mathcal{A}' _{r} $}}.

Note that by construction 
\begin{equation} \label{eq:phi}
\phi _{r} (t) =  g(\widetilde{c}  _{r}(t)),  \quad t \in [0,1],
\end{equation}
if we identify $\widetilde{c} _{r}(t) $ with an element of $S ^{3} $.

Let $D ^{2}_{0} \subset \overline{\mathcal{R}}_{4}$ be an
embedded closed disk, not intersecting the boundary
$\partial \overline{\mathcal{R}}_{4}$, so that $\partial
D ^{2} _{0}  $ is in the gluing normal neighborhood $N$ of
$\partial \overline{\mathcal{R}}_{4}  $, as defined in Part I.

So we have a continuous map $$f: D ^{2} _{0} \to \Omega _{L _{0}} Lag (S ^{2}). $$
And $f(\partial D ^{2} _{0})= L _{0} $, with the right hand side denoting the constant loop at $L _{0} $.
Then by construction, and \eqref{eq:phi} in particular, $f \simeq lag$, where $\simeq$ is a homotopy equivalence, and where
\begin{equation} \label{eq:lagdef}
lag: S ^{2} \to \Omega _{L _{0} } Lag (S ^{2})
\end{equation}
is the composition 
$$S ^{2} \xrightarrow{g'} \Omega _{id} PU (2) \to \Omega _{L_{0} }Lag ^{eq}  (S ^{2}), $$ for $g'$ naturally induced by $g$, and for the second map naturally induced by the map $$PU (2) \to Lag ^{eq} (S ^{2} ), \quad \phi \mapsto \phi (L _{0}).$$

We then deform each $\mathcal{A}'_{r} $ to a connection $\mathcal{A} _{r} $, which is as follows. In the region $R _{+} $ $\mathcal{A} _{r} $ is still flat, but at each end $e _{i} $, $\mathcal{A} _{r} $ is compatible with $\mathcal{A} (L _{0}, L _{0})$, where this is as in Section \ref{section:exact}, and so that $\mathcal{A} _{r} $ is still trivial over the boundary of $\mathcal{S} _{r} $.


Since $\widetilde{\mathcal{S}} _{r}  $ and $\mathcal{A}' _{r} $ are trivial  for $r \in \overline{\mathcal{R}} _{d} - D ^{2} _{0} $, with trivialization induced by the trivialization of $P _{+} $, and since the condition \eqref{eqAreaAd} holds, we may insure that  
\begin{equation} \label{eq:hbarbound2}
   \area (\mathcal{A} _{r}) < \hbar -5\kappa, 
\end{equation}
for $r$ in the complement of $D ^{2} _{0}  $.
In other words $\{\mathcal{A} _{r}\} $ extends to a system of connections corresponding to small data $\mathcal{D} _{0} $ for $P$, as intended. 
\subsection{Restructuring $\{\mathcal{A} _{r} \}$}
Applying Lemma \ref{lemma:gluinglowerbound} we see that the
resulting Hamiltonian structure $\mathcal{H} := \{\mathcal{H} _{r}:= \widetilde{\mathcal{S}} _{r}, \mathcal{S} _{r}, \mathcal{L} _{r}, \mathcal{A} _{r} \}$ is $A _{0} $-admissible. We now further mold this data for the purposes of computation.

First cap off the ends $e _{i} $, $i \neq 0$, of
each $\mathcal{H} _{r}$ as in the paragraph preceding
Lemma \ref{lemma:gluinglowerbound}. 
This gives a Hamiltonian structure $$\mathcal{H} ^{\wedge}:= \{ \widetilde{S} _{r} ^{\wedge},  {S} _{r} ^{\wedge},  \mathcal{L} _{r} ^{\wedge}, \mathcal{A} _{r}^{\wedge} \} _{\mathcal{K} = D ^{2} _{0}} ,$$ satisfying $$\area (\mathcal{A} _{r} ^{\wedge}  ) + \kappa < \hbar,$$ for each $r$. Again by Lemma \ref{lemma:gluinglowerbound} $\mathcal{H} ^{\wedge}$ is $A _{0} $-admissible. 
%

By the classical gluing of holomorphic curves it follows that 
\begin{equation} \label{eq:mu=H^}
[\mu ^{4}_{\Sigma ^{4}} (\gamma _{1} , \ldots , \gamma _{4} )] = [ev (\mathcal{H} ^{\wedge}, A _{0} )].
\end{equation}
It remains to compute the right hand side, to this end we further restructure. 

\begin{figure} 
 \centering 
\includegraphics[width=2in]{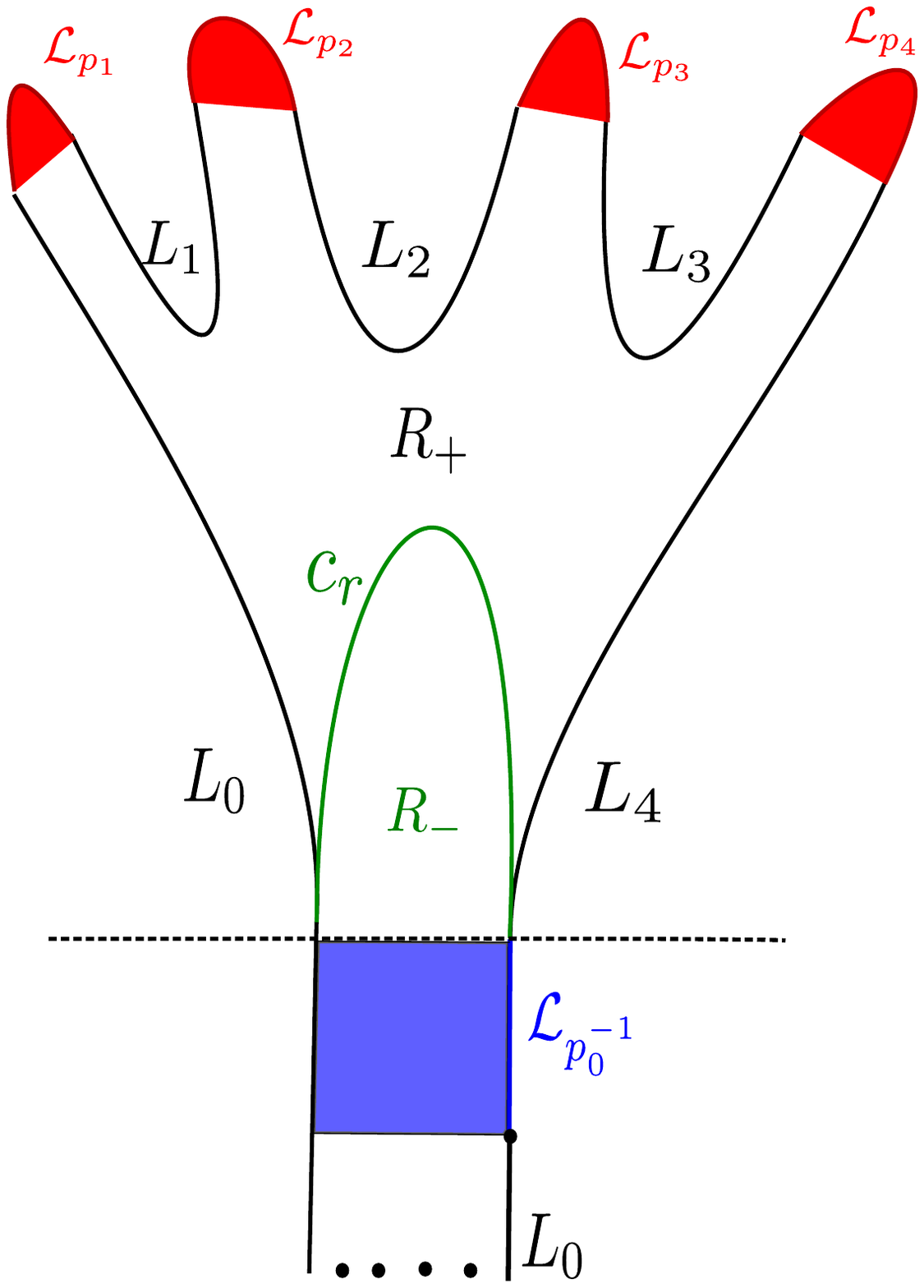}
   \caption {Over the boundary components with black labels
	 $L _{i} $ the  Lagrangian subbundle $\mathcal{L} _{r}
	 ^{\mathfrak {n}} $ is constant with corresponding fiber
	 $L _{i} $. Over the $i$'th red boundary component the
	 Lagrangian subbundle corresponds to the path of
	 Lagrangians $p _{i} $, analogously to Definition
	 \ref{defInducedLoop} further below. Likewise over the
	 right boundary component of the blue region the Lagrangian subbundle
	 corresponds to the path of Lagrangians $p _{0} ^{-1}
	 $. In the red striped regions we have removed the curvature of the connection, the blue striped region we have added it.}
\label{diagram:deformation}
\end{figure}

Let $p _{1}: [0,1] \to Lag (S ^{2}, L _{0}  ) $ be the path generated by $\mathcal{A} (L _{0}, L _{0})$, with $p _{1} $ starting at $L _{0} $, and where generated has the same meaning as in the previous section.
Suppose we have defined $p _{i-1} $, set $L _{i-1}:=p _{i-1} (1)  $ and define $p _{i} $ to be the path in $Lag (S ^{2}, L _{0}  )$ starting at $L _{i-1} $, generated by $\mathcal{A}(L _{0}, L _{0})    $.
Set $p _{0} :=p _{1} \cdot \ldots \cdot {p _{d}}  $, where $\cdot$ is path concatenation in diagrammatic order. 
We may assume that $L _{0} $ is transverse to $L _{4} = p _{0} (1)$ by adjusting the connection $\mathcal{A} (L _{0}, L _{0})$ if necessary. 
Then for each $r$, deform $\mathcal{L} _{r} ^{\wedge}  $ to
the Lagrangian subbundle over $\partial \mathcal{S} _{r}
$ denoted $\mathcal{L} _{r} ^{\mathfrak {n}}  $, which is as
illustrated in Figure \ref{diagram:deformation}.

We simultaneously deform $\mathcal{A} _{r} ^{\wedge}  $ to
an $\mathcal{L} _{r} ^{\mathfrak {n}}  $ exact Hamiltonian
connection $\mathcal{A} _{r} ^{\mathfrak {n}}  $ which
satisfies the following conditions, (referring  to the Figure
\ref{diagram:deformation}):
\begin{itemize}
	\item $\mathcal{A} _{r} ^{\mathfrak {n}} $ is flat in the entire region $R _{+} $ (which includes the red shaded finger regions).
	\item The blue region is contained in the strip end chart
	at the $e _{0} $ end, which is down in the figure.
	\item Along top boundary segment of the blue region,
	contained in the dashed line,  $\mathcal{A} _{r}
	^{\mathfrak {n}}  $ is the trivial connection in the
	distinguished trivialization at the $e _{0}$ end.
	\item At the $e _{0} $ end the connection is unchanged over $[0,1] \times (s, \infty)$, for $s$ large.
\end{itemize}
In order to get such a deformation, we have to introduce curvature in the blue region of Figure \ref{diagram:deformation}.

We name this new Hamiltonian structure by $$\mathcal{H} ^{\mathfrak {n}}:=\{\widetilde{S} _{r} ^{\wedge}, {S}_{r}^{\wedge}, \mathcal{L}_{r} ^{\mathfrak {n}}, \mathcal{A} _{r}^{\mathfrak {n}}\}.$$ 
It can be understood as obtained by capping off the
respective ends $e _{i}$, $i \neq 0$ of $\{\mathcal{H}
_{r}\}$  with $\Theta ({H _{i}})$, where the latter is as in
\eqref{eq:ThetaH}. Note however, that the capping is
modified (in the blue region). 

For each $p _{i} $, $L ^{\pm} (p _{i})$  can be arranged to
be arbitrarily small. Consequently, we may choose the
deformation from $\mathcal{H} ^{\wedge} $ to $\mathcal{H}
^{\mathfrak {n}} $ to be small near boundary (Definition
\ref{def:small}) and hence to be an $A _{0} $-admissible concordance. More specifically, we may choose a concordance 
from $\mathcal{H} ^{\wedge} $ to $\mathcal{H} ^{\mathfrak
{n}} $ such that for the associated family of connections
$\{\mathcal{A} _{r,t} \}$, $t \in [0,1] $, the following is
satisfied:
\begin{itemize}
	\item $$\mathcal{A} _{r,0}= \mathcal{A} _{r} ^{\wedge},
	\quad \mathcal{A} _{r,1}= \mathcal{A} _{r} ^{\mathfrak {n}}.  $$
	\item 
	\begin{equation*} \forall r \, \forall v \in T _{z} \mathcal{S} _{r}: \frac{d}{dt} |R _{\mathcal{A} _{r,t}} (v, jv)| _{+} <0, 
\end{equation*}
for each $z \in \mathcal{S} _{r} $, except for $z$ in the
region which is blue shaded in Figure \ref{diagram:deformation}. 
	\item  
$$\forall t: |\area (\mathcal{A} _{r,t} )-\area (\mathcal{A} _{r} ^{\wedge})| \leq L ^{+} (p _{0} ^{-1}). $$ 
\end{itemize}  
The last condition can be satisfied as the area increase in
this blue region is bounded from above by $L ^{+} (p _{0}
^{-1})$. In fact we can arrange that $$ \frac{d}{dt}\area (\mathcal{A} _{r,t})=0,$$ since the gain of $\area$ in the blue region is exactly equal to the loss of $\area$ in the red regions, but this extra precision is not necessary.

Of course:
\begin{equation*}
   [ev (\mathcal{H} ^{\wedge}, A _{0} )] =
   [ev(\mathcal{H} ^{\mathfrak {n}}, A _{0} )]
\end{equation*}
since the corresponding Hamiltonian structures are $A _{0} $-admissibly concordant.
This finishes our restructuring.
\subsection {Computing $[ev(\mathcal{H} ^{\mathfrak {n}}, A _{0} )]$}
If we stretch the neck along the  
dashed line in Figure \ref{diagram:deformation}, the upper
half of the resulting building gives us a new Hamiltonian
structure:
$$\mathcal{H} ^{0} = \{\widetilde{S} _{r} ^{0}, {S}_{r}^{0}, \mathcal{L}_{r} ^{0}, \mathcal{A} _{r}^{0}\}.$$
By the classical theory of continuation maps in Floer homology we clearly have that 
$$[ev _{\mathfrak {n}}] :=[ev(\mathcal{H} ^{\mathfrak {n}}, A _{0} )] \in FH (L _{0}, L _{0})$$ 
is non-zero iff 
$$ 
[ev _{0} ] := [ev(\mathcal{H} ^{0}, A _{0}]  \in FH (L _{0}, L _{4})
$$ is non-zero.

Let $\mathcal{P}_{L _{0}, L _{4}}  Lag (M)$ denote the space of smooth paths in $Lag (M)$ from $L _{0} $ to $L _{4} $. Let $$f': D _{0} ^{2} \to  \mathcal{P}_{L _{0}, L _{4}}  Lag (S ^{2}), $$ be like $f$ but defined with respect to $\{\mathcal{A} _{r} ^{0} \} $, so that 
\begin{equation} \label{eq:f'}
f' (t)= g(\widetilde{c} _{r}(t))(p _{0}  (t)),
\end{equation}
if we suppose that the holonomy path of $\mathcal{A} ^{0}
_{r} $ over $c _{r}$, in the trivialization over $R _{+} $, generates $p _{0} $ (which can be insured by adjusting $\{\mathcal{A} ^{0} _{r}\}$ or the parametrizations $\{c _{r}\}$).
Here the right hand side of \eqref{eq:f'} means as before: apply an element of $\operatorname {Ham} (S ^{2}, \omega) $ to a Lagrangian to get a new Lagrangian.
In this case $$f'(\partial D _{0} ^{2}) =   p_0 \in
\mathcal{P}_{L _{0}, L _{4}}  Lag (M).$$
In particular, $f'$ represents a class $a \in \pi _{2} (\mathcal{P}_{L _{0}, L _{4}}  Lag (M), p _{0}) $.

In what follows we omit specifying the parameter space $D ^{2} _{0}  $ for $r$, since it will be the same everywhere.
%
Let  $\mathcal{L} _{p}$ be as in the Definition \ref{defInducedLoop}.  
\begin{lemma} \label{lemma:concordant} The $A _{0} $-admissible  Hamiltonian structure $\mathcal{H} _{0} = \{\widetilde{\mathcal{S}}_{r}^{0}, \mathcal{S} _{r} ^{0},  \mathcal{L}
_{r}^{0},
\mathcal{A} ^{0}_{r}    \} $ is $A _{0} $-admissibly concordant to 
$$ \Theta' = \{
   \mathcal{D} \times S ^{2}, \mathcal{D}, \mathcal{L} _{f'} (r), \mathcal{B} _{r} 
\},$$ for certain Hamiltonian connections $\{\mathcal{B} _{r} \}$ (which are not explicitly relevant yet).
\end{lemma}

\begin{proof}
Let $R _{\pm} \subset \mathcal{S}_{r} ^{\wedge}  $ be as before. 
Fix a smooth deformation retraction $$ret _{r}: \mathcal{S} _{r} ^{0} \times I \to
\mathcal{S}_{r}^{0}, $$  of $\mathcal{S} _{r} ^{0} $ onto $R ^{-}  $, smooth in $r$. 
Since $\mathcal{A} _{r} ^{0} $ is flat over $R ^{+} $, 
the pull-back by $ret _{r} $ of the data $\mathcal{H} _{0} $ then  induces an  $A _{0} $-admissible concordance 
between $\mathcal{H} _{0} $ and $$\{
   \mathcal{D} \times S ^{2}, \mathcal{D}, \mathcal{L}_{f'}  (r), \mathcal{B} _{r} = ret _{r} ^{*} \mathcal{A} _{r} ^{0}    
   \},$$ 
once we use smooth Riemann mapping theorem to identify
each $R ^{-} \subset \mathcal{S} _{r} ^{0}   $ with its induced complex structure $j _{r}   $ with $(\mathcal{D}, j _{st}) $,
   smoothly in $r$. 
\end{proof}

\section {Quantum Maslov classes} 
\label{section:QuantumMaslov}
The class of $ev(\Theta', A _{0})$ is related to
what we christen as quantum Maslov classes.  These
are relative analogues of the quantum 
characteristic classes
\cite{citeSavelyevQuantumcharacteristicclassesandtheHofermetric}.
The name quantum Maslov
class is meant to be suggestive,  as the classical Maslov numbers are
relative analogues of Chern numbers, while 
quantum characteristic classes are
directly related (via semi-classical
approximation) to Chern classes,
~\cite{citeSavelyevBottperiodicityandstablequantumclasses}.
\footnote{The latter may be renamed
   ``quantum Chern classes''.}

We will not give extensive detail here since we
don't need the full theory, we present it
because it gives extra perspective.
%
%
The ordinary relative Seidel morphism appears in Seidel's
\cite{citeSeidelFukayacategoriesandPicard-Lefschetztheory} in the exact case
and further developed in
\cite{citeHuLalondeArelativeSeidelmorphismandtheAlbersmap} in the monotone case.  Let $Lag (M)$ denote the space whose components are objects of $Fuk (M)$ in the previous sense, so in particular oriented, spin, Hamiltonian isotopic Lagrangian submanifolds of $M$. We may also denote the component of $L$ by $Lag
(M, L)$.  Then the relative Seidel morphism is  a functor 
$$S: \Pi Lag (M)  \to DF (M),$$
where $\Pi Lag (M)$ is the fundamental groupoid of $Lag (M)$ and $DF (M)$ is the Donaldson-Fukaya category of $M$, see also \cite{citeCorneaBiranLagrangianCobordism}, \cite{citeBiranCorneaLagCobGAFA}  which can be understood as an extension.
%

We sketch how this works. To a path $p$ in $Lag (M)$ from $L _{0} $  to $L _{1} $ we have have an associated Lagrangian subbundle $\mathcal{L} _{p} $ of $\mathcal{D} \times M$ over the boundary, as in Definition \ref{defInducedLoop}.
Extend this to a Hamiltonian structure $$\Theta _{p}
= (\mathcal{D} \times M, \mathcal{D}, \mathcal{L} _{p},
\widetilde{\mathcal{A}} _{0} ^{p}  )$$ where $\widetilde{\mathcal{A}} _{0} ^{p}$ is as in Lemma \ref{lemma:areadisk}. 
Assuming $\Theta _{p} $ is regular, we define $S ([p]) \in DF (L _{0}, L _{1})$ by
\begin{equation*}
   S ([p]) = \sum _{A} [ev (\Theta _{p}, A)],
\end{equation*}
where by monotonicity only finitely many $A$ can have non-zero contribution.

\subsection {Definition of the quantum Maslov
classes} \label{section:SeidelMorph}
Let $M$ be as before, and let $\mathcal{P} (L _{0}, L _{1})$ denote the space of smooth paths in $Lag (M)$
from $L _{0} $ to $L _{1} $, constant in $[0, \epsilon] \cup [1-\epsilon, 1]$ for some $0<\epsilon<1$.
There is then an additive group homomorphism:
\begin{equation} \label{eqPsiPaths}
\Psi: H_*(\mathcal{P} (L _{0}, L _{1}  ), \mathbb{Q}) \to FH (L _{0}, L _{1}
)
\end{equation}
defined analogously to above and to \cite{citeSavelyevQuantumcharacteristicclassesandtheHofermetric} in non-relative context. Although formally we will only need the restriction of $\Psi$ to spherical classes.

This works as follows.
To a smooth cycle $$f: B \to  \mathcal{P} (L _{0}, L _{1}  )
 $$ for $B$ a smooth closed oriented manifold, we may associate a Hamiltonian structure 
\begin{equation*}
   \{\mathcal{D} \times M, \mathcal{D}, \mathcal{L} _{b}\} _{B} , 
\end{equation*}
$\mathcal{L}_{b}:=\mathcal{L} _{f (b)}  $ a Lagrangian subbundle of $M \times \mathcal{D
} $ over
$\partial \mathcal{D} $ determined by $f (b)$ as before. The end of $\mathcal{D}$ here is negative.

Now let $\mathcal{A} _{0} $ be a Hamiltonian connection on
$[0,1] \times M$. And let $\mathcal{A} _{0} (L _{0} ) \subset \{1\} \times M $ denote the $\mathcal{A} _{0} $-transport over $[0,1]$ of $L _{0} \subset \{0\} \times M $.
Suppose that $\mathcal{A} _{0} (L _{0} ) $ is transverse to
$L _{1} $. 

For each $b$ the space of Hamiltonian connections $\mathcal{L} _{b} $-exact with respect to $\mathcal{A} _{0} $, (as in Section \ref{section:exact}) is contractible, c.f. \cite{citeAkveldSalamonLoopsofLagrangiansubmanifoldsandpseudoholomorphicdiscs}.
So we get an induced Hamiltonian structure:
$$\Theta _{f}= \{\mathcal{D} \times M, \mathcal{D}, \mathcal{L} _{b}, \mathcal{A} _{b} \}$$
well defined up to concordance.

%
We may then define $\Psi ([f])$ by:
\begin{equation*}
   \Psi ([f]) = \sum _{A} [ev (\Theta _{f}, A)],
\end{equation*}
where again by monotonicity only finitely many $A$ can give non-zero contribution. 
It is immediate that $\Psi$ is an additive group homomorphism. 
\begin{remark}
We should mention that the morphism $\Psi $ extends to
a certain functor to $DF (M)$,  see
\cite{citeBiranCorneaLagCobGAFA} for a related discussion,
in the degree 0 case. 
\end{remark}

Given the definition above,
$$[ev(\Theta', A _{0})] = \Psi (a)$$
clearly holds, as $A _{0} $ is the only
class that can contribute to $\Psi (a)$, since by the dimension formula \eqref{eq:dimensionMaslov} only a class $A$ with $Maslov ^{vert} (A ^{/} ) = -2 $ can contribute.

\section {Computation of
the quantum Maslov class $\Psi (a)$} \label{section:computeSeidel}

\subsection {Morse theory for the Hofer length
functional} 
Under certain conditions  the spaces of perturbation data for certain problems in
Gromov-Witten theory admit a Hofer like functional. Although these spaces of
perturbations are usually contractible, there may be a gauge group in the
background that we have to respect, so that working equivariantly there is topology. The reader may think of the analogous situation in
Yang-Mills theory \cite{citeAtiyahBottTheYang-MillsequationsoverRiemannsurfaces}.

Without elaborating too much, the basic idea of the computation that we will perform consists of cooling the perturbation data as much as possible (in the sense
of the functional) to obtain a mini-max (for the functional)
data, using which we may write down our moduli spaces
explicitly. This idea was first used in
the author's \cite{citeSavelyevVirtualMorsetheory}.
\subsubsection {Hofer length} \label{section:hoferlength} 
For $p: [0,1] \to \operatorname {Ham} (M,\omega)$  a smooth path, define
\begin{align*}
&  L ^{+} (p) :=   \int _{0} ^{1} \max _{M}  H ^{p} _{t} dt, \\
   &  L ^{-} (p) :=   \int _{0} ^{1} \max _{M}  (-H ^{p} _{t}) dt, \\
   & L ^{\pm} (p) :=   \max \{L ^{+} (p), L ^{-} (p)  \},  \\
\end{align*}
where $H ^{p}: M \times [0,1] \to \mathbb{R} $ generates
$p$, and is normalized by the condition that for each $t$, $H ^{p} _{t}:= H ^{p}| _{M \times \{t\}}   $ has mean 0, that is $\int _{M} H ^{p} _{t} dvol _{\omega} =0  $.
Also define
\begin{equation*}
   L ^{+} _{lag}: \mathcal{P} Lag (M)   \to \mathbb{R},     
\end{equation*}
\begin{equation*}
 L ^{+} _{lag}  (p) :=  \int _{0} ^{1} \max _{p (t)}  H ^{p } _{t} dt, 
\end{equation*}
$p (0) =L$ and where $H ^{p}: M \times [0,1] \to \mathbb{R}  $ is
normalized as above and generates a lift $\widetilde{p} $
of $p$ to $\operatorname {Ham} (M, \omega )$ starting at $id$. By lift we mean that $p (t)=\widetilde{p} (t) (p (0)) $. (That is $H ^{p}  $ generates a path in
$\operatorname {Ham} (M,  \omega) $, which moves $L_0$ along $p$.) Some theory of this latter functional is developed in \cite{citeLagrangianHofer}.
We may however omit the subscript $lag$ from notation, as usually there can be no confusion which functional we mean.

Note that $Lag ^{eq}(S ^{2} ) $ is naturally diffeomorphic to $S ^{2} $ and
moreover it is easy to see that the functional $L ^{+}| _{Lag ^{eq}  (S ^{2} )} $ is proportional to the
Riemannian length functional $L _{met} $ on the path space of $S ^{2} $, with its
standard round metric $met$.

Let now $L _{0}, L _{1} \in Lag ^{eq}  (S ^{2})$ be any transverse pair, and  $$f': S ^{2} \to \mathcal{P} ({L _{0},L _{1}}):= \mathcal{P} _{L _{0},L _{1}}  Lag ^{eq}(S ^{2} ), $$  
be the generator of the group $H _{2} (\mathcal{P} ({L _{0},L _{1}}), \mathbb{Z})$.
The idea of the computation is then this: perturb $f'$ to be transverse to the (infinite dimensional) stable manifolds
for the Riemannian length functional on $$\mathcal{P} ({L _{0},L _{1}}):= \mathcal{P} _{L _{0},L _{1}}  Lag ^{eq}(S ^{2} ), $$ push the cycle down 
by the ``infinite time'' negative gradient flow for this functional, and use the resulting representative to compute $\Psi (a=[f'])$. Although, we will not actually need infinite dimensional topology.
\subsubsection {The ``energy'' minimizing perturbation data} 
\label{secion:pertubationdata}
%
%
Classical Morse theory \cite{citeMilnorMorsetheory} 
tells us that the energy functional $$E (p) = \int _{[0,1]}  \langle \dot p (t), \dot p (t)  \rangle _{met}  \,dt  $$ on $ \mathcal{P} ({L_0, L_1})$ is
Morse non-degenerate with a single critical point in each degree. Consequently
$a$ (as a homology class) has a representative in the
2-skeleton of $ \mathcal{P} ({L_0, L_1})$, for the Morse cell decomposition induced by $E$. This follows by Whitehead's compression lemma which is as follows.
\begin{lemma} [Whitehead, see \cite{citeHatcherAlgebraic}]
   Let $(X,A)$ be a CW pair and let $(Y,B)$ be any pair with $B \neq \emptyset$. For each $n$ such that $X - A$ has cells of dimension $n$, assume that $\pi _{n} (Y,B,y _{0})=0$ for all $y _0 \in B$. Then every map $f:(X,A) \to (Y,B)$ is homotopic relative to $A$ to a map $X \to B$.
\end{lemma}
Suppose that $a$ has a representative 
$f': S ^{2} \to \mathcal{P} _{L_0, L_1} (S ^{2}) $ mapping into the $n$-skeleton $B ^{n} $ for the Morse cell decomposition for $E$, $n>2$.
Apply the lemma above with $(X,A) = (S ^{2}, pt )$, $Y = B ^{n} $ and $B = B ^{n-1} $ as above. Then the quotient $B ^{n} /B ^{n-1} $ is a wedge of $n$-spheres and since $\pi _{2} (S ^{n}) = 0 $ for $n>2$, $f$ can be homotoped into $B ^{n-1}   $ by the Whitehead lemma. Descend this way until we get a representative  mapping into $B ^{2} $.

Furthermore since
$\pi_2 (S ^{1} ) = 0 $  such a representative cannot
entirely lie in the 1-skeleton. 
It follows, since we have a single Morse 2-cell that 
there is a representative $f: S ^{2} \to  \mathcal{P} _{L_0, L_1} (S
^{2}) $, for $a$, s.t. the function ${f} ^{*} E  $  is Morse with a 
maximizer $\max$,  of index 2, and s.t. $\gamma _{0} = f (max)$ is the index $2$
geodesic.  We call such a representative \textbf{\emph{minimizing}}.
\begin{remark}
In principle there maybe more than one such maximizer $\max$, but recall
that we assumed that $a$ is the generator, so by further
deformation we may insure that there is only
one maximizer. The relevant representative $f$, with a single maximizer $\max$ as above, can also be 
constructed by hand.
\end{remark}
It follows that $\max$ is likewise the unique index 2 maximizer of the function ${f} ^{*} L _{met}  $ by the classical relation between the energy functional and length functional. And so $\max$ is the index 2 maximizer of $f ^{*}L ^{+}  $.
\subsubsection {The corresponding minimizing data} \label{section:distinguisheddata} 
\begin{lemma} \label{lemma:taut} There is a minimizing
representative $f _{0} $ for the class $a$ and a taut
Hamiltonian structure $$\Theta _{f _{0}}  =\{\mathcal{D}
\times S ^{2}, \mathcal{D}, \mathcal{L} _{f _{0} (b)},
\mathcal{A} _{b}  \},$$ 
satisfying:
\begin{equation} \label{eq.arealength}
  \forall b: \area (\mathcal{A} _{b}) = L ^{+} (f _{0}  (b)).
\end{equation} 

\end{lemma}
\begin{proof}
Note that a geodesic segment $p: [0,1] \to S ^{2} $ for the
round metric $met$ on $S ^{2} $ has a unique lift
$$\widetilde{p}: [0,1] \to PU (2) \simeq SO (3),$$
$\widetilde{p} (0)=id $ with $\widetilde{p} $ a segment of
a one parameter subgroup, and in this case $$L ^{+} _{lag}
(p)= L ^{+} (\widetilde{p}). $$ It then follows that for
a piecewise geodesic path $p$ in $S ^{2} $, there is
likewise a unique lift $ \widetilde{p}: [0,1] \to PU (2), $ satisfying $$L ^{+} _{lag} (p)= L ^{+} (\widetilde{p}). $$

Now, if $f$ is a minimizing representative of $a$, we may
homotop it to a likewise minimizing representative $f _{0}
$, so that for all $b$ $f _{0}  (b)$ is piecewise geodesic.  This follows by the piecewise geodesic approximation theorem Milnor~\cite[Theorem 16.2]{citeMilnorMorsetheory} of the loop space. 

Let $\mathcal{A} _{0} $ be the trivial Hamiltonian
connection on $[0,1] \times M$. Use the construction of
Lemma \ref{lemma:areadisk}, to get a family of Hamiltonian
connections $\{ \widetilde{\mathcal{A}} ^{f _{0}  (b)} _{0} \}   $. In
this case, since $\mathcal{A} _{0} $ is trivial $$\area
(\widetilde{\mathcal{A}} ^{f _{0}  (b)} _{0}  ) = L ^{+}(f _{0}  (b)). $$  

Set $A _{b} = \widetilde{\mathcal{A}} ^{f _{0}  (b)} _{0}$. It remains to verify that $\Theta _{f _{0} }  =\{\mathcal{D} \times S ^{2}, \mathcal{D}, \mathcal{L} _{f _{0} (b)}, \mathcal{A} _{b} \}$ is taut.
This follows by the following more general lemma.  
\begin{lemma} \label{lemma:tautLagS2} Let $Lag ^{eq}
(\mathbb{CP} ^{n})$ denote the space of oriented Lagrangian
submanifolds of $\mathbb{CP} ^{n} $  Hamiltonian isotopic to
$\mathbb{RP} ^{n} $. Then two loops $p _{1}, p _{2}: S ^{1}
\to  Lag ^{eq} (\mathbb{CP} ^{n}) $ are taut concordant, as
defined in Section \ref{section:introHamrigidity}, iff they are homotopic. 
\end{lemma}
\begin{proof}
Let $\mathcal{L}$ be a sub-fibration of $Cyl \times M$ as in the definition of taut concordance of loops. Let $\mathcal{A}$ be any $PU (n)$ connection on $Cyl \times \mathbb{CP} ^{n}$ which preserves $\mathcal{L}$ (there are no obstructions to constructing this). Then $R _{{A}}$ is a $\lie PU (n)$ valued 2-form, such that for all $v,w \in T _{z} Cyl $ the vector field $X=R _{\mathcal{A}} (z) (v,w)$ is tangent to $\mathcal{L} _{z} $. In particular if $H _{X} $ is the Hamiltonian generating $X$, then since $X$ is an infinitesimal unitary isometry preserving $\mathcal{L} _{z}, $ $H _{X} $ vanishes on $\mathcal{L} _{z} $. It follows by the definition of $\Omega _{\mathcal{A}} $, that it vanishes on $\mathcal{L}$ and so we are done.
\end{proof}
\end{proof}

So given $\{\mathcal{A} _{b} \}$ as in the lemma above, since 
$$\forall b: \area (\mathcal{A}  _{b}  ) = L ^{+} (f _{0}  (b)),$$ we immediately deduce:
\begin{lemma} \label{lemma.morseatmax} The function $area: b \mapsto \area
   (\mathcal{A}  _{b}  )$ has
   a unique maximizer, coinciding with the maximizer $\max$ of ${f} _{0}  ^{*}L
   _{met}  $ and $\area$ is Morse at $\max$ with index $2$.
\end {lemma}
\subsubsection {Finding class $A _{0} $ holomorphic sections for the data}
Let us now rename $f _{0} $ by $f$, $\mathcal{L} _{f _{0} (b) } $ by $\mathcal{L} _{b} $, and $\Theta _{f _{0} }  $ by $\Theta = \{\Theta _{b} \}$.

   As $ f (\max)$ is a geodesic for $met$,  its lift $\widetilde{f} (\max) $ to $SO (3)$ is a rotation around an axis intersecting $L _{0} = f (\max) (0)$ in a pair of points, in particular there is a unique point $$x _{\max} \in \bigcap _{t} (L _{t} = f (\max) (t))   $$ maximizing $H ^{\max} _{t} $ for each $t$. In our case this follows by elementary geometry but there is a more general phenomenon of this form c.f. \cite{citeLagrangianHofer}. 

Define  $$\sigma  _{\max}: \mathcal{D}  \to  \mathcal{D} \times S ^{2}   $$ to be the constant section $z \mapsto x _{\max }.
$ 
Then $\sigma _{\max} $ is a $\mathcal{A}  _{\max}  $-flat
section with boundary on $\mathcal{L}_{\max} $, and is consequently $J (\mathcal{A} _{\max} )$-holomorphic.

\begin{lemma}
   $$[\sigma _{\max}] = A _{0}.$$
\end{lemma}
\begin{proof}
 Set $$T ^{vert} _{z} \mathcal{L} _{\max}:= \{v \in T \mathcal{L} \subset T _{z} (\mathcal{D} \times S ^{2})\,|\, pr _{*} v =0 \}   $$ where $pr: \mathcal{D} \times S ^{2} \to \mathcal{D} $ is the projection. 
Denote by $$Lag (T _{x _{\max}} S ^{2} \simeq Lag (\mathbb{R} ^{2} ) \simeq S ^{1}  $$ the space of oriented linear Lagrangian subspaces of $T _{x _{\max} } S ^{2} $. 
Let $\rho$ be the path in $Lag (T _{x _{\max}} S ^{2})$ defined by 
\begin{equation*}
\rho (t) =  T ^{vert} _{(\zeta (t), x _{\max} )}   \mathcal{L} _{\max}, \quad t \in [0,1]
\end{equation*}
where $\zeta: \mathbb{R} \to \partial \mathcal{D}$ is a fixed parametrization as in Definition \ref{defInducedLoop}. 

By our conventions  for the Hamiltonian
vector field: $$\omega (X _{H}, \cdot ) = - dH (\cdot),$$
$\rho$ is a clockwise oriented path from $$T _{x _{\max} } L _{0} := T ^{vert} _{(\zeta (0), x _{\max} )} \mathcal{L} _{\max} $$ to $$T _{x
_{\max} } L _{1} := T ^{vert} _{(\zeta (1), x _{\max} )} \mathcal{L} _{\max}  $$ for the orientation induced by the complex orientation on $T _{x _{\max} } S ^{2}  $. 

By the Morse index theorem in Riemannian geometry \cite{citeMilnorMorsetheory} and by the condition that $f (\max)$ has Morse index 2, $\rho$  visits initial point $\rho (0)$ exactly twice if we count the start, as this corresponds to the geodesic $f (\max)$ passing through two conjugate points in $S ^{2} $. So the concatenation of $\rho$ with the minimal counter-clockwise path from $T _{x
_{\max} } L _{1}   $ back to $T _{x
_{\max} } L _{0}  $ is a degree $-1$ loop, if $S ^{1}  \simeq  Lag
(\mathbb{R}^{2})$ is given the counter-clockwise orientation. Consequently
   $$Maslov ^{vert} (\sigma _{\max} ^{/}) = -2,  $$ cf. Appendix
\ref{appendix:maslov}, in other words $[\sigma _{\max} ]=A _{0}$.
\end {proof}

\begin{proposition} \label{lemma:onlyelement}
   $(\sigma _{\max}, \max)$ is the sole element of
$
   \overline {\mathcal{M}} (\Theta, A _{0}).
$
\end{proposition}
\begin{proof}
By Stokes theorem, since $\omega$ vanishes on $\sigma _{\max} $, it is immediate:
\begin{equation} \label{eqLength1}
carea (\Theta _{\max}, A _{0})=- \int _{\mathcal{D}} \sigma _{\max}  ^{*}  \widetilde{\Omega} _{\max}
= L ^{+} (f (\max)). 
\end{equation}
   Moreover, since $\Theta=\{\Theta _{b} \} $ is taut $carea (\Theta _{b}, A _{0}  ) =  L ^{+} (f (\max))$.
So by \eqref{eq.arealength} and by Lemmas \ref{lemma:energypositivity1},
   \ref{lemmaInvariantPairing} we have:
\begin{equation*}
 L ^{+} (f (\max))  \leq \area
(\mathcal{A}_{b}) = L ^{+}(f (b)),
\end{equation*}
whenever there is an element 
$$(\sigma, b) \in  
\overline {\mathcal{M}} (\{
   \Theta _{b}  
\}, A _{0}).
$$
But clearly this is impossible unless $b= \max$, 
since $L ^{+} (f (b) ) < L ^{+} (f (\max) )  $
for $b \neq \max$. So to finish the proof of the proposition we just need:
\begin{lemma}
There are no elements $\sigma$ other than $\sigma
_{\max}$ of
the moduli space 
$$
\overline {\mathcal{M}} (
   \Theta _{\max}, A _{0}).
$$ 
\end{lemma}
\begin{proof}
We have by \eqref{eqLength1}, and by \eqref{eq.arealength}
$$0=  \langle [\widetilde{\Omega} _{\max} +
   \alpha _{\widetilde{\Omega}_{\max}}] , [\sigma_{\max} ]
   \rangle, $$
and so given another element $\sigma$   
we have:
\begin{equation*} \label{eq:vanishing}
   0=    \langle [\widetilde{\Omega}   _{\max} +
   \alpha _{\widetilde{\Omega} _{\max}}] , [\sigma] \rangle.
\end{equation*}
It follows that $\sigma $ is necessarily $\widetilde{\Omega}   _{\max}
$-horizontal, since 
$$(\widetilde{\Omega}   _{\max} +
\alpha _{\widetilde{\Omega}   _{\max}  }) (v, J _{\widetilde{\Omega}  _{\max} }v) \geq 0.$$
Since $J _{\widetilde{\Omega}  _{\max}}$ by assumptions preserves the vertical and $\mathcal{A} _{\max} $-horizontal subspaces of $T (\mathcal{D} \times S ^{2} )$, and since the 
inequality is strict for $v$ in the vertical tangent bundle of
$$S ^{2} \hookrightarrow \mathcal{D} \times S ^{2}   \to \mathcal{D},
$$ 
the above inequality is strict whenever $v$ is not horizontal.
So $\sigma$ must be $\mathcal{A} _{\max} $-horizontal. But then $\sigma = \sigma _{\max} $ since $\sigma _{\max} $ is the only flat
section asymptotic to $\gamma _{0} $.
\end {proof}
\end {proof} 
\subsubsection {Regularity}
 It will follow that $$\Psi (a) = \pm [\gamma _{0} ]$$ if we knew
that
$(\sigma_{\max}, \max)$ be a regular element of 
$$ 
\overline {\mathcal{M}} (\{
   \Theta _{b} 
\}, A _{0}).
$$
We won't answer directly if $(\sigma _{\max}, \max)$ is regular, 
although it likely is. But
it is regular after a suitably  small Hamiltonian perturbation of the family
$\{\mathcal{A}  _{r} \}$ vanishing at $\mathcal{A}  _{\max} $.  
We call this essentially automatic regularity. 
\begin{lemma} \label{lemma.perturbedreg} 
There is a family $\{\mathcal{A} ^{reg}_{b}  \}  $ arbitrarily $C ^{\infty}
$-close to $\{\mathcal{A} _{b} \}$ with  $\mathcal{A} ^{reg}_{\max} =
\mathcal{A} _{\max}  $ and such that
\begin{equation} \label{eqModuliReg} 
\overline {\mathcal{M}} (\{
   \mathcal{D} \times S ^{2}, \mathcal{D}, \mathcal{L} _{b}, \mathcal{A} ^{reg}  _{b} 
\}, A _{0}),
\end{equation} 
is regular, with $({\sigma}_{\max}, {\max}  )$  its sole element. In particular $$\Psi (a) = \pm [\gamma _{0} ].$$
\end{lemma}
\begin{proof}
The associated real linear Cauchy-Riemann operator $$D_{\sigma_{\max}}: \Omega ^{0} (\sigma _{\max} ^{*} T
   ^{vert} \mathcal{D} \times S ^{2} _{\max})    \to \Omega ^{0,1} (\sigma _{\max} ^{*} T
^{vert} \mathcal{D} \times S ^{2} _{\max}),$$ 
has no kernel, by Riemann-Roch 
\cite[Appendix
C]{citeMcDuffSalamonJholomorphiccurvesandsymplectictopology}, 
as the vertical Maslov number of $[\sigma _{\max} ]$ is $-2$.
And the Fredholm index of
$({\sigma}_{\max}, {\max}  )$ which is -2,  is -1 times the Morse index of the function $\area$
at $\max$, by
Lemma \ref{lemma.morseatmax}.
Given this, our lemma follows by a direct translation of
\cite[Theorem 1.20]{citeSavelyevMorsetheoryfortheHoferlengthfunctional}, itself
elaborating on the argument in
\cite{citeSavelyevVirtualMorsetheory}.
\end{proof}
To summarize: 
\begin{theorem} \label{thmNonZeroPsi} For $0 \neq a \in H _{2} (\mathcal{P}_ {L _{0}, L _{1}} Lag (S ^{2}), \mathbb{Z})$,
\begin{equation*}
  0  \neq  \Psi (a) \in HF (L _{0}, L _{1}).
\end{equation*}
\end {theorem} 
\begin{proof} We have shown that $0  \neq  \Psi (a) \in HF (L _{0}, L _{1})$, for $a$ the generator of the group $H _{2} (\mathcal{P}_ {L _{0}, L _{1}} Lag (S ^{2}), \mathbb{Z}) $. Since $\Psi$ is an additive group homomorphism the conclusion follows.
   \end{proof}
  %
\section{Finishing up the proof of Lemma \ref{lemma:nontrivialelement}}
The existence of small data $\mathcal{D} $ is proved in
Section \ref{sec:constructionSmallData}.
Given this existence, starting with \eqref{eq:mu=H^}  we showed that $[\mu ^{4} _{\Sigma ^{4}} (\gamma _{1} , \ldots , \gamma _{4} )]$ is non-vanishing in Floer homology iff $$[ev (\mathcal{H} _{0}, A _{0})] \in HF (L _{0}, L _{4}), $$ is non-vanishing. We then use Lemma \ref{lemma:concordant} to identify $[ev (\mathcal{H} _{0}, A _{0})]$ with $[ev (\Theta', A _{0})]$, which is also identified with $\Psi (a)$, for a certain spherical 2-class $a$. Finally, in Section \ref{section:computeSeidel} we compute $\Psi (a)$ and show that it is non-zero. This together with Lemma \ref{lemma:corellator} imply Lemma \ref{lemma:nontrivialelement}.
\qed
\section{Proof of Theorem \ref{thm:Hofer}} \label{section:applicationHofer}
Suppose otherwise, so that $$\min_{f, [f] = a'} \max _{b \in S ^{2} } L ^{+} (f (b)) = U < \hbar,$$
   for $a' = i_*g$ as in the statement of the theorem. Fix $L _{1} \in Lag ^{eq} (S ^{2} ) $ so that $L _{0} $ intersects $L _{1} $ transversally, and so that there is a geodesic path $p _{0} \in \mathcal{P} Lag ^{eq}  (L _{0}, L _{1}  ) $  with $$\kappa := L ^{\pm} (\widetilde{p} _{0} )   < \epsilon = (\hbar - U)/2. $$ Here $\widetilde{p} _{0}  $ is the geodesic lift to $PU (2)$ starting at $id$. Then concatenating $f$ with $p _{0} $ we obtain a smooth family of paths 
\begin{align*}
   & g: S ^{2} \to \mathcal{P} (L _{0},L _{1}) \\
   & g (0) = p _{0},
\end{align*}
and $g$ represents the previously appearing class $a$, that
is the generator of the group $$\pi _{2}(
\mathcal{P} (L _{0},L _{1}), p _{0}). $$ 
Let $$\{\Theta _{b} \}= \{
   \mathcal{D} \times S ^{2}, \mathcal{D}, \mathcal{L} _{b}, \mathcal{A}_{b} 
\} _{\mathcal{K} = S ^{2}} ,$$ be the corresponding
Hamiltonian structure, where $\mathcal{A} _{b}$ is as in
Lemma \ref{lemma:taut}, defined with respect to $g$, and where $\mathcal{L} _{b}:=\mathcal{L} _{g (b)} $. 
In particular, $\{\Theta _{b} \}$ is taut and satisfies:
\begin{equation} \label{eq:hbarlower}
\forall b \in S ^{2}: \area (\mathcal{A} _{b} ) = L ^{+} (g (b))  < \hbar - \kappa.
\end{equation}

By assumption that each $f (b)$ is taut concordant to the
constant loop at $L _{0} $, each $\Theta _{b} $ is taut
concordant to $$\mathcal{H} = (\mathcal{D} \times S ^{2},
\mathcal{D}, \mathcal{L} _{0}, \mathcal{A}),$$ where
$\mathcal{L} _{0} = \mathcal{L} _{p _{0}}  $,
${\mathcal{A}} = \widetilde{\mathcal{A}} _{0} ^{p _{0}}
$, where $\widetilde{\mathcal{A}} _{0} ^{p _{0}} $   is as
in Lemma \ref{lemma:areadisk}, for $\mathcal{A} _{0} $ the
trivial connection. 

Let $\Theta (L, \mathcal{A} _{0})$ be the construction as in \eqref{eq:ThetaL}.
Then for each $b$, $$\Theta _{b} ^{/0}
: = \Theta _{b} \# _{0} \Theta (L _{0},
\mathcal{A} (L _{0}, L _{0})) $$  is taut concordant to
$\mathcal{H} ^{/0} $ (which is defined analogously)  by Lemma \ref{lemma:immediategluing}.
On the other hand, by Lemma \ref{lemma:tautLagS2}
$\mathcal{H} ^{/0} $ is taut concordant to the trivial
Hamiltonian structure $(D ^{2}  \times S ^{2}, D ^{2},
\mathcal{L} _{tr} , \mathcal{A} _{tr} )$, where $\mathcal{L} _{tr} $ the trivial bundle with fiber $L _{0} $ and $\mathcal{A} _{tr} $ the trivial Hamiltonian connection. 
So for each $b$:
\begin{equation} \label{eq:careablah}
\carea (\Theta _{b} ^{/0} , A _{0}) =  \carea (\mathcal{H}
^{/0} , A _{0}) = \hbar.
\end{equation}
%

Now by Theorem \ref{thmNonZeroPsi} $$ev (\{\Theta _{b} \}, A _{0} ) = \Psi (a) \neq 0.$$ And so:
\begin{equation*}
\overline {\mathcal{M}} (\{\Theta _{b}\}, A _{0}) \neq \emptyset,
\end{equation*}
but this contradicts the conjunction of \eqref{eq:hbarlower}, \eqref{eq:careablah},  and Lemma \ref{lemma:gluinglowerbound}.

\qed
\section{Singular and simplicial connections and curvature bounds} \label{sec:qcurvature}
Let $\mathcal{A}$ be a $G$ connection on a principal $G$ bundle $P \to \Delta ^{n} $, and the Finsler norm $\mathfrak{n}$ on $\lie G$ be as in Section \ref{sec:non-metric} of the introduction. 
As previously discussed, a given system $\mathcal{U}$ in particular specifies maps: 
\begin{equation*}
u (m_1, \ldots, m_n, r,n): \mathcal{S} _{r} \to \Delta ^{n},
\end{equation*}
where $r \in \overline {\mathcal{R}} _{n}$, $\mathcal{S} _{r} $ is the fiber of ${\mathcal{S}}  ^{\circ} _{n}    $ over $r$, and where $(m_1, \ldots, m _{n} )$ is the composable chain of morphisms in $\Pi (\Delta ^{n} )$, $m _{i} $ being the edge morphism from the vertex $i-1$ to $i$.
%
Then define
\begin{equation} \label{eq:definitionArea}
   \area _{\mathcal{U}} (\mathcal{A}) = \sup _{r}
   \area _{\mathfrak n}  (u (m_1, \ldots, m_n, r,n) ^{*} \mathcal{A}),
\end{equation}
where $\area _{\mathfrak n} $ on the right hand side is as defined in equation \eqref{eq:areasurface}.
In the case $G= \operatorname {Ham} (M,\omega)$ we take
$$\mathfrak n: \lie \operatorname {Ham} (M, \omega) \to \mathbb{R}$$ to be $$\mathfrak n (H) = |H| _{+} =  \max _{M} H. $$

Let $\omega$ be the area 1 Fubini-Study symplectic 2-form on
$M=\mathbb{CP} ^{1}$. Then the pull-back by the natural map
$$\lie h: \lie PU (2) \to \lie \operatorname {Ham} (\mathbb{CP} ^{1},\omega) \simeq C ^{\infty} _{0} (\mathbb{CP} ^{1} )$$ of the semi-norm: $|H| _{+} =  \max _{M} H$ is the operator norm on $PU (2)$, up to normalization. This will be used to get the specific form of Theorem \ref{thm:lowerboundsingular},  from the more general form here.

%
\subsection {Simplicial connections} \label{section:simplicialconnections}
We now introduce a certain abstraction of simplicial connections, which can partly be understood as simplicial resolutions of singular connections.
Let ${G} \hookrightarrow P \to X$ be a principal $G$ bundle, where $G$ is a Frechet Lie group.
Denote by $X_{\bullet} $ the simplicial set whose set of $n$-simplices, $X _{\bullet} (n) $, consists smooth maps $\Sigma: \Delta ^{n} \to X$, with $\Delta ^{n} $ standard topological $n$-simplex with vertices ordered $0, \ldots, n$. 
And denote by $Simp(X _{\bullet})$ the category with objects $\cup _{n} X _{\bullet} (n)$ and with $hom (\Sigma _{0}, \Sigma _{1})$ commutative diagrams:
\begin{equation*}
\begin{tikzcd}
   \Delta ^{n} \ar [rd, "\Sigma _{0}"] \ar[r, "mor"] & \Delta ^{m} \ar[d, "\Sigma _{1}" ] \\
& X,
 \end{tikzcd}
\end{equation*}
for $mor$ a simplicial face map, that is an injective affine map preserving order of the vertices.
\begin{definition} \label{def:singualCon}
 Define a \textbf{\emph{simplicial $G$-connection $\mathcal{A}$}} on $P$ to be the following data:
 \begin{itemize}
   \item For each $\Sigma: \Delta ^{n} \to X  $  in $X _{\bullet} (n) $ a smooth  $G$-connection $\mathcal{A} _{\Sigma} $ on $\Sigma ^{*}P \to \Delta ^{n} $, (a usual Ehresmann $G$-connection.)
\item For a morphism $mor: \Sigma _{0} \to \Sigma _{1}  $ in $Simp (X _{\bullet} )$, 
we ask that $mor ^{*} \mathcal{A} _{\Sigma _{1}} = \mathcal{A} _{\Sigma _{0} } $. 
\end{itemize}
\end{definition}
\begin{example} \label{example:exampleSimplicial}
   If $\mathcal{A}$ is a smooth $G$-connection on
   $P$, define a simplicial connection by
   $\mathcal{A} _{\Sigma} = \Sigma ^{*}
   \mathcal{A}  $ for every simplex $\Sigma \in X
   _{\bullet} $. We call such a simplicial
   connection \textbf{\emph{induced}}.
\end{example}
If we try to ``push forward'' a simplicial connection to get a ``classical'' connection on $P$ over $X$, then we get a kind of multi-valued singular connection. Multi-valued because each $x \in X$ may be in the image of a number of $\Sigma: \Delta ^{n} \to X $ and $\Sigma$ itself may not be injective, and singular because each $\Sigma$ is in general singular so that the naive push-forward may have blow up singularities. We will call the above the naive pushforward of a simplicial connection. 

\begin{proof} [Proof of Theorem
   \ref{thm:lowerboundsingular} and Corollary \ref{corol:example}] We will prove this by way of a stronger result.
Let $P$ be a Hamiltonian fibration $S ^{2} \hookrightarrow P \to S ^{4}  $.
And let $\mathcal{A}$ a simplicial $\mathcal{G}  = \operatorname {Ham} (S
^{2},  \omega) $ connection on $P$. Let $\sigma ^{1} _{0}
\in S ^{4} _{\bullet} $ be the degenerate $1$-simplex at $x
_{0}$, in other words the constant map: $\sigma ^{1} _{0}:
[0,1] \to x _{0}. $  Let $\kappa$ be the $L ^{\pm} $-length of the holonomy path of $\mathcal{A} _{\sigma ^{1} _{0}}$ over $[0,1]$. 

Let $\Sigma _{\pm} \in S ^{4} _{\bullet} (4)  $ be a
complementary pair as in Section
\ref{section:model}.
The connection $\mathcal{A} $ gives us a simplicial connection  as in Example
\ref{example:exampleSimplicial}.
By inductive procedure as
in Part I, Lemma 5.6, we may find 
perturbation data $\mathcal{D}$ for $P$ so that for this data
\begin{align} 
& \forall r: pr _{1} \mathcal{F} (L ^{0} _{0}, \ldots, L ^{n} _{0},  \Sigma _{\pm},r)  \simeq _{\delta}  
u (m_1, \ldots, m_s, r,n) ^{*} \mathcal{A} _{\Sigma _{\pm} }, \label{eq:F=G} \\
 & \mathcal{A} (L _{0}, L _{0}) \simeq _{\delta}  \mathcal{A} _{\sigma ^{1} _{0}}, \label{eq:F=G2}
\end{align}
where $L ^{i} _{0}  $ are the objects as before,
where $\simeq _{\delta} $ means $\delta$-close in the metrized $C ^{\infty} $ topology, and $\delta$ is as small as we like.  Here we are using notation of Part I as before. Set $$\widetilde{u} _{r}:= \Sigma _{-}  \circ u(m_1, \ldots, m _{4},r,4),$$
so $\widetilde{u} _{r}: \mathcal{S} _{r} \to S ^{4}$.
Set $\widetilde{\mathcal{S}} _{r}:= \widetilde{u} _{r} ^{*} P, $
set 
$\mathcal{A}'_{r}:= pr _{1} \mathcal{F} (L ^{0} _{0}, \ldots, L ^{n} _{0},  \Sigma _{-},r)$
and set $$\{\Theta _{r}\} :=  \{\widetilde{\mathcal{S}} _{r},
   {\mathcal{S} _{r}}, \mathcal{L} _{r}, \mathcal{A}'_{r}\}. $$
\begin{definition}\label{def:perfect}
We say that $\mathcal{A}$ is \textbf{\emph{perfect}}, if 
for every arbitrarily small $\delta$ as above, $\mathcal{A}
(L _{0}, L _{0}  )$, as above, can be chosen  so that the
   corresponding Floer chain complex $CF (L _{0},
   L _{1}, \mathcal{A} (L _{0}, L
   _{0}) )  $ is perfect.
\end{definition}
\begin{theorem} \label{prop:alternative} 
Let $\mathcal{A}$ be a perfect simplicial Hamiltonian
connections on $P$. If 
$P$ is non-trivial as a Hamiltonian bundle
 then 
\begin{equation*}
   (\area _{\mathcal{U}}  (\mathcal{A} _{\Sigma _{+}  })  \geq \hbar -5 \kappa) \lor (\area _{\mathcal{U}}  ( \mathcal{A} _{\Sigma _{-}} ) \geq \hbar - 5\kappa),
\end{equation*}
\end{theorem}
\begin{proof}
Suppose 
\begin{equation} \label{eq:arealesshalf}
\area _{\mathcal{U}}  (\mathcal{A} _{\Sigma _{+} } ) < \hbar - 5 \kappa.
\end{equation}
Then by \eqref{eq:F=G}, \eqref{eq:F=G} and by Lemma
\ref{lemma:gluinglowerbound} $\mathcal{D}$, as defined
above, can be assumed to be small provided $\delta$ is chosen to be sufficiently small. Take the unital replacement as in Lemma \ref{lemma:corellator}. Since we know that $K (P)$ does not admit a section by Theorem \ref{thm:noSection}, the simplex $T$ of the Lemma \ref{lemma:corellator} does not exist. 
Hence again by this lemma 
$$ev (\{\Theta _{r} \}, A _{0}) = [\mu ^{4}
_{\Sigma _{-}} (\gamma _{1},  \ldots,  \gamma _{4} )] \neq 0.$$  
In particular $$\overline {{\mathcal{M}}}(\{\Theta _{r} \}, A _{0} ) \neq \emptyset.$$ 

So by Lemma \ref{lemma:gluinglowerbound} there exists an $r _{0}$ so that
\begin{equation} \label{eq:ineq1}
\area (\mathcal{A}'_{r _{0}}) \geq \hbar - 5 \kappa'.
\end{equation}
where $\kappa'$ denotes the $L ^{\pm} $ length of the
holonomy path in $\operatorname {Ham} (S ^{2},  \omega) $ of $\mathcal{A} (L _{0}, L _{0}) $. 
By \eqref{eq:F=G2} $\kappa' \to \kappa$ as $\delta
\to 0$.
By \eqref{eq:F=G}, \eqref{eq:ineq1}
passing to the limit as $\delta \to 0$ we get:
$$\area _{\mathcal{U}}  ( \mathcal{A} _{\Sigma _{-}} ) \geq \hbar - 5\kappa.$$
\end{proof}
\begin{corollary}
   \label{corol:alternative} 
Let $\mathcal{A}$ be a $PU (2) $ connection 
on a non-trivial principal $PU (2) $
bundle $P \to S ^{4}$. 
 Then 
\begin{equation*}
   (\area _{\mathcal{U}}  (\mathcal{A} _{\Sigma
   _{+}  })  \geq \hbar -5 \kappa) \lor (\area
   _{\mathcal{U}}  ( \mathcal{A} _{\Sigma _{-}} )
   \geq \hbar - 5\kappa).
\end{equation*}
\end{corollary}
\begin{proof}
 A
simplicial $PU (2)$ connection $\mathcal{A}$ on a
principal $PU (2)$ bundle $PU (2) \hookrightarrow
P' \to S ^{4} $ is automatically perfect, when
understood as a Hamiltonian connection on the
associated bundle $S ^{2}  \hookrightarrow P \to S
^{4} $. So that this is an immediate consequence
of the theorem above.
\end{proof}
To prove Corollary \ref{corol:example},
we just note that
for an induced simplicial Hamiltonian connection
$\mathcal{A} $, as defined in
Example \ref{example:exampleSimplicial}, $\mathcal{A} _{\sigma ^{1} _{0}}$
is trivial.  And hence $\mathcal{A}$ is
automatically perfect.  So that this corollary follows by Theorem \ref{prop:alternative}.

\end {proof}

\appendix 

\section {Homotopy groups of Kan complexes} \label{appendix:Kan}
For convenience let us quickly review Kan
complexes just to set notation.  This notation is
also used in Part I. Let $$\Delta ^{n}_{\bullet } (k):= hom _{\Delta} ([k], [n]),  $$  be the standard representable $n$-simplex,
where $\Delta$ is the simplicial category  with
objects ordered finite sets $[n] = \{1, \ldots,
n\}$ and morphisms order preserving set maps. 

Let $\Lambda ^{n} _{k} \subset \Delta ^{n} _{\bullet } $ denote the sub-simplicial set corresponding to the ``boundary'' of $\Delta ^{n} _{\bullet}   $ with the $k$'th face removed, $0 \leq k \leq n$. 
By $k'th$ face we mean the face opposite
to the $k$'th vertex. 
Let $X _{\bullet}$ be an abstract simplicial set.
A simplicial map $$h: \Lambda ^{n} _{k} \subset \Delta ^{n} _{\bullet}  \to X _{\bullet}  $$ will be called a \textbf{\emph{horn}}.
A simplicial set $S _{\bullet}$ is said to be a \textbf{\emph{Kan complex}} if for all $n,k \in \mathbb{N}$ given a diagram  with solid arrows
\begin{equation*} 
\begin {tikzcd}
   \Lambda ^{n} _{k} \ar [r, "h"]  \ar [d, "i"] & S _{\bullet} \\
 \Delta ^{n} _{\bullet} \ar [ur, dotted, "\widetilde{h}"]  &, \\ 
\end{tikzcd}
\end{equation*}
there is a dotted arrow making the diagram commute. 
The map $\widetilde{h} $ will be called \textbf{\emph{the Kan filling}} of the horn $h$. The $k$'th face of $\widetilde{h} $ will be called \textbf{\emph{Kan filled face along $h$}}. 
As before we will denote Kan complexes and $\infty$-categories by calligraphic letters.

Given a pointed Kan complex $(\mathcal{X} ,x) $ and $n \geq 1$ the
\emph{$n$'th simplicial homotopy group} of
$(\mathcal{X},x) $: $\pi _{n} (\mathcal{X}, x ) $ is defined to be the set
of equivalence classes of maps
$$\Sigma: \Delta ^{n} _{\bullet}   \to \mathcal{X},$$ 
such that $\Sigma$ takes $\partial \Delta
^{n}_{\bullet} $ to $x _{\bullet} $, with the
latter denoting the image of $\Delta ^{0}
_{\bullet} \to \mathcal{X}   $, induced by the vertex inclusion $x \to X$.

More precisely, we have a commutative diagram:
\begin{equation*} 
\begin{tikzcd}
   \Delta ^{n} _{\bullet}  \ar [rd, "\Sigma"] \ar[r] & \Delta ^{0} _{\bullet}  \ar[d, "x" ] \\
& \mathcal{X}.
 \end{tikzcd}
\end{equation*}
\begin{example}
   \label{example:singularset} 
When $\mathcal{X} = X _{\bullet} $ is the
simplicial set of singular simplices of a  
topological space $X$, the maps above are in complete correspondence with maps:
\begin{equation*}
\Sigma: \Delta ^{n}  \to X,   
\end{equation*}
taking the topological boundary of $\Delta ^{n} $ to $x$.

\end{example}

For $X _{\bullet }$ general simplicial set, a pair of maps $\Sigma _{1}: \Delta ^{n}
_{\bullet}   \to X _{\bullet}, \Sigma _{2}: \Delta
^{n} _{\bullet}   \to X _{\bullet},$ are
equivalent if there is a diagram, called simplicial homotopy: 
\begin{equation*}
\begin{tikzcd}
 \Delta ^{n} _{\bullet} \ar[rd, "\Sigma _{1}"] \ar
 [d, "i_0"] & \\
 \Delta ^{n} _{\bullet}  \times I _{\bullet}   \ar
 [r,""] & X _{\bullet} \\
 \Delta ^{n} _{\bullet}  \ar [u, "{i_1}"]  \ar
 [ru, "{\Sigma _{2}}" ].  &  
\end{tikzcd}
\end{equation*}
such that $\partial \Delta ^{n}_{\bullet} \times I _{\bullet}$ is taken by $H$
to $x _{\bullet}  $.
The simplicial homotopy groups of a Kan complex $
(\mathcal{X},x)$ coincide with the
classical homotopy groups of the geometric
realization $(|\mathcal{X}|,x)$. 
\begin{proof} [Proof of Lemma \ref{lemma:kanfib}]  
We refer the reader to Part I, Appendix A.2,
for more details on the notions here.  We prove a stronger claim. 
\begin{lemma}
   \label{lemma:innerfibkan} Let $p: \mathcal{Y}  \to \mathcal{X}$
be an inner fibration of quasi-categories
   $\mathcal{Y},\mathcal{X}$,
with $\mathcal{X}$ a Kan complex. And let 
     $K (\mathcal{Y}) \subset \mathcal{Y}
 $ denote  the maximal Kan subcomplex.
Then $p: K (\mathcal{Y})  \to \mathcal{X} $ is a Kan
fibration.
\end{lemma}
The above is probably well known,  but it is simple
to just provide the proof for convenience.   
\begin{proof}
By definition of an inner fibration,  whenever we are given a commutative diagram with
solid arrows and with $0<k<n$:
\begin{equation} 
\begin{tikzcd} \label{diagram:inner2}
   \Lambda ^{n} _{k}  \ar [r, "\sigma"]  \ar[d, hookrightarrow] &
K (\mathcal{Y}) \ar [r, hookrightarrow] &  \mathcal{Y} 
\ar [dl, "p"]  \\
 \Delta ^{n} _{} \ar [r] \ar [urr, "\Sigma ", dashrightarrow] & \mathcal{X}
 _{}, \\
\end{tikzcd} 
\end{equation}
there exists a dashed arrow $\Sigma $ as
indicated, making the whole diagram commutative.  
When $n>2$ the edges, i.e. 1-faces, of $\Sigma$
   are all automatically isomorphisms in $\mathcal{Y} _{}$, as $\Sigma$ extends
$\sigma$, 
and all edges of $\sigma$ are isomorphisms
by definition.  For $n=2$  the edges of $\Sigma$
are either edges of 
$\sigma$, or are compositions of edges of $\sigma$
in the quasi-category $\mathcal{Y}$, and hence again
always invertible. 

It follows that 
$\Sigma$ maps into $K (\mathcal{Y}) \subset \mathcal{Y} $. Since the starting
diagram was arbitrary, we just proved
that $p: K (\mathcal{Y} _{}) \to \mathcal{X} _{} $ is an inner fibration.  In particular the pre-images $p
^{-1} (\Sigma (\Delta  ^{n} _{})) \subset K
(\mathcal{Y} _{ }) 
 $ are quasi-categories, for all $n$, where $\Sigma: \Delta  ^{n}
_{} \to \mathcal{X}$ is any $n$-simplex,  see
Part I, Appendix A.2. But $K (\mathcal{Y}) $ is a Kan
complex, so that also the above pre-images $p
^{-1} (\Sigma (\Delta  ^{n} _{})) $ are Kan
complexes. It readily follows from this that   
$p: K (\mathcal{Y} _{ }) \to \mathcal{X} _{ } $ is a Kan
fibration.
\end {proof}  
The main lemma then follows, since if $p: \mathcal{Y} _{} \to \mathcal{X} _{ } $ is a
categorical fibration,  it is in particular an
inner fibration.

\end{proof}
\section {On the Maslov number and dimension formula} \label{appendix:maslov}
Let ${S}$ be a Riemann surface with boundary and a strip end
structure as previously.

Let $\mathcal{V} \to S$ be a rank $r$ complex vector bundle,
trivialized at the open ends $\{e _{i} \}$, so that we have distinguished bundle charts $
[0,1] \times (0,\infty) \times \mathbb{C} ^{r}  \to \mathcal{V}$ at the positive ends. (Similarly, for negative ends.) 

Let $$\Xi \to \partial S  \subset S $$
be a totally real rank $r$ subbundle of $\mathcal{V}$, which  is constant in the
coordinates $$ [0,1] \times (0, \infty) \times \mathbb{C}^{r},  $$
at the positive ends, again similarly with negative ends. 

For each (positive) end $e _{i} $ and its chart  $e _{i}:
[0,1] \times (0, \infty) \to S$, let
$b ^{j}  _{i}: (0, \infty) \to \partial S  $, $j=0,1$ be the restrictions of $e _{i} $ to $\{i\} \times (0, \infty)$.



We then have a pair of real vector spaces $$\Xi _{i} ^{j} =  \lim _{\tau \mapsto \infty} \Xi|_{b
^{j} _{i} 
 (\tau) }.
 $$

There is a Maslov number $Maslov (\mathcal{V}, \Xi, \{\Xi
^{j}  _{i} \})$ associated to this data, and which we now
briefly describe.
In the case 
$\Xi _{i} ^{0} = \Xi _{i} ^{1}$, let $(\mathcal{V} ^{/}, \Xi
^{/}  )$ be obtained from $(\mathcal{V}, \Xi,
\{\Xi ^{j}  _{i} \})$ by capping off each $e _{i} $ end of $\mathcal{V} \to S$. Here the capping operation is similar to the one in Section
\ref{section:gluing}.
Then $Maslov (\mathcal{V}, \Xi, \{\Xi ^{j}  _{i} \})$
coincides with the boundary Maslov index of $(\mathcal{V} ^{/}, \Xi ^{/}  )$ in the sense of
\cite[Appendix
C3]{citeMcDuffSalamonJholomorphiccurvesandsymplectictopology}.

When $\Xi _{i} ^{0}$ is transverse to $\Xi _{i} ^{1}$ for each $i$,  $Maslov (\mathcal{V}, \Xi, \{\Xi ^{j}
_{i} \})$ is obtained as the Maslov index for the
modified pair $(\mathcal{V} ^{/}, \Xi ^{/}  )$
obtained by again capping off the ends $e _{i} $ via gluing (at each end $e _{i} $) with $$(\mathcal{D} \times \mathbb{C}^{r}, 
 \widetilde{\Xi}, \{\widetilde{\Xi}_{0} ^{j}  \} ),$$ where
 $\mathcal{D}$ is as before.  Here  
 $\widetilde{ {\Xi}} ^{0} _{i} = { {\Xi}} ^{1} _{0} $ and 
$\widetilde{ {\Xi}} ^{1} _{i} = { {\Xi}} ^{0} _{0} $, while
$ \widetilde{\Xi}$ over the boundary of $\mathcal{D}$ is determined by the ``shortest path'' from 
$\widetilde{ {\Xi}} ^{0} _{0}$ to $\widetilde{ {\Xi}} ^{1}
_{0}$, which means the following. 
As $\widetilde{ {\Xi}} ^{0} _{0}$ to $\widetilde{ {\Xi}} ^{1}
_{0}$ are a pair of transverse, totally real subspaces, up
to a complex isomorphism of $\mathbb{C} ^{r} $ (whose choice
will not matter), we may identify them with the subspaces
$\mathbb{R} ^{r} $, and $i \mathbb{R} ^{r} $. After this identification our shortest path is just $e ^{i \theta} \mathbb{R} ^{r}  $, $\theta \in [0, \pi _{2} ]$.

%
Let $D$ be a real linear Cauchy-Riemann operator on
$\mathcal{V}$, which in particular is an operator: 
\begin{equation*}
D: \Omega ^{0} _{\Xi} (S, \mathcal{V}) \to \Omega ^{0,1}
_{\Xi} (S, \mathcal{V}),  
\end{equation*}
where $\Omega ^{0} _{\Xi} (S, \mathcal{V})$ denotes the
space of smooth $\mathcal{V} $-valued $0$-forms (i.e.
sections) satisfying $\theta (\partial S) \subset \Xi $, and
$\Omega ^{0,1} _{\Xi} (S, \mathcal{V})$ denotes the
analogous space of smooth complex anti-linear 1-forms.

Suppose further that $D$ is asymptotically $\mathbb{R} ^{}
$-invariant in the strip end coordinates at the ends.
After standard Sobolev completions,  the Fredholm index of
$D$ is given by:
\begin{equation*}
r \cdot \chi (S) + Maslov (\mathcal{V}, \Xi, \{\Xi _{i} \}).
\end{equation*}
The proof of this is analogous to \cite[Appendix
C]{citeMcDuffSalamonJholomorphiccurvesandsymplectictopology}, we can also reduce it to that statement via a gluing
argument. (This kind of argument appears for instance in
\cite{citeSeidelFukayacategoriesandPicard-Lefschetztheory})
\subsection {Dimension formula for moduli space of sections}
Suppose we have a Hamiltonian structure $\Theta = (\widetilde{S}, S,
\mathcal{L}, \mathcal{A} ) $.  Suppose that either the
corresponding Lagrangian submanifolds $$L _{i} ^{j} =  \lim _{\tau \mapsto \infty} \mathcal{L}|_{b
^{j} _{i} (\tau) },$$ intersect transversally (identifying
the corresponding fibers) or coincide. (Similarly for
negative ends.) 

Let $A \in H _{2} ^{sec} (\widetilde{S}, \mathcal{L} ) $,
with the latter as in Section \ref{def:HamiltonianStructure}. And let $\mathcal{M} (\Theta, A)$ be as in Section \ref{def:HamiltonianStructure}. Define $$Maslov ^{vert} (A) $$ to be the 
Maslov number of the triple $(\mathcal{V}, \Xi, \{\Xi _{i} \})$ determined by the pullback by $\sigma \in \mathcal{M} (A)$ of the 
vertical tangent bundle of $\widetilde{S}$, $\mathcal{L}$.
Then the expected dimension of $\mathcal{M} (A)$ is:
\begin{equation} \label{eq:dimensionMaslov}
r \cdot \chi (S) + Maslov ^{vert}  (A).
\end{equation}

\bibliographystyle{siam}     
\bibliography{C:/Users/yasha/texmf/bibtex/bib/link}

\def\cprime{$'$}
\begin{thebibliography}{10}

\bibitem{citeAkveldSalamonLoopsofLagrangiansubmanifoldsandpseudoholomorphicdiscs}
{\sc M.~Akveld and D.~Salamon}, {\em {Loops of Lagrangian submanifolds and
  pseudoholomorphic discs.}}, Geom. Funct. Anal., 11 (2001), pp.~609--650.

\bibitem{citeAlbersPSS}
{\sc P.~Albers}, {\em A {Lagrangian} {Piunikhin}-{Salamon}-{Schwarz} morphism
  and two comparison homomorphisms in {Floer} homology}, Int. Math. Res. Not.,
  2008 (2008), p.~56.
\newblock Id/No rnm134.

\bibitem{citeAtiyahBottTheYang-MillsequationsoverRiemannsurfaces}
{\sc M.~F. Atiyah and R.~Bott}, {\em {The Yang-Mills equations over Riemann
  surfaces.}}, Philos. Trans. R. Soc. Lond., A, 308 (1983), pp.~523--615.

\bibitem{citeBiranCorneaLagQuantHomology}
{\sc P.~{Biran} and O.~{Cornea}}, {\em {A Lagrangian quantum homology}},
  (2009), pp.~1--44.

\bibitem{citeBiranCorneaLagCobGAFA}
\leavevmode\vrule height 2pt depth -1.6pt width 23pt, {\em {Lagrangian
  cobordism and Fukaya categories}}, {Geom. Funct. Anal.}, 24 (2014),
  pp.~1731--1830.

\bibitem{citeCorneaBiranLagrangianCobordism}
{\sc O.~Cornea and P.~Biran}, {\em {Lagrangian Cobordism}}, arXiv:1109.4984.

\bibitem{citeGuilleminLermanEtAlSymplecticfibrationsandmultiplicitydiagrams}
{\sc V.~Guillemin, E.~Lerman, and S.~Sternberg}, {\em Symplectic fibrations and
  multiplicity diagrams}, Cambridge University Press, Cambridge, 1996.

\bibitem{citeReeseHarverSingular}
{\sc R.~{Harvey} and H.~B. jun. {Lawson}}, {\em {A theory of characteristic
  currents associated with a singular connection}}, {Bull. Am. Math. Soc., New
  Ser.}, 31 (1994), pp.~54--63.

\bibitem{citeHatcherAlgebraic}
{\sc A.~T. {Hatcher}}, {\em {Algebraic topology.}}, Cambridge: Cambridge
  University Press, 2002.

\bibitem{citeHuLalondeArelativeSeidelmorphismandtheAlbersmap}
{\sc S.~{Hu} and F.~{Lalonde}}, {\em {A relative Seidel morphism and the Albers
  map.}}, {Trans. Am. Math. Soc.}, 362 (2010), pp.~1135--1168.

\bibitem{citeLagrangianHofer}
{\sc H.~{Iriyeh} and T.~{Otofuji}}, {\em {Geodesics of Hofer's metric on the
  space of Lagrangian submanifolds.}}, {Manuscr. Math.}, 122 (2007),
  pp.~391--406.

\bibitem{citeFilteredHomotopyTypeFaria}
{\sc J.~F. Martins}, {\em On the homotopy type and the fundamental crossed
  complex of the skeletal filtration of a {CW}-complex}, Homology Homotopy
  Appl., 9 (2007), pp.~295--329.

\bibitem{citeMcDuffSalamonIntroductiontosymplectictopology}
{\sc D.~McDuff and D.~Salamon}, {\em Introduction to symplectic topology},
  Oxford Math. Monographs, The Clarendon Oxford University Press, New York,
  second~ed., 1998.

\bibitem{citeMcDuffSalamonJholomorphiccurvesandsymplectictopology}
\leavevmode\vrule height 2pt depth -1.6pt width 23pt, {\em $J$--holomorphic
  curves and symplectic topology}, no.~52 in American Math. Society Colloquium
  Publ., Amer. Math. Soc., 2004.

\bibitem{citeMilnorMorsetheory}
{\sc J.~Milnor}, {\em Morse theory}, Based on lecture notes by M. Spivak and R.
  Wells. Annals of Mathematics Studies, No. 51, Princeton University Press,
  Princeton, N.J., 1963.

\bibitem{citeYangMillsBlackHoles}
{\sc F.~Naeimipour, B.~Mirza, and F.~Jahromi}, {\em Yang–mills black holes in
  quasitopological gravity}, The European Physical Journal C, 81 (2021).

\bibitem{citeItamar}
{\sc I.~R. Rauch}, {\em {On The Hofer girth of the sphere of great circles}},
  \url{https://arxiv.org/pdf/2009.05256.pdf}.

\bibitem{citeRiehlAmodelstructureforquasi-categories}
{\sc E.~Riehl}, {\em A model structure for quasi-categories},
  \url{https://emilyriehl.github.io/files/topic.pdf}.

\bibitem{citeSavelyevGlobalFukayaCategoryI}
{\sc Y.~Savelyev}, {\em {Global Fukaya category I}}, Int. Math. Res. Not., (to
  appear), arXiv:1307.3991,
  \url{http://yashamon.github.io/web2/papers/fukayaI.pdf}.

\bibitem{citeSavelyevQuantumcharacteristicclassesandtheHofermetric}
\leavevmode\vrule height 2pt depth -1.6pt width 23pt, {\em {Quantum
  characteristic classes and the Hofer metric}}, Geom. Topol., 12 (2008),
  pp.~2277--2326.

\bibitem{citeSavelyevVirtualMorsetheory}
\leavevmode\vrule height 2pt depth -1.6pt width 23pt, {\em {Virtual Morse
  theory on $\Omega $Ham$(M, \omega)$.}}, J. Differ. Geom., 84 (2010),
  pp.~409--425.

\bibitem{citeSavelyevBottperiodicityandstablequantumclasses}
\leavevmode\vrule height 2pt depth -1.6pt width 23pt, {\em {Bott periodicity
  and stable quantum classes}}, Sel. Math., New Ser., 19 (2013), pp.~439--460.

\bibitem{citeSavelyevMorsetheoryfortheHoferlengthfunctional}
\leavevmode\vrule height 2pt depth -1.6pt width 23pt, {\em {Morse theory for
  the Hofer length functional}}, Journal of Topology and Analysis, 06, Issue
  No. 2 (2014).

\bibitem{citeSeidelFukayacategoriesandPicard-Lefschetztheory}
{\sc P.~Seidel}, {\em {Fukaya categories and Picard-Lefschetz theory}}, {Zurich
  Lectures in Advanced Mathematics. Z\"urich: European Mathematical Society
  (EMS). vi, 326~p. EUR~46.00 }, 2008.

\bibitem{citeSingularSibnerSibner}
{\sc L.~M. {Sibner} and R.~J. {Sibner}}, {\em {Classification of singular
  Sobolev connections by their holonomy}}, {Commun. Math. Phys.}, 144 (1992),
  pp.~337--350.

\bibitem{citeZingerPseudocycles}
{\sc A.~Zinger}, {\em Pseudocycles and integral homology}, Trans. Am. Math.
  Soc., 360 (2008), pp.~2741--2765.

\end{thebibliography}
\end{document}